\def\no{\if01}
\def\iftwelvept{\no}

\def\ifusepdf{\no}
\def\ifpsfont{\no}

\iftwelvept
\documentclass[leqno,12pt]{amsart}
\else
\documentclass[leqno,11pt]{amsart}
\fi
\usepackage{amssymb}
\usepackage{amscd}
\usepackage{latexsym}
\usepackage{verbatim}
\usepackage[all]{xy}

\ifusepdf
\usepackage{hyperref}
\else\fi
\ifpsfont
\usepackage[T1]{fontenc}
\usepackage{times}
\else\fi

\iftwelvept
\else\fi

\theoremstyle{plain}
\newtheorem{Theorem}{Theorem}[section]

\newtheorem{Proposition}[Theorem]{Proposition}
\newtheorem{Lemma}[Theorem]{Lemma}
\newtheorem{Corollary}[Theorem]{Corollary}

\theoremstyle{definition}

\newtheorem{Definition}[Theorem]{Definition}
\newtheorem{Remark}[Theorem]{Remark}
\newtheorem{Example}[Theorem]{Example}

\newcommand{\ZZ}{{\mathbb{Z}}}
\newcommand{\QQ}{{\mathbb{Q}}}

\newcommand{\CC}{{\mathbb{C}}}

\newcommand{\NN}{{\mathbb{N}}}
\newcommand{\NNNN}{\operatorname{N}}

\newcommand{\hhhh}{\mathfrak{h}}

\newcommand{\Sh}{\operatorname{Sh}}

\newcommand{\CCC}{{\mathcal{C}}}
\newcommand{\MSp}{\operatorname{MSp}_{\textup{\'et}}}
\newcommand{\ESh}{\operatorname{Sh}_{\textup{\'et}}}
\newcommand{\NSp}{\operatorname{MSp}}
\newcommand{\EM}{\textup{M}_{\textup{\'et}}}

\newcommand{\Et}{\textup{\'Et}}
\newcommand{\AMod}{\operatorname{AMod}}
\newcommand{\sss}{\mathfrak{S}}

\newcommand{\PR}{\operatorname{Pr}^{\textup{L}}}

\newcommand{\Rep}{\operatorname{Rep}}
\newcommand{\PRep}{\operatorname{PRep}}
\newcommand{\OO}{{\mathcal{O}}}

\newcommand{\LLL}{\mathcal{L}}

\newcommand{\DM}{\mathsf{DM}}

\newcommand{\Hom}{\operatorname{Hom}}
\newcommand{\Ext}{\operatorname{Ext}}

\newcommand{\Ker}{\operatorname{Ker}}
\newcommand{\Coker}{\operatorname{Coker}}
\newcommand{\Comp}{\operatorname{Comp}}
\newcommand{\Spec}{\operatorname{Spec}}

\newcommand{\Gal}{\operatorname{Gal}}

\newcommand{\Cor}{\operatorname{Cor}}
\newcommand{\SP}{\operatorname{Sp}}

\newcommand{\Mod}{\operatorname{Mod}}
\newcommand{\PMod}{\operatorname{PMod}}

\newcommand{\SSS}{\mathcal{S}}

\newcommand{\colim}{\operatorname{colim}}
\newcommand{\Cat}{\textup{Cat}_{\infty}}

\newcommand{\Map}{\operatorname{Map}}

\newcommand{\Fun}{\operatorname{Fun}}

\newcommand{\End}{\operatorname{End}}

\newcommand{\Aff}{\operatorname{Aff}}
\newcommand{\GL}{\textup{GL}}
\newcommand{\Grp}{\operatorname{Grp}}

\newcommand{\wCat}{\widehat{\textup{Cat}}_{\infty}}

\newcommand{\CAlg}{\operatorname{CAlg}}

\newcommand{\calD}{\mathcal{D}}
\newcommand{\QC}{\operatorname{QC}}

\newcommand{\Ind}{\operatorname{Ind}}

\newcommand{\Aut}{\operatorname{Aut}}

\newcommand{\XX}{\mathcal{X}}
\newcommand{\Sym}{\operatorname{Sym}}
\newcommand{\Proof}{{\sl Proof.}\quad}
\newcommand{\QED}{{\unskip\nobreak\hfil\penalty50\quad\null\nobreak\hfil
{$\Box$}\parfillskip0pt\finalhyphendemerits0\par\medskip}}

\begin{document}

\title{Mixed motives and Quotient stacks: Abelian varieties}

\author{Isamu Iwanari}



\address{Mathematical Institute, Tohoku University, Sendai, Miyagi, 980-8578,
Japan}

\email{iwanari@math.tohoku.ac.jp}

\maketitle

\section{Introduction}

In this paper, we investigate the structure of the motivic Galois group
of mixed motives tensor-generated by an abelian variety.
In our papers \cite{Tan}, \cite{Bar} and \cite{DTD},
we developed tannakian constructions for symmetric
monoidal stable $\infty$-categories
by means of derived algebraic geometry and proved Tannaka duality
results for symmetric monoidal stable $\infty$-categories.
We apply this theory to the symmetric monoidal stable $\infty$-category
of mixed motives and find consequences for the structures of
motivic Galois groups.

\vspace{2mm}

Let $\mathbf{DM}_{\textup{gm}}(k, \QQ)$ be the triangulated
category of mixed motives, which was constructed independently
by Hanamura, Levine, and Voevodsky. Here $k$ is a base perfect field
and $\QQ$ means rational coefficients.
We briefly recall the original approach to Galois groups of mixed motives.
It has been conjectured that
$\mathbf{DM}_{\textup{gm}}(k, \QQ)$
admits a $t$-structure which satisfies certain properties,
that is called a motivic $t$-structure (cf. \cite{bei}, \cite{Ha}).
If a motivic $t$-structure exists, then
the heart $\mathcal{MM}$ is a tannakian category
of mixed motives equipped with the realization functors of Weil
cohomology theories.
A motivic Galois group of mixed motives is defined to be
its Tannaka dual of $\mathcal{MM}$, i.e., the
pro-algebraic group over a suitable field corresponding to
$\mathcal{MM}$.
We refer the reader to e.g. \cite{A},
\cite{LHB}, \cite{SMG} for
topics on motivations and structures of motivic Galois groups.
In the case of mixed Tate motives,
there is the motivic $t$-structure on
the triangurated category (e.g., $k$ is a number field),
and
we have
the motivic Galois group of mixed Tate motives, which has various applications
and has received much attention.

Meanwhile,
the existence of a motivic $t$-structure of $\mathbf{DM}_{\textup{gm}}(k, \QQ)$
remains inaccessible (except for the mixed Tate case), and it is eventually related with Beilinson-Soul\'e
vanishing conjecture, Bloch-Beilinson-Murre filtration
and Grothendieck's standard conjectures.
In \cite{Tan}, \cite{Bar},
we developed a ``higher approach''.
Let us briefly summarize our previous results obtained in
\cite{Tan}, \cite{Bar}.
Let $\DM^\otimes$ be the symmetric monoidal
stable presentable
$\infty$-category of mixed motives (see Section 3 or \cite{Tan}, \cite{Bar}, \cite{DTD}). By an $\infty$-category, we means an $(\infty,1)$-category
(cf. \cite{Ber}) in this Introduction, but we use the theory of
quasi-categories from the next section (cf. \cite{Jo}, \cite{HTT}).
Roughly speaking, $\DM^\otimes$ is an $\infty$-categorical enhancement
of $\mathbf{DM}_{\textup{gm}}(k, \QQ)$.
That is, if $\DM^\otimes_{\textup{gm}}$ denotes the stable subcategory
of $\DM^\otimes$ spanned by compact objects,
the homotopy category of $\DM^\otimes_{\textup{gm}}$
is equivalent to $\mathbf{DM}_{\textup{gm}}(k, \QQ)$
as symmetric monoidal triangulated categories.
Let $\mathbf{K}$ be a field of characteristic zero
and $\mathcal{D}^\otimes(\mathbf{K})$ the symmetric monoidal unbounded derived $\infty$-category
of $\mathbf{K}$-vector spaces.
For a Weil cohomology theory
with $\mathbf{K}$-coefficients,
we have the (homological)
realization functor $\mathsf{R}:\mathsf{DM}^\otimes\to \calD(\mathbf{K})^\otimes$,
which is a symmetric monoidal
colimit-preserving functor.
For example, the realization functor
of de Rham cohomology carries the motif $M(X)$ of a smooth projective variety $X$
to the dual of chain complex computing de Rham cohomology of $X$
(there are also \'etale, singular (Betti),
de Rham realizations, etc.; see Section 4 and \cite[Section 5]{Tan} based on mixed Weil theories).
On the basis of the new framework of $\infty$-categories (e.g., \cite{HTT}) and derived algebraic geometry \cite{DAGn} \cite{HAG2},
we have constructed a derived
group scheme $\mathsf{MG}_{\mathsf{R}}$.
The notion of derived affine group schemes is a direct
generalization of affine group
schemes in derived algebraic geometry (see \cite[Appendix]{Tan}).
The derived affine group scheme $\mathsf{MG}_{\mathsf{R}}$ over $\mathbf{K}$
has the following properties (see \cite{Tan}, \cite{Bar} for the details):

\begin{itemize}

\item $\mathsf{MG}_{\mathsf{R}}$ represents the automorphism group
of the realization functor $\mathsf{R}$. We can obtain an ordinary pro-algebraic group $MG_{\mathsf{R}}$
from $\mathsf{MG}_{\mathsf{R}}$ by a truncation procedure, and
$MG_{\mathsf{R}}$ coarsely represents the automorphism group
of $\mathsf{R}$
(the representability is one important reason for the
need to work in the framework of $\infty$-categories).
Namely, $MG_{\mathsf{R}}$ is the coarse moduli space
for $\mathsf{MG}_{\mathsf{R}}$.
We shall refer to $\mathsf{MG}_{\mathsf{R}}$
and $MG_{\mathsf{R}}$ as the derived motivic Galois group of mixed motives
(with respect to $\mathsf{R}$) and the motivic Galois group of mixed motives
respectively. We remark that one can also construct the affine derived scheme representing the comparison torsor between singular and de Rham realization
functors.

\item If a motivic $t$-structure exists, then the Tannaka
dual of the heart $\mathcal{MM}$
with $\mathbf{K}$-coefficients is isomorphic to $MG_{\mathsf{R}}$.
That is, the symmetric monoidal
category of finite dimensional representations of
$MG_{\mathsf{R}}$ is equivalent to $\mathcal{MM}$ (after the base change to
$\mathbf{K}$).

\item If the same construction is applied to the stable $\infty$-category
of mixed Tate motives $\mathsf{DTM}^\otimes\subset \mathsf{DM}^\otimes$,
then the associated pro-algebraic group $MTG$ coincides with the motivic Galois
group of mixed Tate motives constructed by Bloch-Kriz, and others.

\end{itemize}

The guiding principle behind this construction is that
$\mathsf{DM}^\otimes$ equipped with
a realization functor should form a ``tannakian $\infty$-category'',
and the coarse moduli space of the Tannaka dual of $\mathsf{DM}^\otimes$
should be the Tannaka dual of the heart of a motivic $t$-structure;
especially
we focus on automorphism groups of fiber functors.
The following table shows the principle of correspondences:

\vspace{2mm}

\begin{center}
\begin{tabular}{c  l}
\hline
Category & Group \\
\hline
$\mathsf{DM}^\otimes$ & $\mathsf{MG}_\mathsf{R}$ \\
A conjectural heart & $MG_\mathsf{R}$ \\
\hline
\end{tabular}
\end{center}

\vspace{2mm}

We have one more important category to introduce, and
that is the category $NM$ of Grothendieck's numerical motives
 (cf. \cite{A}, Section 3).
The category $NM$ was introduced by Grothendieck with the aim of
proving Weil conjectures. By a theorem of
Jannsen, it is
a semi-simple abelian category.
According to the perspective of mixed motives due to
Beilinson and Deligne, it has been conjectured
that there is a deep relation between $NM$ and $\mathsf{DM}$
(cf. \cite{Be}, \cite{Aquo}, \cite{Nek}): The abelian (and tannakian) subcategory of the conjectural tannakian category $\mathcal{MM}$ spanned by semi-simple objects should be equivalent to $NM$, and moreover
there is a (weight) filtration of every object in $\mathcal{MM}$ whose
graded quotients are semi-simple objects. (We remark that
$NM$ is conjecturally tannakian.)
We now turn our attention to motivic Galois groups.
In the light of the above conjectural relation,
it is conjectured that
the structure of motivic Galois group for mixed motives
is given by
\[
MG \simeq UG \rtimes  MG_{pure}
\]
such that $MG_{pure}$ is a
pro-reductive group scheme which is a Tannaka dual of the (conjecturally tannakian) category of Grothendieck's numerical motives, and $UG$ is a pro-unipotent group scheme
which encodes the data of extensions in $\mathcal{MM}$
and should be described in terms of motivic complexes; see \cite[1.3.3]{Aquo}, \cite[5.3.1]{LHB}.
This decomposition is known (only) for the case of mixed Tate motives:
If $MTG$ denotes the motivic Galois group for mixed Tate motives,
there is a decomposition $MTG \simeq UG_{Tate}\rtimes \mathbb{G}_m$,
such that $\mathbb{G}_m$ is one dimensional algebraic torus that is the Tannaka
dual of numerical Tate motives, and $UG_{Tate}$ is a pro-unipotent group.

\vspace{2mm}

Let $\DM_X^\otimes\subset \DM^\otimes$
be the symmetric monoidal stable presentable
subcategory generated (as a symmetric monoidal stable presentable
$\infty$-category) by the motives $M(X)$
and the dual $M(X)^\vee$ of an abelian variety $X$ over a field of characteristic zero
(see Definition~\ref{memdef}).
Let $NM_X$ be the tannakian category which is the symmetric monoidal
abelian subcategory of $NM$ generated (as a symmetric monoidal abelian
category) by the numerical motive of $X$ and its dual (in this case, Grothendieck's standard conjectures hold, and thus $NM_X$ is tannakian; see \cite{A}).
We denote by
 $MG_{pure}(X)$ its Tannaka dual, i.e., the motivic Galois group of
numerical (pure) motives $NM_X$.
Let $MG_{\textup{\'et}}(X)$ be the motivic Galois group
of $\DM_X^\otimes$ with respect to an \'etale realization functor $\mathsf{R}_{\textup{\'et}}:\DM_X^\otimes\to \mathcal{D}^\otimes(\QQ_l)$,
that is constructed from $\DM_X^\otimes$
in the same way as $MG_{\mathsf{R}}$ from $\DM^\otimes$.
It is important to note that an existence of a motivic $t$-structure
on $\DM_X^\otimes$ is still inaccessible.
The following our main result proves the conjectural structure of motivic
Galois group of $\DM_X^\otimes$
(for many cases of abelian varieties). See Theorem~\ref{unconditional}, Corollary~\ref{unconditional2} for details.

\begin{Theorem}
\label{main4th}
Let $X$ be an abelian variety over a number field
$k$. Let $MG_{pure}(X)_{\QQ_l}$ be the reductive algebraic group
over $\QQ_l$ which is the Tannaka dual of the ($\QQ_l$-linear)
tannakian category of numerical motives generated by
the motives of $X$ (see Section 3 and 4).
Suppose that $X$ satisfies either (i) $\End(X\otimes_k\bar{k})=\ZZ$ or (ii) $X$ is one dimensional, or (iii) $X$ is a simple CM abelian variety of prime dimension
 (with some more conditions, see 
 Theorem~\ref{unconditional} for the precise
formulation).
Then there exists an exact sequence of pro-algebraic group schemes
\[
1\to UG_{\textup{\'et}}(X)\to MG_{\textup{\'et}}(X) \to MG_{pure}(X)_{\QQ_l}\to 1.
\]
Moreover, $UG_{\textup{\'et}}(X)$ is a connected pro-unipotent group scheme over $\QQ_l$ which is constructed from a motivic algebra $A_X$, which is
a commutative differential graded algebra.
We refer the reader to Section~\ref{decomposition} for further details.
\end{Theorem}

\begin{Corollary}
\label{main5th}
\begin{enumerate}
\renewcommand{\labelenumi}{(\roman{enumi})}
\item There is an isomorphism of affine group schemes
\[
MG_{\textup{\'et}}(X) \simeq UG_{\textup{\'et}}(X)\rtimes MG_{pure}(X)_{\QQ_l}.
\]

\item $UG_{\textup{\'et}}(X)$ is the unipotent radical of $MG_{\textup{\'et}}(X)$, i.e., the maximal normal unipotent closed subgroup.

\item $MG_{pure}(X)_{\QQ_l}$ is the reductive quotient.
\end{enumerate}
\end{Corollary}

Corollary~\ref{main5th} especially
says that one can extract the Tannaka dual of $NM_X$
from $\mathsf{DM}_X^\otimes$ in the group-theoretic
fashion: (i) take a derived affine group scheme $\mathsf{MG}_{\textup{\'et}}(X)$
which represents
$\Aut(\mathsf{R}_{\textup{\'et}})$, (ii) then take its coarse moduli space
$MG_{\textup{\'et}}(X)$, (iii) finally, the reductive quotient of
$MG_{\textup{\'et}}(X)$, i.e., the quotient by the unipotent radical
is the Tannaka dual of numerical motives generated by $X$.

\vspace{1mm}

Our proof of Theorem~\ref{main4th}
is based on the interplay between
 two main ingredients:

\begin{itemize}

\item We apply results from derived Tannaka duality and techniques
from derived algebraic geometry: $\DM_X^\otimes$ is
a fine $\infty$-category in the sense of \cite{DTD}, and there exist
a derived quotient stack $Z=[\Spec A_X/\GL_n]$ and an equivalence
$\DM_X^\otimes\simeq \QC^\otimes(Z)$ where $\QC^\otimes(Z)$
denotes the symmetric monoidal stable $\infty$-category of quasi-coherent
complexes. Then the derived motivic Galois group is obtained from
$Z$ by the construction of the based loop space.

\item We use results on images of Galois representations of abelian varieties (cf. Mumford-Tate conjecture on abelian varieties, see e.g., Introduction of
\cite{Pi}).
Arguably, the simplest but (highly)
nontrivial result in this direction is Serre's theorem,
which states that
the image of $l$-adic Galois representation $\Gal(\bar{k}/k)\to \GL_2(\QQ_l)$
of an elliptic curve without complex multiplication over a number field $k$
is dense. To apply them, we also construct an \'etale realization functor endowed with a Galois action.

\end{itemize}

\vspace{1mm}

Let us give some instructions to the reader.
In Section 2 we recall some generalities 
concerning $\infty$-categories, $\infty$-operads,
and spectra, etc.
In Section 3, applying the Tannakian characterization in \cite{DTD}
to $\DM_X^\otimes$
we discuss the consequences.
We also include some basis definitions
and facts about Chow and numerical motives, and mixed motives.
In Section 4 we construct and study the motivic Galois group for $\mathsf{DM}_X^\otimes$.
That is to say, the main results in Section 4 are Theorem~\ref{main4th}
and Corollary~\ref{main5th} in this introduction.
In the final Section, we construct an $l$-adic realization functor
from the symmetric monoidal $\infty$-category 
of mixed motives with $\ZZ$-coefficients to the derived
$\infty$-category of $\ZZ_l$-modules, which is endowed with
Galois action; see Proposition~\ref{realization}, Remark~\ref{realrem}. This is used in Section 4. But
the construction of a realization functor has a
different nature from the main objectives of this paper.
Thus we treat this issue in the final Section.
The author would like to thank S. Yasuda and S. Mochizuki for
enlightening questions and valuable comments on the case of abelian schemes.
The author is partially supported by Grant-in-Aid for Scientific Research,
Japan Society for the Promotion of Science.

\section{Notation and Convention}

We fix some notation and convention.

\subsection{$\infty$-categories}
\label{NC1}
In this paper we use the theory of {\it quasi-categories}.
A quasi-category is a simplicial set which
satisfies the {\it weak Kan condition} of Boardman-Vogt.
The theory of quasi-categories from higher categorical viewpoint
has been extensively developed by Joyal and Lurie.
Following \cite{HTT} we shall refer to quasi-categories
as {\it $\infty$-categories}.
Our main references are \cite{HTT}
 and \cite{HA}
(see also \cite{Jo}, \cite{DAGn}).
We often refer to a map $S\to T$ of $\infty$-categories
as a {\it functor}. We call a vertex in an $\infty$-category $S$
(resp. an edge) an {\it object} (resp. a {\it morphism}).
When $S$ is an $\infty$-category, by $s\in S$ we
mean that $s$ is an object of $S$.
For the rapid introduction to $\infty$-categories, we
refer to \cite[Chapter 1]{HTT}, \cite{Gro}, \cite[Section 2]{FI}.
We remark also that
there are several alternative theories such as Segal categories,
complete Segal spaces, simplicial categories, relative categories,
etc. For reader's convenience,
we list some of notation concerning $\infty$-categories:

\begin{itemize}

\item $\Delta$: the category of linearly ordered finite sets (consisting of $[0], [1], \ldots, [n]=\{0,\ldots,n\}, \ldots$)

\item $\Delta^n$: the standard $n$-simplex

\item $\textup{N}$: the simplicial nerve functor (cf. \cite[1.1.5]{HTT}). But
we
do not often distinguish notationally between ordinary categories and their nerves.

\item $\mathcal{C}^{op}$: the opposite $\infty$-category of an $\infty$-category $\mathcal{C}$

\item Let $\mathcal{C}$ be an $\infty$-category and suppose that
we are given an object $c$. Then $\mathcal{C}_{c/}$ and $\mathcal{C}_{/c}$
denote the undercategory and overcategory respectively (cf. \cite[1.2.9]{HTT}).

\item $\operatorname{Cat}_\infty$: the $\infty$-category of small $\infty$-categories in a fixed Grothendieck universe $\mathbb{U}$
(cf. \cite[3.0.0.1]{HTT}). We employ
the ZFC-axiom together with the universe axiom
of Grothendieck.
We have a sequence of universes
$(\NN\in)\mathbb{U}\in \mathbb{V}\in \mathbb{W}\in \ldots$.
If $x$ belongs to $\mathbb{V}$ (resp. $\mathbb{W}$),
we call $x$ large (resp. super-large).

\item $\wCat$: $\infty$-category of large $\infty$-categories.

\item $\SSS$: $\infty$-category of small spaces (cf. \cite[1.2.16]{HTT})

\item $\textup{h}(\mathcal{C})$: homotopy category of an $\infty$-category (cf. \cite[1.2.3.1]{HTT})

\item $\Fun(A,B)$: the function complex for simplicial sets $A$ and $B$


\item $\Map(A,B)$: the largest Kan complex of $\Fun(A,B)$ when $A$ and $B$ are $\infty$-categories,


\item $\Map_{\mathcal{C}}(C,C')$: the mapping space from an object $C\in\mathcal{C}$ to $C'\in \mathcal{C}$ where $\mathcal{C}$ is an $\infty$-category.
We usually view it as an object in $\mathcal{S}$ (cf. \cite[1.2.2]{HTT}).

\item $\Ind(\mathcal{C})$; $\infty$-category of Ind-objects
in an $\infty$-category $\mathcal{C}$
(cf. \cite[5.3.5.1]{HTT}, \cite[6.3.1.13]{HA}).

\end{itemize}

\subsection{Stable $\infty$-categories}

We shall employ the theory of stable $\infty$-categories developed
in \cite{HA}.
The homotopy category of a stable $\infty$-category
forms a triangulated category in a natural way.
For generalities
we refer to \cite[Chap. 1]{HA}.
We denote the suspension functor
and the loop functor by $\Sigma$ and $\Omega$ respectively.
For a stable $\infty$-category $\mathcal{C}$
and two objects $C,C'\in \mathcal{C}$,
we write $\Ext^i_{\mathcal{C}}(C,C')$
for $\pi_0(\Map_{\mathcal{C}}(C,\Sigma^i (C')))$.
If no confusion seems to arise, we also use the shift $[-]$ instead of $\Sigma$ and $\Omega$ when we treat (co)chain complexes.

\subsection{Symmetric monoidal $\infty$-categories}
We use the theory of symmetric monoidal $\infty$-categories
developed in \cite{HA}.
We refer to \cite{HA} for its generalities.

\begin{itemize}

\item $\Mod_A$: the stable $\infty$-category of $A$-module spectra
for a commutative ring spectrum $A$. We usually write $\Mod_A^\otimes$
for $\Mod_A$ equipped with
the symmetric monoidal structure given by smash product $(-)\otimes_A(-)$.
See \cite[4.4]{HA}.
When $R$ is the Eilenberg-Maclane spectrum $HK$
of an ordinary algebra $A$, then we write $\Mod_{K}^\otimes$
for $\Mod_{HA}^\otimes$.
We denote by $\PMod_A^\otimes$ the symmetric monoidal full sucategory of $\Mod_A^\otimes$ that consists of dualizable objects.
 For the definition of dualizable objects,
see e.g. \cite{Tan}, \cite{Bar}, \cite{DTD}. 
In the literature dualizable objects
are also called strongly dualizable objects.

\item $\CAlg(\mathcal{M}^\otimes)$: $\infty$-category of commutative
algebra objects in a symmetric
monoidal $\infty$-category $\mathcal{M}^\otimes$.
See \cite[2.13]{HA}.

\item $\CAlg_R$: $\infty$-category of commutative
algebra objects in the symmetric monoidal $\infty$-category $\Mod_R^\otimes$
where $R$ is a commutative ring spectrum.
We write $\CAlg$ for the $\infty$-category of commutative
algebra objects in $\Mod_\mathbb{S}^\otimes$ (i.e., $E_\infty$-ring spectra)
where $\mathbb{S}$ is the sphere spectrum.
If $A$ is an ordinary commutative ring,
then we denote by $HA$ the Eilenberg-MacLane spectrum that belongs to
$\CAlg$. For simplicity, we write $\CAlg_A$ for $\CAlg_{HA}$. We remark that the full subcategory $\CAlg_{HA}^{\textup{dis}}$
of $\CAlg_{HA}$ spanned by
discrete objects $M$ (i.e., $\pi_n(M)=0$ when $n\neq 0$)
is naturally categorical equivalent to the nerve of the category of usual commutative
$A$-algebras.
The inclusion is given by the Eilenberg-MacLane functor $A\mapsto HA$.
We abuse notation and often write $A$ for $HA$.
When $A$ is an ordinary commutative ring that contains the field $\QQ$,
$\CAlg_{A}$ is equivalent to the $\infty$-category obtained from
the category of commutative differential graded $A$-algebras
by inverting quasi-isomorphisms (cf. \cite[8.1.4.11]{HA}).
Therefore, unless stated otherwise we refer to an object
in $\CAlg_{A}$ (i.e., an $E_\infty$-ring spectrum over $A$) as a commutative differentail graded $A$-algebra.

\item $\Mod_A^\otimes(\mathcal{M}^\otimes)$: symmetric monoidal
$\infty$-category of
$A$-module objects,
where $\mathcal{M}^\otimes$
is a symmetric monoidal $\infty$-category such that (1)
the underlying $\infty$-category admits a colimit for any simplicial diagram, and (2)
its tensor product functor $\mathcal{M}\times\mathcal{M}\to \mathcal{M}$
preserves
colimits of simplicial diagrams separately in each variable.
Here $A$ belongs to $\CAlg(\mathcal{M}^\otimes)$
(cf. \cite[3.3.3, 4.4.2]{HA}).

\end{itemize}

 Let $\mathcal{C}^\otimes$ be a symmetric monoidal $\infty$-category.
We usually denote, dropping the superscript $\otimes$,
by $\mathcal{C}$ its underlying $\infty$-category.
If no confusion likely arises, we omit the superscript $(-)^\otimes$.

Let $(\PR)^\otimes$ be the symmetric monoidal $\infty$-category
of presentable $\infty$-categories, see \cite[6.3.1.4, 6.3.1.6]{HA}
(see also \cite[Section 2]{DTD}). In $\PR$,
a morphism is a colimit-preserving functor.
A symmetric monoidal presentable $\infty$-category
$\CCC^\otimes$ whose tensor product $\CCC\times \CCC\to \CCC$
preserves (small) colimits separately in each variable can be
naturally viewed as a commutative algebra object in $(\PR)^\otimes$.
Thus $\Mod_{R}^\otimes$ belongs to $\CAlg((\PR)^\otimes)$.
We refer to $\CAlg((\PR)^\otimes)_{\Mod_{R}^\otimes/}$
as the $\infty$-category of $R$-linear symmetric monoidal
(stable) presentable $\infty$-categories. When $R=HA$,
we use the word ``$A$-linear'' instead of $HA$-linear.

\subsection{Model categories and $\infty$-categories}
\label{modelinf}
Our references of model categories are \cite{Ho} and \cite[Appendix]{HTT}.
Let $\mathbb{M}$ be a combinatorial model category.
We can obtain a presentable $\infty$-category from $\mathbb{M}$ by
inverting weak equivalences. 
For example, we can apply
Dywer-Kan hammock localization and take the simplicial nerve (after
a fibrant replacement). Similarly, we can obtain
a symmetric monoidal presentable $\infty$-category
from a combinatorial symmetric monoidal model category.
We refer the reader to \cite[1.3.4, 4.1.3]{HA}, (or \cite[Section 2]{DTD})
for details.

\subsection{Derived stacks and quasi-coherent complexes}
We use derived stacks and (derived) quotient stacks
and quasi-coherent complexes on them.
For our convention, we refer the reader to \cite[Section 4]{Bar}
or \cite[Section 2.1, 2.3]{DTD}.

\section{Fine tannakian $\infty$-category of mixed abelian motives}

The purpose of this Section
is to deduce a tannakian type representation of
symmetric monoidal stable $\infty$-category of mixed motives generated
by the motif of an abelian scheme (or more generally, Kimura finite motif)
from results of \cite{DTD}.

\subsection{}
\label{Modelmot}

We review briefly the symmetric monoidal stable
$\infty$-category of mixed motives.
We use a model category-theoretic device
and define a $\QQ$-linear symmetric monoidal stable presentable
$\infty$-category $\mathsf{DM}^\otimes$.
Let $S$ be a smooth quasi-projective scheme over a perfect field $k$.
We work with rational coefficients,
though this assumption is not needed
to the construction of the stable $\infty$-category of mixed motives.
We first recall a $\QQ$-linear
category $\operatorname{Cor}$, see \cite[Lec. 1]{MVW}, \cite{CD3}. Objects
in $\operatorname{Cor}$
are smooth separated $S$-schemes of finite type
which we regard as formal symbols.
We denote by $\LLL(X)$ the object in $\textup{Cor}$ corresponding
to a smooth scheme $X$.
For $X$
and $Y$, we let $\Hom_{\textup{Cor}}(\LLL(X),\LLL(Y))$ be the $\QQ$-vector
space $c_0(X\times_SY/X)$ of finite $S$-correspondences
(see \cite[9.1.2]{CD3}).
The composition is determined by
intersection product
(see \cite[page 4]{MVW}, \cite[9.1]{CD3}).
Let $\textup{Sm}_{/S}$ be the
category of (not necessarily connected)
smooth separated $S$-schemes of finite type.
Then there is a functor $\LLL:\textup{Sm}_{/S}\to \Cor$
which carries $X$ to $\LLL(X)$ and sends $f:X\to Y$ to the graph
$\Gamma_f\in \Hom_{\operatorname{Cor}}(\LLL(X),\LLL(Y))$.
Let $\LLL(X)\otimes \LLL(Y)=\LLL(X\times_SY)$ and
define $\gamma_{X,Y}:\LLL(X)\otimes \LLL(Y)\to \LLL(Y)\otimes \LLL(X)$ to be
the isomorphism induced by the graph of the flip
$X\times_S Y\to Y\times_S X$.
These data makes $\Cor$ a symmetric monoidal category,
which shall call as the $\QQ$-linear category of finite $S$-correspondences.
Next let $\widehat{\operatorname{Cor}}$
be the nerve of
functor category consisting of $\QQ$-linear functors from $\Cor^{op}$ to $\textup{Vect}_{\QQ}$;
$\Fun_{\QQ}(\operatorname{Cor}^{op},\operatorname{Vect}_\QQ)$.
Here $\textup{Vect}_\QQ$ denotes the category of $\QQ$-vector spaces.
By enriched Yoneda's lemma, $\operatorname{Cor}$ can be viewed as the full subcategory
of $\widehat{\operatorname{Cor}}$, and every object
of $\operatorname{Cor}$ is compact in $\widehat{\operatorname{Cor}}$. Since $\widehat{\operatorname{Cor}}$ admits small
colimits, by \cite[5.3.5.11]{HTT} we have
a fully faithful left Kan extension $\Ind(\operatorname{Cor})\to \widehat{\operatorname{Cor}}$.
Since $\Ind(\operatorname{Cor})$ admits filtered colimits,
it is idempotent complete.
Day convolution product defines
the symmetric monoidal structure on $\Ind(\textup{Cor})$
whose tensor product preserves filtered colimits separately in each
variable,
and $\LLL(X)\otimes \LLL(Y)=\LLL(X\times_S Y)$ (cf. \cite[6.3.1.13]{HA}).
We say that a $\QQ$-linear functor $F:\Cor^{op}\to \textup{Vect}_\QQ$ if a Nisnevich
sheaf
is the composite $\textup{Sm}_{/S}^{op}\stackrel{\LLL}{\to}\Cor^{op}\stackrel{F}{\to} \textup{Vect}_\QQ$ is a sheaf with respect to Nisnevich topology
on $\textup{Sm}_{/S}$ (see e.g. \cite{MVW} for the definition).
For example, $\LLL(X)$ and its direct summands are sheaves if $X\in \textup{Sm}_{/S}$.
Let $\Sh\subset \widehat{\Cor}$ be the full subcategory
of sheaves.
It is a Grothendieck abelian category (cf. \cite[2.4]{CD1}).

Let us recall from  \cite[2.4]{CD1} the descent structure of $\Sh$.
Let $\mathcal{G}_{\Sh}$ be the set of $\{\LLL(X)\}_{X\in \textup{Sm}_{/S}}$.
Let $\mathcal{H}_{\Sh}$ be the set of complexes obtained as
cones of $\LLL(\mathcal{X})\to \LLL(X)$ where $X\in \textup{Sm}_{/S}$
and $\XX\to X$ is any Nisnevich hypercovering of $X$.
The pair $(\mathcal{G}_{\Sh},\mathcal{H}_{\Sh})$
is a weakly flat descent structure in the sense of \cite{CD1}
(see \cite[3.3]{CD1}).
Then $\textup{Comp}(\Sh)$ admits a symmetric monoidal
model structure, described in \cite[2.5, 3.2]{CD1}, in which weak equivalences are quasi-isomorphisms,
and cofibrations are $\mathcal{G}_{\textup{Sh}}$-cofibrations.
We call this model structure the $\mathcal{G}_{\Sh}$-model structure.
Next put $\mathcal{T}$ be the set of complexes of sheaves
with transfers obtained as cones of $p_*:\LLL(X\times_S\mathbb{A}_S^1)\to \LLL(X)$
for any $X\in \textup{Sm}_{/S}$, where $\mathbb{A}_S^1$
is the affine line over $S$, and $p:X\times_S\mathbb{A}_S^1\to X$
is the natural projection.
Invoking \cite[4.3, 4.12]{CD1} we take the left Bousfield localization of the above model structure
on
$\textup{Comp}(\Sh)$ with respect to $\mathcal{T}$, in which
weak equivalences are called $\mathbb{A}^1$-local equivalences,
and fibrations are called $\mathbb{A}^1$-local fibrations.
Let $\QQ(1)$ be $\Ker(\mathcal{L}(\mathbb{G}_{m,S})\to \LLL(S))[-1]$
where $\mathbb{G}_{m,S}=\Spec \OO_S[t,t^{-1}]$.
Consider the symmetric monoidal category $\SP_{\QQ(1)}^\sss(\textup{Comp}(\Sh))$
of symmetric $\QQ(1)$-spectra (we shall refer the reader to \cite{HoG} for the generalities of symmetric spectra).
By \cite[7.9]{CD1} (see also \cite{HoG}),
$\SP_{\QQ(1)}^\sss(\textup{Comp}(\Sh))$ has a symmetric monoidal
model structure such that weak equivalences (resp. fibrations)
are termwise
$\mathbb{A}^1$-local equivalences (resp. termwise $\mathbb{A}^1$-local fibrations).
We refer to this model structure as the $\mathbb{A}^{1}$-local projective
model structure.
Following \cite[7.13]{CD1} and \cite{HoG}
we define the stable model structure on $\SP_{\QQ(1)}^\sss(\textup{Comp}(\Sh))$
to be the left Bousfield localization with respect to
$\{s_n^{\LLL(X)}:F_{n+1}(\LLL(X)\otimes \QQ(1))\to F_{n}(\LLL(X))\}_{X\in \textup{Sm}_{/S}}$ of the $\mathbb{A}^1$-local projective model structure
(see \cite[7.7]{HoG} for the notation $s_n^{\LLL(X)},\ F_n$, etc.).
We refer to a weak equivalence (resp. fibration)
in the stable model structure
as a stable equivalence (resp. stable fibration).
We let $\mathsf{Sp}_{\QQ(1)}^\sss(\textup{Comp}(\Sh))^\otimes$
be the symmetric monoidal $\infty$-category obtained from
(the full subcategory of cofibrant objects of) 
$\SP_{\QQ(1)}^\sss(\textup{Comp}(\Sh))$
by inverting stable equivalences (See Section~\ref{modelinf}, \cite[1.3.4, 4.1.3]{HA} or \cite[Section 2]{DTD}).
Set
$\mathsf{DM}^\otimes:=\mathsf{Sp}_{\QQ(1)}^\sss(\textup{Comp}(\Sh))^\otimes$
(we use the notation $\mathsf{DM}^\otimes(k)$ in \cite{Tan}).
We abuse notation and often write $\QQ(1)$ also
for the image of $\QQ(1)$
in $\mathsf{DM}^\otimes$ (and its homotopy category).

Let $\mathsf{DM}^{\textup{eff},\otimes}$ be the symmetric monoidal
$\infty$-category
obtained from the model category $\Comp(\textup{Sh})$ by inverting $\mathbb{A}^1$-local
equivalences. 
There is a natural symmetric monoidal
functor $\Sigma^{\infty}:\mathsf{DM}^{\textup{eff},\otimes}\to \mathsf{DM}^\otimes$ obtained from the left Quillen adjoint functor
$\textup{Comp}(\Sh)\to \SP_{\QQ(1)}^\sss(\textup{Comp}(\Sh))$.
For a smooth $S$-scheme $X$, we denote by $M(X)$
the image of $\mathcal{L}(X)$ and refer to it as the motif of $X$.
We denote by $M:\textup{Sm}_{/S} \to \DM^{\otimes}$
the (symmetric monoidal) functor that sends $X$ to $M(X)$.

Let $\DM_\vee^\otimes$ be a symmetric monoidal full subcategory of
$\DM^\otimes$
spanned by
dualizable objects.
When the base scheme $S$ is $\Spec k$, compact objects and dualizable objects
coincides in $\DM^\otimes$ (since we work with rational coefficients and
alteration).

\subsection{}

\begin{Definition}
\label{memdef}
Let $S$ be a smooth quasi-projective scheme over a perfect field $k$.
Let $X$ be an abelian scheme of relative dimension $g$ over $S$.
Let $\mathsf{DM}_X$ be the smallest
presentable stable subcategory of $\mathsf{DM}$,
which contains $M(X)$ and its dual $M(X)^\vee$ and are closed under
tensor products.
We refer to $\mathsf{DM}_X$ as the $\infty$-category of {\it mixed abelian motives}
generated by $X$.
Let $\mathsf{DM}_{\textup{gm},X}$ be the smallest
stable subcategory of $\mathsf{DM}$,
which contains $M(X)$ and $M(X)^\vee$ and are closed under
retracts and tensor products.
The stable $\infty$-categories $\mathsf{DM}_X$ and $\mathsf{DM}_{\textup{gm},X}$ inherit symmetric monoidal structures respectively in the obvious way.
\end{Definition}

\begin{Theorem}
\label{representationtheorem}
There exist a derived quotient stack $[\Spec A_X/\GL_{2g}]$
given by an action of $\GL_{2g}$ on a commutative differentail graded
$\QQ$-algebra $A_X$,
and 
an $\QQ$-linear symmetric monoidal equivalence
\[
\QC^\otimes([\Spec A_X/\GL_{2g}])\simeq \mathsf{DM}_X^\otimes.
\]
Moreover, the composite $\QC^\otimes(B\GL_{2g})\to \QC^\otimes([\Spec A_X/\GL_{2g}])\simeq \mathsf{DM}_X^\otimes$ carries the standard representation
of $\GL_{2g}$ to $M_1(X)[-1]$, where the first functor is the pullback functor
along the natural morphism $[\Spec A_X/\GL_{2g}]\to B\GL_{2g}$.
Here $\GL_{2g}$ denotes the general linear group scheme over $\QQ$.
Here $M_1(X)[-1]$ denotes a certain shifted
direct summand of $M(X)$, see Section~\ref{HCNM}. 
\end{Theorem}

We shall refer the underlying commutative differential graded
algebra $A_X$ as the motivic algebra for $X$.

\begin{Remark}
The symmetric monoidal stable $\infty$-category
$\QC^\otimes([\Spec A_X/\GL_{2g}])$ is equivalent to
$\Mod_{A_X}^\otimes(\QC^\otimes(B\GL_{2g}))$.
\end{Remark}

\begin{Remark}
The symmetric monoidal $\infty$-category $\QC^\otimes(B\GL_{2g})$
is equivalent to the symmetric monoidal $\infty$-category
obtained from the category of chain complexes of $\QQ$-linear
representations of $\GL_{2g}$
by inverting quasi-isomorphisms. See Section~\ref{verytech} for
details.
We remark that
in this paper, by a representation of an affine group scheme $G$ over
a field $k$
we means a $k$-vector space $V$ equipped with a $G$-action
(in the scheme-theoretic sense).
\end{Remark}

\subsection{Chow motives}
\label{HCNM}
Let us review
a symmetric monoidal $\QQ$-linear (additive) category $CHM^\otimes$
of homological (relative) Chow motives.
Our reference is \cite{sch} \cite{DM}, whereas in {\it loc. cit.} the cohomological theory is presented. But we shall adopt
the homological theory.

Let $S$ be a smooth quasi-projective variety over a perfect field $k$.
Let $\textup{SmPr}_{/S}$
denote the category of schemes that are smooth and projective over $S$.
Let $CHM'$ be the $\mathbb{Q}$-linear
category whose objects
are formal symbols $(X,i)$
where $X$ belongs to $\textup{SmPr}_{/S}$, and $i$ is an integer.
For two objects
$(X,i)$ and $(Y,j)$,
we define $\Hom_{CHM'}((X,i),(Y,j))$ to be the Chow group
$\textup{CH}^{j-i+d}(X\times_SY)_\QQ:=\textup{CH}^{j-i+d}(X\times_SY)\otimes_\ZZ\QQ$ when $Y$ is purely $d$-dimensional over $S$.
If $Y=\sqcup Y_s$ where each $Y_s$ is connected,
then $\Hom_{CHM'}((X,i),(Y,j))=\oplus \Hom_{CHM'}((X,i),(Y_s,j))$.
Composition is defined by
\begin{eqnarray*}
\Hom_{CHM'}((X_1,i_1),(X_2,i_2))_\QQ&\times& \Hom_{CHM'}((X_2,i_2),(X_3,i_3))_\QQ \\ &\to& \Hom_{CHM'}((X_1,i_1),(X_3,i_3))_\QQ
\end{eqnarray*}
which carries $(U,V)$ to $(p_{1,3})_*(p_{1,2}^*U \cdot p_{2,3}^*V)$
where $p_{i,j}:X_1\times_S X_2\times_S X_3\to X_i\times_SX_j$ is the
natural projection.
The symmetric monoidal structure is given by $(X,i)\otimes (Y,j)=(X\times_SY,i+j)$.
Define the category $CHM$ of Chow motives to be the idempotent completion of $CHM'$.

There is a natural functor $h:\textup{SmPr}_{/S}\to CHM'\to CHM; X\mapsto h(X)=(X,0)$
(any morphism $f:X\to Y$ induces the graph $\Gamma_f$
in $X\times_SY$).
We usually regard objects in $CHM'$ as objects in $CHM$.
We put $L=(S,1)$ and
let $L^{-1}$ be $(S,-1)$.
The symmetric monoidal category $CHM$ is rigid, that is,
 every object is dualizable.

Let $X$ be an abelian scheme over a smooth quasi-projective $k$-scheme $S$.
Suppose that $X$ is relatively $g$-dimensional.
If $\Delta$ denotes the diagonal
of $X \times_SX$, there is a decomposition
$[\Delta] =\Sigma_{i=0}^{2g}\pi_i$ in $\textup{CH}^g(X\times_SX)_{\QQ}$
such that $\pi_i\circ \pi_i=\pi_i$ for any $i$, and $\pi_i\circ \pi_j=0$
for $i\neq j$,
due to Manin-Shermenev, Deninger-Murre, and K\"unnemann
(see e.g. \cite[Section 3]{Kun}).
Let $h_i(X)$ be the direct summand of $h(X)$ in $CHM$
corresponding
to the idempotent morphism $\pi_i$.
We have a natural decomposition
\[
h(X)\simeq \oplus_{i=0}^{2g}h_i(X)
\]
in $CHM$, which is called the Chow-K\"unneth decomposition.
Each $h_i(X)$ is defined to be $\Ker(\textup{id}-\pi_i)$.
If $[\times n]:X\to X$ denotes the multiplication by $n$, then
$[n]$ acts on $h_i(X)$ as the multiplication by $n^i$.
The idempotent morphisms $\pi_0$ and $\pi_{2g}$
are determined by $X\times_Se(S)$ and $e(S)\times_SX$ respectively,
where $e:S\to X$ is a unit of the abelian scheme.
In $CHM$, $h_0(X)$ is a unit object, and $h_{2g}(X)$ is isomorphic to
$L^{g}$.
Moreover, the $i$-fold symmetric power $\Sym^i(h_1(X)) \simeq h_i(X)$.
Here we define $\Sym^i(M)$ to be the kernel of $(\textup{id}-\frac{1}{i!}\Sigma_{\sigma\in \sss_i}\sigma):M^{\otimes i}\to M^{\otimes i}$,
$\sss_i$ is the symmetric group.
For $X,Y$ in $\textup{SmPr}_{/S}$, there is a natural isomorphism
$\textup{Hom}_{\textup{h}(\mathsf{DM})}(M(X),M(Y))\simeq \textup{CH}^d(X\times_SY)_{\QQ}$
where $Y$ is relatively $d$-dimensional (cf. \cite[11.3.8]{CD3}).
Through this comparison
the composition
\[
\textup{Hom}_{\textup{h}(\mathsf{DM})}(M(X),M(Y))\times \textup{Hom}_{\textup{h}(\mathsf{DM})}(M(Y),M(Z))\to \textup{Hom}_{\textup{h}(\mathsf{DM})}(M(X),M(Z))
\]
can be identified with the composition in Chow motives since the comparison commutes with flat pullbacks, intersection product \cite{KY},
and proper push-forwards.
Let $X$ be an abelian scheme over $S$ of relative dimension $g$.
The decomposition
$[\Delta] =\Sigma_{i=0}^{2g}\pi_i$ in $\textup{CH}^g(X\times_SX)_{\QQ}$
induces
a decomposition $M(X)\simeq \oplus_{i=0}^{2g}M_i(X)$
such that the multiplication
$[\times n]:X\to X$ acts on $M_i(X)$ as the multiplication by $n^i$.
By \cite[3.3.1]{Kun}, we have a natural equivalence
$\wedge^{i}(M_1(X)[-1])\simeq M_i(X)[-i]$ for any $i\ge0$.
Here $\wedge^{i}(M)$ is defined to be the kernel of
the kernel of $(\textup{id}-\frac{1}{i!}\Sigma_{\sigma\in \sss_i}\textup{sgn}(\sigma)\sigma):M^{\otimes i}\to M^{\otimes i}$
in the (idempotent complete) homotopy category of $\mathsf{DM}$, where $\textup{sgn}(-)$
indicates the signature.
The $0$-th and $2g$-th components
$M_0(X)$ and $M_{2g}(X)[-2g]$ are isomorphic
to the unit and the Tate object $\QQ(g)$ respectively
(notice that
an isomorphism $M_{2g}(X)\simeq \QQ(g)[2g]$ amounts to
the isomorphism $h_{2g}(X)\simeq L^g$).

\vspace{3mm}

{\it Proof of Theorem~\ref{representationtheorem}.}
Note first that $M_1(X)[-1]$ is a $2g$-dimesional
wedge-finite object in $\mathsf{DM}^\otimes$
in the sense of \cite{DTD}.
That is, $\wedge^{2g}(M_1(X)[-1])\simeq M_{2g}(X)[-2g]$
is invertible in $\mathsf{DM}^\otimes$ and 
$\wedge^{2g+1}(M_1(X)[-1])\simeq 0$.
Thus by \cite[Theorem 4.1]{DTD} (see also \cite[Section 4.2]{DTD})
we have a derived quotient stack $[\Spec A_X/\GL_{2g}]$
and $\QC^\otimes([\Spec A_X/\GL_{2g}])\simeq \mathsf{DM}_X^\otimes$.
The composite $\chi:\QC^\otimes (B\GL_{2g})\to \QC^\otimes([\Spec A_X/\GL_{2g}])\simeq\mathsf{DM}_X^\otimes$
sends the standard representation of $\GL_{2g}$
to $M_1(X)[-1]$.
We remark that 
if $\omega$ denotes a lax symmetric
monoidal
right adjoint functor of $\chi$, then
$A_X\in \CAlg(\QC^\otimes(B\GL_{2g}))$ is obtained as the image $\omega(1_{\mathsf{DM}_X})$.
Here $1_{\mathsf{DM}_X}$ denotes the unit.
\QED

According to \cite[Proposition 4.7]{DTD}
we have the following explicit presentation:

\begin{Proposition}
\label{superexplicit}
Let $Z$ be the set of isomorphism classes of all (finite-dimensional)
irreducible representations of $\textup{GL}_{2g}$. For $z\in Z$, we denote by $V_z$ the corresponding
irreducible representation.
Let $\mathsf{1}_{\mathsf{DM}}$ be a unit of $\mathsf{DM}^\otimes$.
Then there exist
equivalences
\[
A_X \simeq \bigoplus_{z\in Z} V_{z}\otimes \mathsf{Hom}_{\mathsf{DM}}(\chi(V_{z}),\mathsf{1}_{\mathsf{DM}})
\]
in $\QC(B\textup{GL}_{2g})$.
Here $\mathsf{Hom}_{\mathsf{DM}}(-,-)\in \Mod_{\QQ}$ is the hom complex
in $\mathsf{DM}$ (cf. Remark~\ref{homcomp}).
The action of $\textup{GL}_{2g}$ 
on the right hand side is given by the action on $V_{z}$ and the trivial
action on the hom complexes.

\end{Proposition}

\begin{Remark}
\label{homcomp}
Let $\mathsf{Hom}_{\mathsf{DM}}(-,-)$ denote the hom complex which belongs
to $\Mod_\QQ$. Namely, for any $D\in \mathsf{DM}$, we have
the adjoint pair 
\[
D\otimes s(-):\Mod_k \rightleftarrows \CCC:\mathsf{Hom}_{\CCC}(D,-)
\]
where $s$ is the ``$\QQ$-linear structure'' functor $\Mod_\QQ^\otimes\to \mathsf{DM}^\otimes$,
and 
the existence of the right adjoint functor $\mathsf{Hom}_{\mathsf{DM}}(D,-)$
is implied by the adjoint functor theorem and the fact that
$D\otimes s(-)$ preserves small colimits.
If $D=M(X)$, then
$\mathsf{Hom}_{\mathsf{DM}}(M(X),-)$ is a motivic complex of $X$, i.e.,
$H^n(\mathsf{Hom}_{\mathsf{DM}}(M(X),\QQ(i)))$
is the motivic cohomology $H^n(X,\QQ(i))$.
\end{Remark}

\begin{Remark}
\label{superexplicitrem}
Representation theory of $\GL_{2g}$ can be described in terms of
Young diagrams and Schur functor (see e.g. \cite[Section 6]{Ful},
\cite[Section 8]{Ful2}).
For a Young diagram $\lambda$ with $d$ boxes, let $c_\lambda$ be
 the Young symmetrizer in the group algebra
 $\QQ[\sss_d]$ of the symmetric group $\sss_d$
(strictly speaking, we obtain $c_\lambda$
from a Young tableaux whose underlying diagram is $\lambda$).
It gives rise to a direct summand $\QQ[\sss_d]c_\lambda$ of $\QQ[\sss_d]$.
Let $V$ be the standard representation of $\GL_{2g}$, i.e., the
$2g$-dimensional
vector space $V$ equipped with the action of $\Aut(V)=\GL_{2g}$.
Every irreducible representation of $\textup{GL}_{2g}$
is obtained from $V$ by Schur-Weyl's construction:
If $V_z$ is an irreducible
representation, then it is equivalent to 
\[
(V^{\otimes d}\otimes_{\QQ[\sss_d]}\QQ[\sss_d]c_{\lambda})\otimes (\wedge^{2g}V^{\vee})^{\otimes j}
\]
for some $\lambda$ and $j\ge0$.
The (irreducible) representation $V^{\otimes d}\otimes_{\QQ[\sss_d]}\QQ[\sss_d]c_{\lambda}$ is the retract
of $V^{\otimes d}$ that is induced by the retract $\QQ[\sss_d]c_\lambda$
of $\QQ[\sss_d]$ (viewed as a left $\QQ[\sss_d]$-module).
In addition, $V^{\otimes d}\otimes_{\QQ[\sss_d]}\QQ[\sss_d]c_{\lambda}$ is an irreducible representation of $\GL_{2g}$ for any Young diagram (tableaux)
$\lambda$.
By taking highest weights into account,
the set $Z$
can be identified with
the set $\{(\lambda_1,\ldots,\lambda_{2g})\in \ZZ^{\oplus 2g};\lambda_1\ge \lambda_2\ge \ldots \ge \lambda_{2g} \}$:
for an irreducible representation $W$ of $\textup{GL}_{2g}$,
the highest weight is defined as the natural action of the diagonal torus subgroup $T\simeq \mathbb{G}_m^{\times 2g}$ in $\textup{GL}_{2g}$
on the (one dimensional) invariant subspace
$W^{U}$, where $U$ is the subgroup scheme of
the upper triangular invertible matrices with 1's on the diagonal.
When $\lambda_{2g}\ge 0$, $(\lambda_1,\ldots ,\lambda_{2g})$
corresponds to the irreducible representation
$V^{\otimes d}\otimes_{\QQ[\sss_d]}\QQ[\sss_d]c_{\lambda}$
where $\lambda$ has the underlying Young diagram
associated to $(\lambda_1,\ldots,\lambda_{2g})$.
When $0>\lambda_{2g}$,
$(\lambda_1,\ldots,\lambda_{2g})$
corresponds to the irreducible representation
$(V^{\otimes d}\otimes_{\QQ[\sss_d]}\QQ[\sss_d]c_{\lambda})\otimes (\wedge^{2g} V^\vee)^{\otimes (-\lambda_{2g})}$
where $\lambda$ has the underlying Young diagram
$(\lambda_1-\lambda_{2g},\lambda_2-\lambda_{2g},\ldots, \lambda_{2g}-\lambda_{2g})$.
Since $\chi$ is symmetric monoidal and preserves colimits, every $\chi(V_z)$ can be obtained
from $M_1(X)[-1]$ by Schur-Weyl's construction.
Namely, $\chi$ sends $V^{\otimes d}\otimes_{\QQ[\sss_d]}\QQ[\sss_d]c_{\lambda}$ to $(M_1(X)[-1])^{\otimes d}\otimes_{\QQ[\sss_d]}\QQ[\sss_d]c_\lambda$,
and it carries $\wedge^{2g} V^\vee$ to $\QQ(-g)$.
\end{Remark}

\begin{Remark}
\label{tradrel}
We will discuss the symmtric monoidal functor
$\QC^\otimes(B\GL_{2g})\to \DM_X^\otimes$
from the viewpoint of numerical motives.
Suppose that the base scheme is $\Spec k$ with $k$ a perfect field.
In the construction of the category of Chow motives, if we replace
the Chow group $\textup{CH}^{\dim Y}(X\times_kY)_{\QQ}$ by the quotient
$\textup{CH}^{\dim Y}(X\times_kY)_{\QQ}/\sim_{\textup{num}}$ by numerical equivalence,
we obtain another symmetric monoidal category $NM^\otimes$
and the natural symmetric monoidal functor
$CHM^\otimes \to NM^\otimes$
where $\sim_{\textup{num}}$ indicates the numerical equivalence.
We refer to $NM$ as the category of numerical
motives.
By a theorem of Jannsen \cite{Jan}, $NM$ is a semisimple abelian category.

Let $X$ be a smooth projective variety over $k$ and $\mathfrak{h}(X)$ the image in $NM$.
Let $NM_X^\otimes$ be the smallest symmetric monoidal abelian subcategory
of $NM^\otimes$ which contains $\mathfrak{h}(X)$ and its dual $\mathfrak{h}(X)^{\vee}$
(consequently, it is closed under duals).
Suppose that $X$ is an abelian variety over $k$.
The Chow-K\"unneth decomposition induces
the decomposition $\mathfrak{h}(X)\simeq \oplus_{i=0}^{2g}\mathfrak{h}_i(X)$
in $NM$, which we call the motivic decomposition.
We here remark that $\hhhh_i(X)$ is a direct summand
of $\hhhh_1(X)^{\otimes i}$.

In what follows we will modify
the symmetric monoidal structure of $NM_X$ defined above.
We change only the commutative constraint.
In $NM$, the structure morphism of the commutative constraint
\[
\iota:\hhhh_1(X)^{\otimes i}\otimes \hhhh_1(X)^{\otimes j}\to \hhhh_1(X)^{\otimes j}\otimes \hhhh_1(X)^{\otimes i}
\]
is induced by the flip $X^{\times i}\times X^{\times j}\to X^{\times j}\times X^{\times i}; (a,b)\mapsto (b,a)$.
We let $(-1)^{ij}\iota:\hhhh_1(X)^{\otimes i}\otimes \hhhh_1(X)^{\otimes j}\to \hhhh_1(X)^{\otimes j}\otimes \hhhh_1(X)^{\otimes i}$
be a modified commutative constraint.
This modification is extended to retracts of $\hhhh_1(X)^{\otimes i}$
($i \in \ZZ$) in the obvious way.
Unless otherwise stated, from now on we equip $NM_X$ with this modified
symmetric monoidal structure.
When $NM_X$ is equipped with this (modified) symmetric monoidal
structure, there is a natural isomorphism
$\hhhh_i(X)\simeq \wedge^{i}\hhhh_1(X)$, where the latter denotes the wedge
product.

If $\mathsf{DM}_X^\otimes$ (more precisely, its full subcategory consisting
of dualizable objects) admits a motivic $t$-structure,
then heart $MM_X$ is a tannakian category.
Moreover, the fundamental conjecture predicts
that the full subcategory of $MM_X$
spanned by semisimple objects can be identified with the symmetric monoidal
category $NM_X$
of Grothendieck
numerical motives (generated by $\hhhh_1(X)$).
From this conjectural perspective, we should obtain
a symmetric monoidal functor $\mathcal{D}^\otimes(NM_X)\to \mathsf{DM}_X^\otimes$, where we informally denote by $\mathcal{D}^\otimes(NM_X)$
a derived $\infty$-category of $NM_X$.
We can think that
$\QC^\otimes(B\GL_{2g})\to \DM_X^\otimes$ plays a role similar to
a desired functor $\mathcal{D}^\otimes(NM_X)\to \mathsf{DM}_X^\otimes$.
That is, $\QC^\otimes(B\GL_{2g})$ can be viewed as the derived $\infty$-category of ``framed numerical motives''.
\end{Remark}

\subsection{Kimura finite Chow motives}

We conclude this Section by treating Kimura finite Chow motives.
Theorem~\ref{representationtheorem} has a direct generalization
to the Kimura finite case.
Note first that the argument of Theorem~\ref{representationtheorem}
or a direct use of \cite[Theorem 4.1]{DTD} shows the following:

\begin{Corollary}
\label{Kimu}
Let $M$ be a wedge-finite object
in $\mathsf{DM}^\otimes$ (see \cite{DTD}, this condition says that
there is a natural number $d$
such that the exterior-product $\wedge^d M$ is invertible, and $\wedge^{d+1}M\simeq 0$).
Then there exist a derived quotient stack $[\Spec A_M/\GL_{d}]$
and an equivalence $\QC^\otimes([\Spec A_M/\GL_d])\simeq \mathsf{DM}^\otimes_M$. Here $\mathsf{DM}^\otimes_M$ is the smallest presentable stable subcategory
of $\mathsf{DM}^\otimes$,
which contains $M$ and its dual $M^\vee$ and is closed under tensor product.
\end{Corollary}

Suppose that $S$ is irreducible. 
Let $M$ be an object in $\mathsf{DM}^\otimes$.
The notion of wedge-finiteness is closely related to
Kimura finiteness:
An object $M$ in $CHM$ is said to be Kimura finite if
there is a decomposition $M^+ \oplus M^-$
such that the exterior-product
$\wedge^n M^+$ is zero and symmetric product $\textup{Sym}^nM^-$ is zero
for $n>>0$.
If one regards the Chow motif $M$ as an object in
$\mathsf{DM}^\otimes$, then by \cite[Proposition 6.1]{DTD}
$M^+[2n]\oplus M^-[2m+1]$ is wedge-finite in $\mathsf{DM}^\otimes$
for any $n,m\in \ZZ$.

\begin{Example}
The Tate object $\QQ(1)$ is wedge-finite.
Let $\mathsf{DTM}^\otimes$ be the smallest stable
$\infty$-category which contains $\QQ(1)$ and $\QQ(1)^\vee=\QQ(-1)$
and is closed under tensor product and colimits. We call this
category the symmetric monoidal stable $\infty$-category of mixed Tate motives.
By Corollary~\ref{Kimu} there exist $[\Spec A_{\textup{Tate}}/\mathbb{G}_m]$
and $\QC^\otimes([\Spec A_{\textup{Tate}}/\mathbb{G}_m])\simeq \mathsf{DTM}^\otimes$. It is essentially
a representation theorem proven by Spitzweck \cite{Spi}.
See also \cite{Bar}.
\end{Example}

\begin{Remark}
Consider the case when $X=E$ is an elliptic curve with no complex multiplication.
In \cite{Pat}, \cite{KT},
Patashnik and Kimura-Terasoma
have constructed categories of mixed elliptic motives.
Both constructions have built upon the construction of Bloch-Kriz \cite{BK}
and give delicate (and clever) handmade constructions of differential graded algebras
(in \cite{KT}, a
{\it quasi}-DGA has been constructed).
Meanwhile,
(quasi-differential graded, triangulated or abelian) categories 
constructed in \cite{Pat} and \cite{KT} were not compared
with $\mathsf{DM}_{\textup{gm},E}$ 
(they seem by no means easy to compare), whereas in the introduction of
\cite{KT}
the authors hope that the their triangulated
category coincides with the homotopy category of
$\mathsf{DM}_{\textup{gm},E}$.
The comparison between these constructions and ours might be interesting.
The construction of $A_E$ here is somewhat abstract
since we use $\infty$-categorical setting, but it
has an explicit description using motivic cohomology (see Proposition~\ref{superexplicit}, Remark~\ref{superexplicitrem}) and simultaneously it
is adequate for homotopical operations.
\end{Remark}

\begin{Remark}
Let
$\mathsf{R}:\DM^\otimes \to \Mod^\otimes_{\mathbf{K}}$
be
the (singular, \'etale, de Rham, etc.) realization fucntor that is a
symmetric monoidal colimit-preserving functor (see the next Section).
For applications described in the next Section,
we need that
\[
\QC^\otimes(B\GL_d)\to \QC^\otimes ([\Spec A_M/\GL_d])\stackrel{\sim}{\to} \DM_M^\otimes\stackrel{\mathsf{R}}{\to} \Mod_{\mathbf{K}}^\otimes
\]
is the forgetful functor.
This amounts to the condition that 
$\mathsf{R}(M)$ is concentrated in degree zero. Namely, $H_i(\mathsf{R}(M))=0$ for $i\neq 0$.
\end{Remark}

\section{Motivic Galois group of mixed abelian motives}

In this Section, applying the results in Section 3
we will construct the derived motivic Galois group and (underived)
motivic Galois group for $\mathsf{DM}_X$ and prove main results of this paper.
In this Section, the base scheme is $S=\Spec k$ and $k$ is a subfield of $\CC$.

\subsection{}
Let $\mathsf{R}:\mathsf{DM}^\otimes\to \Mod_{\QQ}^\otimes$
be the realization functor associated to singular cohomology
constructed in \cite[5.12]{Tan}.
It is a symmetric monoidal colimit-preserving functor.
For a smooth scheme $X$ over $k$, it carries $M(X)$ to the dual
of a chain complex computing singular cohomology of $X(\CC)$.
Let $\mathsf{DM}_{X}^\otimes \hookrightarrow \mathsf{DM}^\otimes \stackrel{\mathsf{R}}{\to} \Mod_{\QQ}^\otimes$ be the composite where the second functor is the realization functor. We abuse notation and write $\mathsf{R}$ for the composite.
Consider the composite
\[
\QC^\otimes(B\textup{GL}_{2g})\stackrel{\chi}{\to} \mathsf{DM}_{X}^\otimes\stackrel{\mathsf{R}}{\to} \Mod_{\QQ}^\otimes.
\]
By the relative version of adjoint functor theorem \cite[8.3.2.6]{HA},
we have a lax symmetric monoidal right adjoint functor $U$ of $\mathsf{R} \circ \chi$.
Let $\mathsf{1}_{\Mod_{\QQ}}$ be a unit of $\Mod_{\QQ}^\otimes$
and $A'= U(\mathsf{1}_{\Mod_{\QQ}})\in \CAlg(\QC^\otimes(B\textup{GL}_{2g}))$.

\begin{Lemma}
\label{dj}
The stack $[\Spec A'/\textup{GL}_{2g}]$ is isomorphic to
$\Spec \QQ$.
\end{Lemma}

\Proof
Let $\textup{GL}_{2g}=\Spec B$.
We will show that $A'$ has the underlying commutative algebra object
$B$ equipped with the natural action of $\textup{GL}_{2g}$
that
arises from the multiplication of $\textup{GL}_{2g}$.
It suffices to prove that
the composite $\QC^\otimes(B\GL_{2g})\to  \Mod_{\QQ}^\otimes$
is a forgetful functor.
To this end, by \cite[Theorem 3.1]{DTD}
we are reduced to showing that $\mathsf{R}(M_1(X)[-1])$
is equivalent to a vector space placed in degeree zero.
Since $\mathsf{R}(M_1(X))$ is placed in homological degree one, thus
our assertion is clear.
\QED

{\it The construction of a based loop space and motivic Galois groups.}
According to Lemma~\ref{dj} there exists a morphism
\[
\Spec \QQ \simeq [\Spec A'/\textup{GL}_{2g}]\to [\Spec A_X/\textup{GL}_{2g} ].
\]
We refer this morphism as the Betti point.
Then if $\Delta_+$ denotes the category of finite (possibly empty)
linearly ordered sets,
the \v{C}ech nerve $\NNNN(\Delta_+)^{op} \to \Sh(\CAlg_{\mathbb{Q}}^{\textup{\'et}})$
associated to the Betti point gives rise to a derived affine group scheme
$\mathsf{MG}(X)$
over $\mathbb{Q}$
as the group object $\NNNN(\Delta)^{op}\hookrightarrow \NNNN(\Delta_+)^{op} \to \Sh(\CAlg_{\mathbb{Q}}^{\textup{\'et}})$.
In fact, it is a group object of derived affine schemes.
This construction of $\mathsf{MG}(X)$
can be thought of as a $\textup{GL}_{2g}$-equivariant version of bar construction.
For the notion of group objects, \v{C}ech nerves and derived group schemes, we refer the reader
to \cite[6.1.2, 7.2.2.1]{HTT}, \cite[Appendix]{Tan}, \cite[Section 4]{Bar}.
The underlying derived affine scheme of $\mathsf{MG}(X)$
is the fiber product
\[
\Spec \QQ\times_{[\Spec A_X/\textup{GL}_{2g}]}\Spec \QQ
\]
associated to the Betti point.
Let $\Grp(\widehat{\mathcal{S}})$ is the $\infty$-category of group objects
of $\widehat{\mathcal{S}}$.
From the functorial point of view,
as explained in \cite[Appendix]{Tan}, the derived affine group
scheme $\mathsf{MG}(X)$ can be viewed as a functor $\CAlg_{\QQ}\to \Grp(\widehat{\mathcal{S}})$. It is easy to see that $\mathsf{MG}(X)$
is an derived affine group scheme; we put $\mathsf{MG}_X=\Spec B_X$.
We call this derived group scheme the {\it derived motivic Galois
group} for $\mathsf{DM}_X^\otimes$ with respect to
the singular realization.

{\it Automorphism group of the realization functor.}
We abuse notation
and write
$\mathsf{R}:\mathsf{DM}^\otimes_{\textup{gm},X}\to \PMod_{\QQ}^\otimes$
for the restriction of $\mathsf{R}:\mathsf{DM}^\otimes_{X}\to \Mod_{\QQ}^\otimes$.
The automorphism group functor $\Aut(\mathsf{R})$ of 
$\mathsf{R}:\mathsf{DM}^\otimes_{\textup{gm},X}\to \PMod_{\QQ}^\otimes$
is a functor $\CAlg_{\QQ}\to \Grp(\widehat{\mathcal{S}})$
which is informally given by $R\mapsto \Aut(f_R \circ \mathsf{R})$,
where $f_R:\PMod_{\QQ}^\otimes\to \PMod_{R}^\otimes$ is the base change
$(-)\otimes_{\QQ}R$,
and $\Aut(f_R \circ \mathsf{R})$ is $\Map_{\Map^\otimes(\mathsf{DM}_{\textup{gm},X}^\otimes,\PMod_R^\otimes)}(f_R \circ \mathsf{R},f_R \circ \mathsf{R})$.
Roughly speaking,
the group structure of $\Aut(f_R \circ \mathsf{R})\in
\Grp(\widehat{\mathcal{S}})$ is determined by
the composition of symmetric monoidal natural equivalences.
We shall refer the reader to \cite[Section 3]{Tan} for the precise definition.

\begin{Theorem}
\label{motivicGalexp}
The automorphism group functor $\Aut(\mathsf{R})$
of
$\mathsf{R}:\mathsf{DM}^\otimes_{\textup{gm},X}\to \PMod_{\QQ}^\otimes$
is representable by the derived affine group scheme $\mathsf{MG}(X)$.
\end{Theorem}

\Proof
We apply \cite[4.6, 4.9]{Bar} and Theorem~\ref{representationtheorem}
to $\mathsf{DM}^\otimes_{\textup{gm},X}\to \PMod_{\QQ}^\otimes$.
\QED

\begin{Remark}
In \cite{Tan} we have constructed derived motivic Galois groups
which represent the automorphism group functors
of the realization functors in a much more general situation by the abstract machinery of
tannakization.
Here we give an explicit construction of
$\mathsf{MG}(X)$ by means of equivariant bar constructions.
Consequently, Theorem~\ref{motivicGalexp} reveals
the structure of derived motivic Galois groups by means of
$\textup{GL}_{2g}$-equivariant bar constructions in the case of mixed abelian motives.
\end{Remark}

The natural morphism $\Spec \QQ\to \Spec \QQ\times_{B\textup{GL}_{2g}}[\Spec A_X/\textup{GL}_{2g}]\simeq \Spec A_X$ gives rise to its \v{C}ech nerve
as in the case of $\Spec \QQ\to [\Spec A_X/\textup{GL}_{2g}]$, and thus
we have a derived affine group scheme $\mathsf{UG}(X)$ whose underlying derived scheme
is given by
$\Spec (\QQ\otimes_{A_X}\QQ)$.

\begin{Proposition}
\label{derivedpullback}
There is a pullback square of derived group schemes
\[
\xymatrix{
\mathsf{UG}(X) \ar[r] \ar[d] & \mathsf{MG}(X)  \ar[d] \\
\Spec \QQ \ar[r] & \textup{GL}_{2g}.
}
\]
\end{Proposition}

\Proof
There is a pullback square of derived stacks
\[
\xymatrix{
\Spec A_X \ar[r] \ar[d] & [\Spec A_X/\textup{GL}_{2g}]  \ar[d] \\
\Spec \QQ \ar[r] & B\textup{GL}_{2g}.
}
\]
All $\Spec A_X$, $[\Spec A_X/\textup{GL}_{2g}]$
and $B\textup{GL}_{2g}$ are pointed; they are equipped with morphisms
from the final object
$\Spec \QQ$. The above pullback square can be viewed as
a pullback square of pointed derived stacks.
Applying the based loop space functor to this pullback square, we have the desired
pullback square of derived group schemes.
\QED

\begin{Remark}
Let $H_*(A_X)$ be the $\ZZ$-graded algebra associated to $A_X$. 
Using the adjunction and Proposition~\ref{superexplicit}
one can compute the graded
algebra structure of $H_*(A_X)$ in
terms of cup product on motivic cohomology. It is straightforward and left
to the reader.
There is a convergent spectral sequence
\[
E_2^{p,q}:=\textup{Tor}_p^{H_*(A_X)}(\QQ,\QQ)_q \Longrightarrow H_{p+q}(\QQ\otimes_{A_X}\QQ)
\]
where $\QQ$ is placed in degree zero, cf. \cite{HA}.
\end{Remark}

For any commutative differential graded algebra $A$ let $\tau A$
be the quotient
of $A$
by the differential graded ideal generated by
elements of negative 
cohomological degrees.
As explained in \cite[Section 5]{Tan},
the rule $A\mapsto \tau A$ determines
a functor $\CAlg_{\QQ}\to \CAlg_{\QQ}$
which we denote by $\tau$ again.
(In explicit terms, $\CAlg_{\QQ}\to \CAlg_{\QQ}$ sends
$A$ to $\tau A'$ where $A'$ is a cofibrant model of $A$ in the model category
of commutative differential graded algebras \cite{Hin}.) 
We then obtain a pro-algebraic group $MG(X)=\Spec H^0(\tau B_X)$ over $\QQ$
from $\mathsf{MG}(X)$.
The pro-algebraic group $MG(X)$ is the {\it coarse moduli space}
of $\mathsf{MG}(X)$ in the sense that
the natural morphism $\Spec \tau B_X\to MG(X)$
is universal among morphisms to pro-algebraic groups.
See \cite[Definition 7.15]{Bar}, \cite[Section A.4]{Tan}
for the notion of coarse moduli spaces.
(Remark that $\Spec \tau B_X\to MG(X)$ can be viewed as the morphism
$\mathsf{MG}(X)\to MG(X)$ when we regard them as objects in $\Fun(\CAlg^{\textup{dis}}_{\QQ},\Grp(\widehat{\mathcal{S}}))$ since for any $R\in \CAlg^{\textup{dis}}_{\QQ}$ every
$\Spec R \to \mathsf{MG}(X)$ factors though $\Spec \tau B_X$ uniquely up
to a contractible space of choices.)
We refer to $MG(X)$ as the {\it motivic Galois
group} for $\mathsf{DM}_X^\otimes$ with respect to
the singular realization.
By the same argument as in \cite[Section 5, Theorem 5.17]{Tan}, we have
the following coarse representability:
 
\begin{Proposition}
\label{yakimawashi}
Let $K$ be a $\QQ$-field.
Let $\overline{\Aut(\mathsf{R})}(K)$ be the group of
equivalence classes of automorphisms of $\mathsf{R}$,
that is, $\pi_0(\Aut(\mathsf{R})(K))$.
Let $MG(X)(K)$ be the group of $K$-valued points of $MG(X)$.
Then there is a natural isomorphism of groups
\[
MG(X)(K)\simeq \overline{\Aut(\mathsf{R})}(K).
\]
The isomorphisms are functorial among $\QQ$-fields in the natural way.
\end{Proposition}

We will not use the following result which can be proved as in \cite[Proposition 5.19]{Tan}, but
we describe it as a comparison with the traditional line.

\begin{Proposition}
Let $\mathbf{DM}_{\textup{gm},X}$ be the homotopy category of $\mathsf{DM}_{\textup{gm},X}$.
Suppose that $\mathbf{DM}_{\textup{gm},X}$ admits a motivic $t$-structure,
that is, a non-degenerate $t$-structure such that the realization functor
is $t$-exact, and the tensor operation $\mathbf{DM}_{\textup{gm},X}\times \mathbf{DM}_{\textup{gm},X}\to \mathbf{DM}_{\textup{gm},X}$ is $t$-exact.
Then
the heart of the motivic $t$-structure is equivalent to
the category of finite dimensional representations of $MG(X)$
as a symmetric monoidal abelian category.
\end{Proposition}

\subsection{}
\label{decomposition}
We will study the structure of the pro-algebraic group
scheme $MG(X)$.
We slightly change the situation.
We replace the singular realization by the $l$-adic \'etale realization
$\mathsf{R}_{\textup{\'et},\QQ_l}:\mathsf{DM}_{\textup{gm},X}^\otimes\subset \mathsf{DM}^\otimes \to \Mod_{\QQ_l}^\otimes$ constructed
in Section 5; see Proposition~\ref{realization}.
By the same construction as in Lemma~\ref{dj}, we have a base point
called the $l$-adic \'etale point
\[
\Spec \QQ_l\to [\Spec A_X/\textup{GL}_{2g}]
\]
associated to the $l$-adic \'etale realization and the derived affine group
scheme $\mathsf{MG}_{\textup{\'et}}(X)=\Spec B_X=\Spec \QQ_l\times_{[\Spec A_X/\textup{GL}_{2g}]}\Spec \QQ_l$ over $\Spec \QQ_l$
representing the automorphism group functor
$\Aut(\mathsf{R}_{\textup{\'et},\QQ_l}):\CAlg_{\QQ_l}\to \Grp(\widehat{\SSS})$.
Also, we obtain the version $\mathsf{UG}_{\textup{\'et}}(X)=\Spec (\QQ_l\otimes_{A_X}\QQ_l)$ of $\mathsf{UG}(X)$
associated to the base point $\Spec \QQ_l\to \Spec A_X\simeq \Spec \QQ_l\times_{B\textup{GL}_{2g}}[\Spec A_X/\textup{GL}_{2g}]$.
Let us denote by $MG_{\textup{\'et}}(X)$ and $UG_{\textup{\'et}}(X)$
the pro-algebraic group schemes over $\Spec \QQ_l$ which are
associated to $\mathsf{MG}_{\textup{\'et}}(X)$ and $\mathsf{UG}_{\textup{\'et}}(X)$ respectively. Explicitly, $UG_{\textup{\'et}}(X)=\Spec H^0(\tau(\QQ_l\otimes_{A_X}\QQ_l))$ and
$MG_{\textup{\'et}}(X)=\Spec H^0(\tau B_X)$.

\vspace{2mm}

Let us review the motivic Galois group of numerical motives generated by
an abelian variety $X$.
The category $NM_X^\otimes$ (see Remark~\ref{tradrel})
is equipped with
the realization functor of singular cohomology
which is a symmetric monoidal exact functor
$R_B:NM_X^{\otimes}\to \textup{Vect}_{\QQ}^\otimes$,
where $\textup{Vect}_{\QQ}^\otimes$ is the symmetric monoidal
category of $\QQ$-vector spaces.
The functor $R_B$ sends $\hhhh(X)$ to
the dual of the vector space $H^*(X(\CC),\QQ)$
of the singular cohomology of the complex manifold $X(\CC)$.
This functor is faithful and it
makes $NM_X^\otimes$ a $\QQ$-linear neutral
tannakian category (cf. \cite{A}, in this case all Grothendieck standard
conjectures hold).
Its Tannaka dual $MG_{pure}(X)$ is a reductive algebraic subgroup 
in $\textup{GL}_{2g}$.
For a tannakian category equipped with a fiber functor,
by its {\it Tannaka dual} we mean the pro-algebraic (affine) group scheme
that represents the automorphism group of the fiber functor.
We have $MG_{pure}(X)\simeq \Aut(R_B)$ and the closed immersion
$MG_{pure}(X) \hookrightarrow  \textup{GL}_{2g}$
is determined by the action of
$MG_{pure}(X)\simeq \Aut(R_B)$ on $R_B(\hhhh_1(X))\simeq \QQ^{\oplus 2g}$.
For example, if $X$ is an elliptic curve with no complex multiplication
(i.e., $\End(X\otimes_{k}\CC)=\ZZ$), its Tannaka 
dual is $\operatorname{GL}_2$.
The object $\hhhh_1(X)$ corresponds to the natural
action on $\operatorname{GL}_2=\operatorname{GL}(H^1(X(\CC),\QQ)^{\vee})$ on $H^1(X(\CC),\QQ)^{\vee}$.
Let $MG_{pure}(X)$ be the Tannaka dual of $NM_X^\otimes$
and put $MG_{pure}(X)_{\QQ_l}=MG_{pure}(X)\otimes_{\QQ}\QQ_l$.

\vspace{2mm}

We are in a position to state our results:


\begin{Theorem}
\label{unconditional}
Let $X$ be an abelian variety of dimension $g$
over a number field $k$.
Suppose either (i) or (ii) or (iii):
\begin{enumerate}
\renewcommand{\labelenumi}{(\roman{enumi})}
\item $\End(X\times_k\bar{k})=\ZZ$. Assume neither that (a) $2g$ is a $n$-th power for any odd number $n>1$, nor (b) $2g$ is of the form $\begin{pmatrix} 2n \\ n \end{pmatrix}$
for any odd number $n>1$.

\item $X$ is an elliptic curve.
If $X$ has complex multiplication, we suppose that the complex
multiplication is defined over $k$.

\item $X$ is a simple abelian variety of prime dimension $\dim X=p$ (including $1$)
which has complex multiplication, i.e., $\End(X)\otimes_\ZZ\QQ$ is a
CM-field of degree $2p$.
Suppose further that $X$ is absolutely simple,
i.e., it is also simple after the base change to an algebraic closure.

\end{enumerate}
Then there exists an exact sequence of pro-algebraic group schemes
\[
1\to UG_{\textup{\'et}}(X)\to MG_{\textup{\'et}}(X) \to MG_{pure}(X)_{\QQ_l} \to 1.
\]
Moreover, $UG_{\textup{\'et}}(X)$ is a connected pro-unipotent group scheme over $\QQ_l$.
\end{Theorem}

\begin{Corollary}
\label{unconditional2}
\begin{enumerate}
\renewcommand{\labelenumi}{(\roman{enumi})}
\item There is an isomorphism of affine group schemes \[
MG_{\textup{\'et}}(X) \simeq UG_{\textup{\'et}}(X)\rtimes MG_{pure}(X)_{\QQ_l}.
\]

\item $UG_{\textup{\'et}}(X)$ is the unipotent radical of $MG_{\textup{\'et}}(X)$, i.e., the maximal normal unipotent closed subgroup.

\item $MG_{pure}(X)_{\QQ_l}$ is the reductive quotient.
\end{enumerate}
\end{Corollary}

\Proof
In virtue of a theorem of Hochschild-Mostow \cite{HM} (see also \cite[Section 1]{Pr}),
there is a section of $MG_{\textup{\'et}}(X) \to MG_{pure}(X)_{\QQ_l}$
which is unique up to conjugation by $UG_{\textup{\'et}}(X)$.
It gives rise to $MG_{\textup{\'et}}(X) \simeq UG_{\textup{\'et}}(X)\rtimes MG_{pure}(X)_{\QQ_l}$. Other claims are clear.
\QED

We will prove Theorem~\ref{unconditional}.
From now on, for simplicity
we write $\mathsf{MG}$, $\mathsf{UG}$, $MG$, $UG$, $MG_{pure}$ and $\textup{GL}_{2g}$
for $\mathsf{MG}_{\textup{\'et}}(X)$, $\mathsf{UG}_{\textup{\'et}}(X)$,
$MG_{\textup{\'et}}(X)$, $UG_{\textup{\'et}}(X)$
$MG_{pure}(X)$ and $\textup{GL}_{2g}\times_{\Spec \QQ}\Spec \QQ_l$
respectively.

\vspace{2mm}

{\it Proof of the exactness in Theorem~\ref{unconditional}.}
The natural sequence $\mathsf{UG}\to \mathsf{MG}\to \textup{GL}_{2g}$ in Proposition~\ref{derivedpullback}
induces $UG \stackrel{s}{\to} MG\stackrel{t}{\to} \textup{GL}_{2g}$.
We first claim that for any $\QQ_l$-field $K$, the
sequence of the groups of $K$-valued points
\[
1\to UG(K)\to MG(K)\to \textup{GL}_{2g}(K)
\]
is exact.
To see this, consider the sequence of mapping spaces
\[
\Map(\Spec K,\Spec A_X)\to \Map(\Spec K, [\Spec A_X/\textup{GL}_{2g}])\to \Map(\Spec K,B\textup{GL}_{2g}),
\]
where $\Spec A_X\to [\Spec A_X/\textup{GL}_{2g}]\to B\textup{GL}_{2g}$ is given in the proof in Proposition~\ref{derivedpullback},
and the mapping spaces are taken in
$\Fun(\CAlg_{\QQ_l},\widehat{\SSS})$.
Let $\Spec \QQ_l\to [\Spec A_X/\textup{GL}_{2g}]$ be the morphism arising from
the $l$-adic \'etale realization functor.
For each $\Spec K\to \Spec \QQ_l$, the composition determines
$K$-valued base points of $[\Spec A_X/\textup{GL}_{2g}]$,
$\Spec A_X$ and $B\textup{GL}_{2g}$.
Again by Proposition~\ref{derivedpullback} and the homotopy exact sequence,
we have an exact sequence
\begin{eqnarray*}
\pi_2(\Map(\Spec K,B\textup{GL}_{2g}))\to \pi_1(\Map(\Spec K,\Spec A_X))\simeq \pi_0(\mathsf{UG}(K)) \\
\to \pi_1(\Map(\Spec K, [\Spec A_X/\textup{GL}_{2g}])) \simeq \pi_0(\mathsf{MG}(K))\to \pi_1(\Map(\Spec K,B\textup{GL}_{2g}))
\end{eqnarray*}
where the base points of the homotopy groups
are defined above.
Note that
\[
\pi_1(\Map(\Spec K,B\textup{GL}_{2g}))\simeq \textup{GL}_{2g}(K)
\]
and $\pi_2(\Map(\Spec K,B\textup{GL}_{2g}))\simeq 1$.
Proposition~\ref{yakimawashi} implies that $\pi_0(\mathsf{MG}(K))\simeq MG(K)$.
The similar argument (cf. \cite[Theorem 5.18]{Tan})
also shows $\pi_0(\mathsf{UG}(K))\simeq UG(K)$.
Thus we obtain the desired exact sequence.

Next we will prove that a sequence of pro-algebraic group
schemes $1\to UG\stackrel{s}{\to} MG \stackrel{t}{\to} \textup{GL}_{2g}$ is exact.
We first observe that $s:UG\to MG$ is injective.
Let $\Ker(s)=\Spec P$ be the affine (pro-algebraic) group scheme of the kernel of $s$.
Then $P$ is of the form of a filtered colimit $\varinjlim_{\lambda} P_{\lambda}$
where each $P_i$ is a finitely generated commutative Hopf
subalgebra of $P$. Since we work over the field $\QQ_l$ of
characteristic zero, each algebraic group $\Spec P_i$ is reduced, and thus so is $\Spec P$.
As proved above, the group $\Ker(s)(K)$ of $K$-valued points
is trivial for any $\QQ_l$-field $K$.
Hence we conclude that the unit 
morphism
$\Spec \QQ_l\to \Ker(s)$ is an isomorphism, that is, $\Ker(s)$ is trivial.
We then prove that the injective homomorphism
$UG\to \Ker(t)$ is a surjective morphism of affine
group schemes, where $\Ker(t)$ is the kernel of $t$.
Suppose that $UG\subset \Ker(t)$ is a proper closed subgroup scheme.
Since $\Ker(t)$ is also reduced, it gives rise to
a contradiction that $UG(K)\to \Ker(t)(K)$ is bijective for any $\QQ_l$-field
$K$. Thus we see that $UG\simeq \Ker(t)$.

Finally, we will prove that the image of
$t:MG\to \textup{GL}_{2g}$ is isomorphic $MG_{pure}$.
The morphism $t$ factors into a sequence of homomorphisms $MG\stackrel{t'}{\to} G\stackrel{t''}{\to} \textup{GL}_{2g}$ of affine group schemes
such that $t'$ is surjective and $t''$ is a closed immersion. It suffices to show that
$G\simeq MG_{pure}$. For this purpose, consider the action of the absolute Galois group $\Gamma=\Gal(\bar{k}/k)$ on $\mathsf{R}_{\textup{\'et},\QQ_l}:\mathsf{DM}_{\textup{gm},X}^\otimes\to \PMod_{\QQ_l}^\otimes$.
According to Proposition~\ref{yakimawashi} and Proposition~\ref{realization},
it gives rise to a homomorphism
\[
\Gamma\to MG(\QQ_l)\simeq \pi_0(\Aut(\mathsf{R}_{\textup{\'et},\QQ_l})(\QQ_l)).
\]
Let $\textup{fVect}_{\QQ}^\otimes(\textup{GL}_{2g})$
be
the symmetric monoidal abelian category of
finite dimensional (rational) representations of $\GL_{2g}$.
Let $\iota:\textup{fVect}_{\QQ}^\otimes(\textup{GL}_{2g}) \to \mathsf{DM}_{\textup{gm},X}^\otimes$ be the symmetric monoidal functor
defined as the composite $\textup{fVect}_{\QQ}^\otimes(\textup{GL}_{2g})
\hookrightarrow \QC^\otimes(BGL_{2g})\stackrel{\chi}{\to} \DM_X^\otimes$
where the right functor is given in Theorem~\ref{representationtheorem}
and $\textup{fVect}_{\QQ}^\otimes(\textup{GL}_{2g})
\hookrightarrow \QC^\otimes(B\GL_{2g})$ is the natural
inclusion into the full subcategory of ``objects placed in degree zero''.
Then the composition with $\textup{fVect}_{\QQ}^\otimes(\textup{GL}_{2g}) \stackrel{\iota}{\to} \mathsf{DM}_{\textup{gm},X}^\otimes\stackrel{\mathsf{R}_{\textup{\'et},\QQ_l}}{\longrightarrow} \PMod_{\QQ_l}^\otimes$ induces
\[
\Gamma \to MG(\QQ_l)\to \textup{GL}_{2g}(\QQ_l).
\]
Here $\textup{GL}_{2g}(\QQ_l)$ denotes the group of $\QQ_l$-valued points
on $\textup{GL}_{2g}$, which we identify
with the automorphism group of the composite $\textup{fVect}^\otimes_{\QQ}(\textup{GL}_{2g})\to \textup{fVect}_{\QQ_l}^\otimes\subset \PMod_{\QQ_l}^\otimes$, where $\textup{fVect}_{\QQ_l}^\otimes$
is the symmetric monoidal category of finite dimensional $\QQ_l$-vector spaces
(note here that by the classical tannaka duality $\textup{fVect}^\otimes_{\QQ}(\textup{GL}_{2g})\to \textup{fVect}_{\QQ_l}^\otimes$ is the composition of the forgetful functor and the base change).
The standard representation of $\textup{GL}_{2g}$
maps to $M_1(X)[-1]$
(cf. Theorem~\ref{representationtheorem})
and the realization functor sends it to $\mathsf{R}_{\textup{\'et},\QQ_l}(M_1(X)[-1])\simeq \QQ_l^{\oplus 2g}$ by Proposition~\ref{realization},
where $\QQ_l^{\oplus 2g}$ is considered to be a complex placed in
degree zero. The image of $g\in MG(\QQ_l)$
in $\textup{GL}_{2g}$ can be viewed as the action of $g$
on $\mathsf{R}_{\textup{\'et},\QQ_l}(M_1(X)[-1])\simeq \QQ_l^{\oplus 2g}$.
Note that an invertible sheaf on $X\times_k X$ that gives rise to a polarization on $X$ induces an anti-symmetric morphism
$M_1(X)[-1]\otimes M_1(X)[-1] \to \QQ(1)$.
It yields a symplectic form $Q:\mathsf{R}_{\textup{\'et},\QQ_l}(M_1(X)[-1])\otimes \mathsf{R}_{\textup{\'et},\QQ_l}(M_1(X)[-1]) \to \mathsf{R}_{\textup{\'et},\QQ_l}(\QQ(1))\simeq \QQ_l$.
By replacing $\mathsf{R}_{\textup{\'et},\QQ_l}(M_1(X)[-1])\simeq \QQ_l^{\oplus 2g}$ by an appropriate isomorphism if necessary, $Q$ induces
a standard symplectic form $\QQ_l^{\oplus 2g}\otimes \QQ_l^{\oplus 2g}\simeq \mathsf{R}(M_1(X)[-1])\otimes \mathsf{R}(M_1(X)[-1]) \to \QQ_l$.
Taking account into the compatibility of symmetric monoidal
natural transformations, we see that
the image of $t$ is contained in $\textup{GSp}_{2g}\subset\textup{GL}_{2g}$, where $\textup{GSp}_{2g}$ is the general symplectic group, that is,
subgroup that consists of symplectic similitudes.
The group of symmetric monoidal
natural equivalences of $\textup{fVect}^\otimes_{\QQ}(\textup{GL}_{2g})\to \textup{fVect}_{\QQ_l}^\otimes$ is $\textup{GL}_{2g}(\QQ_l)$
which is naturally identified with the group of automorphisms of $\mathsf{R}_{\textup{\'et},\QQ_l}(M_1(X)[-1])\simeq \QQ_l^{\oplus 2g}$.
Therefore the image of $g\in \Gamma$ in $\textup{GL}_{2g}(\QQ_l)$ corresponds to
the Galois action of $g$ on $\mathsf{R}_{\textup{\'et},\QQ_l}(M_1(X)[-1])\simeq \QQ_l^{\oplus 2g}$ (induced by $\Gamma\to \Aut(\mathsf{R}_{\textup{\'et},\QQ_l})$ in Proposition~\ref{realization}).
Note that $\Gamma$-module $\mathsf{R}_{\textup{\'et},\QQ_l}(M_1(X)[-1])$
is the Tate module $T_l(X)\otimes_{\ZZ_l}\QQ_l$, which is a dual of the \'etale cohomology $H^1_{\textup{\'et}}(\overline{X},\ZZ_l)\otimes_{\ZZ_l}\QQ_l$ where
$\overline{X}$ is the base change $X\times_k\bar{k}$ to the algebraic
closure $\bar{k}$ of $k$.
Namely, the above composite $\Gamma\to \textup{GL}_{2g}(\QQ_l)$
is the Galois representation
$\Gamma\to \Aut(T_l(X)\otimes_{\ZZ_l}\QQ_l)$.
In the case of (i),
by virtue of a theorem on Galois representation \cite[Chap. IV]{Se1}
\cite{Tank}
\cite[Theorem 5.14]{Pi}
 (the affirmative solution of Mumford-Tate conjecture in this case),
$\textup{GSp}_{2g}$
is the smallest algebraic subgroup
of $\textup{GL}_{2g}$ which contains the image of the Galois representation. Hence we deduce that $G=\textup{GSp}_{2g}$. For the case (i),
$\textup{GSp}_{2g}=MG_{pure}$. This prove the case (i).
Next we consider the case (ii).
By (i), we may suppose that $X$ has complex multiplication
defined over $k$. This case is reduced to the case (iii).
Next consider the case (iii).
Set $D:=\End(X)\otimes_\ZZ\QQ$, $V:=\mathsf{R}_{\textup{\'et},\QQ_l}(M_1(X)[-1]))$, and let $\textup{GL}_D(V)$ be the algebraic subgroup of $\textup{GL}_{2p}=\textup{GL}(V)$ consisting of $D$-linear automorphisms ($D$ acts on $V$
in the obvious way).
As in the observation in (i), taking account into the compatibility
with the polarization and $\End(X)\otimes_\ZZ\QQ$ (regarded as
morphisms in the homotopy category of $\mathsf{DM}$)
we see that the image of $MG\to \textup{GL}_{2p}$
is contained in $\textup{GSp}_{D,2p}:=\textup{GSp}_{2p}\cap \textup{GL}_D(V)$.
According to \cite{Tank2} the Mumford-Tate group $MT(X)_{\QQ_l}=MT(X)\times_{\Spec \QQ}\Spec \QQ_l$
of $X$ coincides with
$\textup{GSp}_{D,2p}$ (we use the condition).
On the other hand, as a consequence of complex multiplications
(see e.g. \cite{Yu})
the Zariski closure $G_l$ of the image of the Galois representation
$\Gamma\to MG\to \textup{GL}_{2p}$ is $MT(X)_{\QQ_l}$ through the comparison isomorphism
between $l$-adic \'etale cohomology and singular cohomology. The motivic Galois group $MG_{pure}$ has the relation $G_l\subset MG_{pure}\subset \textup{GSp}_{D,2p}$.
Hence we conclude that the image of $MG\to \textup{GL}_{2p}$ is
$MG_{pure}=\textup{GSp}_{D,2p}$.
\QED

We will present two proofs of the second assertion in Theorem~\ref{unconditional}; $UG$ is a connected pro-unipotent group scheme. For this, we do not use the condition
that $X$ is an abelian variety.
The first proof is more direct than the second one.
But the second proof is powerful. It yields
explicit structures of $\mathsf{UG}$ and $UG$ (as schemes).

\vspace{3mm}

{\it Proof of the second assertion in Theorem~\ref{unconditional}.}
Let $K$ be a field of characteristic zero.
Let $A$ be a commutative differential graded $K$-algebra
$A$ such that $H^i(A)=0$ for $i< 0$, which we regard as an object in $\CAlg_{K}$.
We first recall that the natural morphism $H^0(A)\to A$ in $\CAlg_{K}$
is decomposed into a sequence 
\[
H^0(A)=A(0)\to A(1)\to A(2)\to \cdots \to A(i)\to \cdots
\]
in $\CAlg_{K}$
such that the natural morphism $\varinjlim_i A(i)\to A$ is an equivalence
(here
$H^0(A)$ is viewed as a discrete commutative algebra object which belongs to
$\CAlg_{K}$).
Moreover, for each $i\ge0$ the morphism $A(i)\to A(i+1)$
fits in a pushout diagram of the form
\[
\xymatrix{
\Sym^*(C)\ar[d] \ar[r] & A(i) \ar[d] \\
K \ar[r] & A(i+1)
}
\]
in $\CAlg_K$, where $\Sym^*(C)$ is a free commutative algebra
associated to some $C\in \Mod_{K}$.
This is the well-known fact in the theory of commutative
differential graded
algebras and rational homotopy theory (see e.g. \cite{Hin}).
Here we refer the reader to the proof of \cite[VIII, 4.1.4]{DAGn}
(strictly speaking, in {\it loc. cit.} the case of $H^0(A)=K$ is treated, but the
argument reveals that our claim holds).
The base point $\Spec \QQ_l\to \Spec A$ induces 
the base points of derived affine schemes in the sequence
\[
\Spec A\to \cdots \to Z(i)=\Spec A(i)\to\cdots \to Z(1)=\Spec A(1)\to Z(0)=\Spec A(0)
\]
by the composition.

Let $\pi_n(\Spec A_X):\CAlg_{K}^{\textup{dis}}\to \Grp$
be a group-valued functor
which to any $R\in \CAlg_{K}^{\textup{dis}}$
assigns the homotopy group
\[
\pi_n(\Map_{\Aff_{K}}(\Spec R,\Spec A_X))
\]
with respect to the base point $\Spec R\to \Spec K \to \Spec A_X$.
Here $\Aff_{K}$ is nothing but $\CAlg_{K}^{op}$.
Next by using above sequence
we will show that for $\pi_n(\Spec A_X)$ is represented by 
a pro-unipotent group scheme over $K$ for $n>0$.
To this end,
we define $\pi_n(Z(i))$ in a similar way.
We first prove that $\pi_n(Z(i))$ is represented by
a pro-unipotent group scheme over $\QQ_l$.
We use the induction on $i$.
The case $i=0$ is obvious since the mapping space
$\Map_{\Aff_{K}}(\Spec R,Z(0))$
is discrete, i.e., $0$-truncated.
Suppose that our claim holds for the case $i\le r$.
Consider the pullback diagram
\[
\xymatrix{
Z(r+1) \ar[r] \ar[d] & \Spec K \ar[d] \\
Z(r) \ar[r] & \Spec \Sym^*C
}
\]
in $\Aff_{K}$.
By homotopy exact sequence, we have an exact sequence
\begin{eqnarray*}
\pi_{n+1}(Z(r))(R)&\stackrel{\eta}{\to}& \pi_{n+1}(\Spec \Sym^*C)(R) \\
&\to& \pi_n(Z(r+1))(R)\to \pi_n(Z(r))(R) \stackrel{\mu}{\to} \pi_n(\Spec \Sym^*C)(R).
\end{eqnarray*}
Note that by the assumption $\pi_{n+1}(Z(r))$ and $\pi_n(Z(r))$
are pro-unipotent group schemes over $K$.
Since $\pi_{n+1}(Z(r))(R)$ is commutative and $K$ is characteristic zero,
there is a $K$-verctor space $V$
such that the functor $\pi_{n+1}(Z(r))$
is given by $R\mapsto \Hom_K(V,R)$.
(Namely, in charateristic zero every commutative pro-unipotent group
scheme is of the form $\mathbb{G}_a^{I}=\Spec K[I]$ for some set $I$.)
Here $\Hom_{K}(-,-)$ indicates the group of homomorphism of $K$-vector spaces.
Unwinding the definitions, $\pi_{n+1}(\Spec \Sym^*C)$
is a commutative pro-unipotent group scheme given by
$R\mapsto \Hom_K(H^{n+1}(C),R)$ as a group-valued functor.
The homomorphisms $\Hom_K(V,R)\to \Hom_K(H^{r+1}(C),R)$
which are functorial with respect to $R$ determines 
a homomorphism $\xi:H^{r+1}(C)\to V$ of $K$-vector spaces.
The cokernel $\Coker(\eta)$
of $\pi_{n+1}(Z(r))\to \pi_{n+1}(\Spec \Sym^*C)$
is a commutative pro-unipotent group
scheme given by $R\mapsto \Hom_K(\Ker(\xi),R)$.
We will show that the kernel $\Ker(\mu)$ of
$\mu:\pi_n(Z(r)) \to \pi_n(\Spec \Sym^*C)$ is a pro-unipotent group scheme.
Note that $\pi_n(\Spec \Sym^*C)$ is also a pro-algebraic group scheme;
we put $\pi_n(\Spec \Sym^*C)=\Spec D$
(moreover, it is commutative and pro-unipotent).
Put $\Spec Q=\Spec \varinjlim_{\lambda\in \Lambda}Q_\lambda=\pi_n(Z(r))$
where $\{Q_\lambda\}_{\lambda\in \Lambda}$ is a filtered diagram of finitely generated commutative Hopf subalgebras
in $Q$ such that each $\Spec Q_\lambda$ is a unipotent algebraic group scheme.
The kernel $\Ker(\mu)$ is $\Spec K\times_{\Spec D}\Spec Q=\Spec Q/I$
which is the fiber over the unit $\Spec K\to \Spec D$.
Hence $\Ker(\mu)$ is $\varprojlim_{\lambda\in \Lambda}\Spec Q_\lambda/(Q_\lambda\cap I)$. Here $\Spec Q_\lambda/(Q_\lambda\cap I)$ is a closed subgroup scheme
in $\Spec Q_\lambda$ which is the image of
$\Ker(\mu)$ in $\Spec Q_\lambda$.
Note that a closed subgroup scheme of a unipotent algebraic group
is unipotent. We then conclude that
$\Ker(\mu)$ is a pro-unipotent group scheme.
Using the surjective map
$\pi_n(Z(r+1))(Q/I)\to \Ker(\mu)(Q/I)$ we have its section $\Ker(\mu)\to \pi_n(Z(r+1))$. Thus $\pi_n(Z(r+1))\simeq \Coker(\eta)\times \Ker(\mu)$ as a set-valued functor. It follows that $\pi_n(Z(r+1))$ is an affine group scheme.
Moreover, it is an extention of pro-unipotent group scheme
by a pro-unipotent group scheme. Consequently, $\pi_n(Z(r+1))$
is represented by a pro-unipotent group scheme over $K$.
Finally, it remains to prove that the natural
map $\pi_n(\Spec A_X)\to \lim_i\pi_n(Z(i))$ is
an equivalence.
Taking account of Milnor exact
sequence associated to $\Spec A_X\simeq \varprojlim_iZ(i)$
we are reduced to showing that $\varprojlim^1\pi_n(Z(i))(R)=0$ for $n\ge 2$. 
Since $\pi_n(Z(i))$ is commutative and unipotent,
we have $\pi_n(Z(i))(R)\simeq \Hom_K(V(i),R)$ for some $K$-vector
space $V(i)$ ($\pi_n(Z(i))$ is given by $\Hom_K(V(i),-)$).
It is enough to show that $\varprojlim^1\Hom_K(V(i),R)=0$
for any $K$-algebra $R$.
By definitions, $\varprojlim^1\Hom_K(V(i),R)\simeq\pi_0(\varprojlim_i\Map_{\Mod_{K}}(V(i),R[1]))=\pi_0(\Map_{\Mod_{K}}(\varinjlim_iV(i),R[1]))=\Ext^1(\varinjlim_i V(i),R)$.
Since $\varinjlim_i V(i)$ is obviously
a $K$-vector space, $\Ext^1(\varinjlim_i V(i),R)$ is trivial.
Then the isomorphism $\pi_n(\Spec A)\simeq \varprojlim_i\pi_n(Z(i))$
shows that $\pi_n(\Spec A)$ is represented by
a pro-unipotent group scheme over $K$.

Now we will show that there is a natural equivalence
$\pi_1(\Spec \tau A_X)\simeq UG$.
Note first the natural equivalences
$\pi_1(\Spec \tau A_X)\simeq\pi_1(\Spec A_X)\simeq \pi_0(\Spec {\QQ_l}\otimes_{A_X}{\QQ_l})=\pi_0(\mathsf{UG})$ as group-valued functors $\CAlg_{\QQ_l}^{\textup{dis}}\to \Grp$.
By compositions with the natural inclusions $\Grp\to \Grp(\widehat{\SSS})$ and $\CAlg_{\QQ_l}^{\textup{dis}}\to \CAlg_{\QQ_l}$
let us consider $\pi_1(\Spec A_X)$ and $\mathsf{UG}$
to be functors
$\CAlg_{\QQ_l}^{\textup{dis}}\to \Grp(\widehat{\SSS})$.
Thus the natural morphism $\mathsf{UG}\to \pi_1(\Spec A_X)$
given by $\mathsf{UG}(R)\to \pi_0(\mathsf{UG}(R))$
is universal among morphisms functors taking values in discrete groups.
On the other hand, $\mathsf{UG}\to UG$ is universal among morphisms
to affine group schemes over $\QQ_l$
in $\CAlg_{\QQ_l}^{\textup{dis}}\to \Grp(\widehat{\SSS})$.
Since $\pi_1(\Spec \tau A_X)$
is an affine group scheme, we see that $\pi_1(\Spec \tau A_X)\simeq UG$.
This implies that $UG$ is a pro-unipotent group scheme over $\QQ_l$.

Finally, we remark that $UG$ is connected.
If we set $UG=\Spec R$, then $R$ is a filtered colimit
$\varinjlim_{\lambda} R_{\lambda}$
such that $R_{\lambda}$ is a commutaive Hopf subalgebra and $\Spec R_{\lambda}$
is a unipotent algebraic group. Since we work over the field of characteristic
zero, thus $\Spec R_{\lambda}$ is connected
(if otherwise, the (nontrivial) quotient by the identity component is unipotent
and we have a contradiction).
Therefore $R$ has non idempotent element. This means $UG$ is connected.
\QED

Next we will give the second proof. We use a version of
Hochschild-Kostant-Rosenberg theorem (see \cite{HKR}).
The following generalization is taken from \cite{BN}.

\begin{Proposition}
\label{BNHKR}
Let $K$ be a field of characteristic zero
and let $A \in \CAlg_K$. Let $Y=\Spec A$.
Let $\mathbb{L}_{A/K}$ be the cotangent complex which belongs to
$\Mod_A$.
Then there is an equivalence
\[
LY:=\Spec (A\otimes_{A\otimes_K A}A)\simeq \Spec (\textup{Sym}_A^*(\mathbb{L}_{A/K}[1])),
\]
where $\textup{Sym}_A^*(\mathbb{L}_{A/K}[1])$ is the free commutative differential
graded algebra over $A$. That is, $\textup{Sym}_A^*(-)$ is the left adjoint
of the forgetful functor $\CAlg_{A}\to \Mod_A$.
\end{Proposition}

Here $\Spec (A\otimes_{A\otimes_K A}A)$ is the free loop space of $\Spec A$.
Since $A\otimes_{A\otimes_K A}A$ is the tensor product
$A\otimes S^1$ with the circle $S^1$, we regard $Y\times_{Y\times Y}Y=\Spec (A\otimes_{A\otimes_K A}A)$
as the derived scheme that represents the functor $\CAlg_K\to \mathcal{S}$
informally given by $R \mapsto \Map_{\mathcal{S}}(S^1, Y(R))$, where $Y(R)$ denotes the space of $R$-valued points, i.e.,
$\Map_{\CAlg_K}(A,R)$.

\vspace{1mm}

{\it Sketch of the proof of Proposition~\ref{BNHKR}.}
We will sketch the argument by Ben-Zvi-Nadler (see \cite{BN} for details).
In \cite[Section 3, Proposition 4.4]{BN}, $A$ is supposed to be connective, i.e., $H_i(A)=0$ for $i<0$.
But arguments in {\it loc. cit.} work for every $A\in \CAlg_K$.
Let $\widehat{\CAlg}_K$ be the $\infty$-category of commutative algebra
objects in $A$-module objects in the enlarged universe $\mathbb{V}$.
By Yoneda embedding we have the fully faithful inclusion
$\Aff_K=(\CAlg_K)^{op}\subset \textup{PSh}:=\Fun(\CAlg_K,\widehat{\mathcal{S}})$.
Here $\widehat{\mathcal{S}}$ is the $\infty$-category of spaces that
belong to $\mathbb{V}$.
Let us regard $S^1$ as the constant functor $\CAlg_K\to \widehat{\mathcal{S}}$
that takes the value $S^1$.
We take the limit $\mathcal{O}(S^1):=\varprojlim _{\Spec R\to S^1}R$ in $\widehat{\CAlg}_K$
where $\Spec R\to S^1$ run through $\mathbb{V}$-small $\infty$-category
$(\Aff_K)_{/S^1}$.
We have $\mathcal{O}(S^1)\simeq C^*(S^1,K) \simeq  K\oplus K[-1]$ where the right hand side is the square zero extesion by $K[-1]$.
Since $S^1 \simeq \varinjlim_{\Spec R\to S^1}\Spec R$ in $\textup{PSh}$, there is a natural morphism $S^1\to \Spec \mathcal{O}(S^1)$
(which is called the affinization of $S^1$ in \cite{BN}).
Let $\Map(\Spec \mathcal{O}(S^1),\Spec A)$ denote the functor
$\CAlg_K\to \widehat{\mathcal{S}}$
informally given by $R\mapsto \Map_{\textup{PSh}}(\Spec \mathcal{O}(S^1)\otimes_KR,\Spec A)$.
The composition induces a homotopy equivalence
$\Map_{\textup{PSh}}(\Spec \mathcal{O}(S^1)\otimes_KR,\Spec A) \simeq \Map_{\textup{PSh}}(S^1\times \Spec R,\Spec A)$ for any $R, A\in \CAlg_K$.
Namely, $\Map(\Spec \mathcal{O}(S^1),\Spec A)\simeq LY$.
For any $f:\Spec R\to Y$, we note the equivalence from the universal property of
cotangent complex:
\[
\Map_{\Mod_R}(\mathbb{L}_{A/K}\otimes_A R,R[-1])\simeq \Map_{\Aff_K}(\Spec R\oplus R[-1],Y)\times_{\Map_{\Aff_K}(\Spec R,Y)}\{f\}
\]
where 
$\Map_{\Aff_K}(\Spec R\oplus R[-1],Y)\to\Map_{\Aff_K}(\Spec R,Y)$
is induced by the first projection $R\oplus R[-1]\to R$.
The left hand side can be identified with 
\[
\Map_{(\Aff_K)_{/\Spec A}}(\Spec R, \Spec (\textup{Sym}_A^*(\mathbb{L}_{A/K}[1]))).
\]
Thus we obtain the desired equivalence $LY\simeq \Spec (\textup{Sym}_A^*(\mathbb{L}_{A/K}[1])))$ over $Y$.
\QED

Now consider the motivic algebra $A=A_X$.
Let $\Spec K\to \Spec A_X$ be the base point
which is induced by either the Betti point or the
l-adic \'etale point $\Spec K\to [\Spec A_X/\GL_{2g}]$
($K=\QQ$ or $\QQ_l$).
The derived affine group scheme
$\mathsf{UG}$ is the based loop space $\Omega_*Y\simeq \Spec K\otimes_{A_X}K\simeq LY\times_Y\Spec K$. The projection $LY\to Y$
is induced by the evaluation at a fixed point on $S^1$.
Then
we obtain equivalences of derived schemes
\[
\mathsf{UG}\simeq \Spec (\textup{Sym}_A^*(\mathbb{L}_{A/K}[1])))\otimes_AK\simeq \Spec (\textup{Sym}_K^*(\mathbb{L}_{A/K}[1]\otimes_AK)))\simeq \Spec (\textup{Sym}_K^*(\mathbb{L}_{K/A}))
\]
where $\mathbb{L}_{K/A}$ is the cotangent complex of $\Spec K\to \Spec A_X=\Spec A$. the third equivalence follows from an exact triangle
\[
\mathbb{L}_{A/K}\otimes_AK\to \mathbb{L}_{K/K}\simeq 0 \to \mathbb{L}_{K/A}\to \mathbb{L}_{A/K}\otimes_AK[1].
\]
The cotangent complex $\mathbb{L}_{K/A}$ can be described as indecomposable elements.
For an argmented cofibrant commutative differential graded algebra $A\to K$,
the complex of indecomposable elements is defined to
be $I/I^2$ where $I$ is the kernel of $A\to K$.
Here by a commutative differential graded algebra $A$ we mean an ``actual''
chain complex endowed with commutative algebra structure (cf. \cite[X, 2.3.10, 2.3.11]{DAGn}), and we think of $I/I^2$ as an object in $\Mod_K$.
There is an equivalence $I/I^2[1]\simeq \mathbb{L}_{K/A}$ in $\Mod_K$ (see \cite[2.3.11]{DAGn}).

\begin{Proposition}
\label{UGEX}
There exist equivalences of derived schemes
\[
\mathsf{UG}\simeq \Spec (\textup{Sym}_K^*(\mathbb{L}_{K/A})))\simeq \Spec (\textup{Sym}_K^*(I/I^2[1])).
\]
The coarse moduli space $UG$ is isomorphic to
$\Spec (\textup{Sym}_K^*(H^1(I/I^2))$ as schemes.
\end{Proposition}

\Proof
We have already proved the first claim.
To see the second claim, for any discrete $K$-algebra $R\in \CAlg^{\textup{dis}}_K$
(i.e., $H_i(R)=0$ for $i\neq 0$), there is a natural homotopy equivalence
$\Map_{\Aff_K}(\Spec R,\Spec (\textup{Sym}_K^*\tau_{\le 0}(I/I^2[1])))\simeq \Map_{\Aff_K}(\Spec R,\Spec (\textup{Sym}_K^*(I/I^2[1])))$ where $\tau_{\le0}$
is the left adjoin truncation $\Mod_{K}\to \Mod_{K,\le0}$.
Here $\Mod_{K,\le0}$ is the full subcategory of $\Mod_K$
that consists of those objects $M$ such that $H_i(M)=0$ for $i>0$.
Moreover, the functor $\CAlg_K^{\textup{dis}}\to \mathcal{S}$ given by $R\mapsto \pi_0(\Map_{\Aff_K}(\Spec R,\Spec (\textup{Sym}_K^*\tau_{\le0}(I/I^2[1]))))$
is represented by $\Spec (\textup{Sym}_K^*H^0(\tau_{\le0}(I/I^2[1]))\simeq \Spec (\textup{Sym}_K^*(H^1(I/I^2))$.
Thus $\Spec (\textup{Sym}_K^*(H^1(I/I^2)))$ is an excellent coarse moduli space
of $\mathsf{UG}$ in the sense of \cite[Definition 7.15]{Bar}.
Therefore we deduce that $UG\simeq \Spec (\textup{Sym}_K^*(H^1(I/I^2)))$.
\QED

{\it The second proof of the second assertion in Theorem~\ref{unconditional}.}
It follows from Proposition~\ref{UGEX} and Lazard's theorem for pro-algebraic groups.
According to Proposition~\ref{UGEX} the underlying scheme of pro-algebraic group $UG$ is an (infinite dimensional)
affine space $\Spec (\textup{Sym}^*(H^1(I/I^2)))$.
Then thanks to
Lazard's theorem for pro-algebraic groups due to Chalupnik-Kowalski \cite{CK},
$UG$ is a pro-unipotent group scheme.
The connectedness of $UG$ is the same as the first proof.
\QED

\begin{Remark}
In the light of Mumford-Tate conjecture
it is desirable to have an argument that
uses Hodge realization.
\end{Remark}

\section{\'Etale Realization}

\subsection{}
In this Section, we construct an $l$-adic \'etale realization functor
from $\mathsf{DM}^\otimes$ at the level of symmetric monoidal $\infty$-categories, which equips with an action of the absolute Galois group.

 We fix some notation. Some are different from
 those used in Section 2--4 since we here use
 various coefficients. Let $R$ be a commutative ring.
We here suppose that the base scheme $S$ is the Zariski spectrum of
a perfect field $k$ and let $\Cor$ be
the $\ZZ$-linear category
of finite $S$-correspondences (cf. Section~\ref{Modelmot}).
Let $\Sh_{\textup{\'et}}(\Cor,R)$ be the abelian category of \'etale
sheaves with transfers with $R$-coefficients.
Namely, an object in $\Sh_{\textup{\'et}}(\Cor,R)$
is an additive functor $F:\Cor^{op}\to R\textup{-mod}$
such that $\textup{Sm}_{/k}^{op} \to \Cor^{op}\stackrel{F}{\to}R\textup{-mod}$ is an \'etale sheaf (cf. \cite[Lec. 6]{MVW}).
Here $k$ is a perfect field and $R\textup{-mod}$ is the abelian
category of $R$-modules.
If $\iota:\mathcal{L}(\Spec k)\to \mathcal{L}(\mathbb{G}_m)$ is a morphism
induced by
the unit morphism $\Spec k\to \mathbb{G}_m$,
we let $R(1)$ be $\Coker(\iota)[-1]$
in $\Comp(\Sh_{\textup{\'et}}(\Cor,R))$.
We denote by $R(X)\in \Comp(\Sh_{\textup{\'et}}(\Cor,R))$ the image of
$X\in \textup{Sm}_{/k}$ under $\mathcal{L}$, which we regard as a complex placed in degree zero.
Consider the 
symmetric monoidal presentable category
$\textup{Sp}_{R(1)}^\otimes(\Comp(\Sh_{\textup{\'et}}(\Cor,R)))$.
For simplicity, we shall write $\MSp^\otimes(R)$ for
$\textup{Sp}_{R(1)}^\otimes(\Comp(\Sh_{\textup{\'et}}(\Cor,R)))$,
but we often omit the superscript $\otimes$.
We equip $\MSp(R)$ with the stable
symmetric monoidal model structure described in \cite[Example 7.15]{CD1}
(but we here use the localization by \'etale hypercoverings
instead of Nisnevich hypercoverings).
We abuse notation and write $R(X)$ (resp. $R(1)$)
also for the images of $R(X)$
(resp. $R(1)$) in $\MSp(R)$ and the associated $\infty$-category.
Let $\Sh(\Cor,R)$ be the Nisnevich version of $\ESh(\Cor,R)$, that is,
it consists of Nisnevich sheaves of $R$-modules with transfers.
Similarly, let $\NSp(R)$ be the Nisnevich version of $\MSp(R)$.
If no confusion seems to arise,
we use the notation $R(X)$ and $R(1)$ also in the Nisnevich case.
Let $\mathsf{DM}^\otimes(R)$ be the symmetric monoidal stable
presentable $\infty$-category obtained from $\NSp(R)^c$
by inverting weak equivalences. When $R=\QQ$, $\mathsf{DM}^\otimes(\QQ)$ is $\mathsf{DM}^\otimes$ in the previous Sections.
Let $\mathsf{DM}_{\textup{gm}}(R)$
be the smallest stable idempotent complete subcategory which consists of $\{R(X)\otimes R(n)\}_{X\in \textup{Sm}_{/k}, n\in \ZZ}$.
Remark that $\mathsf{DM}_{\textup{gm}}(R)$
is closed under the tensor operation
$\otimes:\mathsf{DM}(R)\times \mathsf{DM}(R)\to \mathsf{DM}(R)$
and it inherits a symmetric monoidal structure.
We denote by $\textup{SmPr}_{/k}$ the category of projective smooth schemes
over $k$.
The main goal of this section is the following:

\begin{Proposition}
\label{realization}
Let $l$ be a prime number which is different from
the characteristic of $k$.
Then there exists a symmetric monoidal exact functor
\[
\mathsf{R}_{\textup{\'et},\ZZ_l}:\mathsf{DM}^\otimes_{\textup{gm}}(\ZZ)\to  \Mod_{\ZZ_l}^\otimes
\]
which has the following properties:
\begin{enumerate}
\item for any $X\in \textup{SmPr}_{/k}$ and $m\in \ZZ$, $n\in \ZZ$, there is
a natural isomorphism $H^n(\mathsf{R}_{\textup{\'et},\ZZ_l}(\ZZ(X)^{\vee}\otimes\ZZ(m))) \simeq H^{n}_{\textup{\'et}}(\overline{X},\ZZ_l(m))$, where
$\ZZ(X)^\vee$ is a dual of $\ZZ(X)$, $\overline{X}$ is $X\times_k\bar{k}$ with an algebraic closure $\bar{k}$ of $k$,
and $H^{n}_{\textup{\'et}}(\overline{X},\ZZ_l(m))$ is the $l$-adic \'etale
cohomology $\varprojlim_j H^{n}_{\textup{\'et}}(\overline{X},\mu_{l^j}^{\otimes m})$,

\item there is an action of $\Gamma=\Gal(\bar{k}/k)$ on $\mathsf{R}_{\textup{\'et},\ZZ_l}$,
that is a morphism $\Gamma\to \Aut(\mathsf{R}_{\textup{\'et},\ZZ_l})$
as objects in $\Grp(\widehat{\SSS})$,
such that it induces the action on $H^n(\mathsf{R}_{\textup{\'et},\ZZ_l}(\ZZ(X)^{\vee}\otimes\ZZ(m)))$ which coincides with the Galois action on $H^{n}_{\textup{\'et}}(\overline{X},\ZZ_l(m))$.
\end{enumerate}
Moreover, there is also a symmetric monoidal exact functor $\mathsf{R}_{\textup{\'et},\QQ_l}:\mathsf{DM}^\otimes_{\textup{gm}}(\QQ)\to  \PMod_{\QQ_l}^\otimes$ satisfying the same properties; $H^n(\mathsf{R}_{\textup{\'et},\QQ_l}(\QQ(X)^{\vee}\otimes\QQ(m)))=H^{n}_{\textup{\'et}}(\overline{X},\QQ_l(m))$
and the Galois actions coincide.
\end{Proposition}

\begin{Remark}
\label{realrem}
For any projective smooth scheme $X$ over $k$,
$\ZZ(X)$ admits a dual $\ZZ(X)^{\vee}$ (cf. \cite[11.3.4 (4)]{CD3}.
$\Aut(\mathsf{R}_{\textup{\'et},\ZZ_l})$ is the automorphism
group space of $\mathsf{R}_{\textup{\'et},\ZZ_l}$, that is,
a group object
in $\mathcal{S}$ which is
$\Aut(\mathsf{R}_{\textup{\'et},\ZZ_l})(\QQ_L)$ in Setion~\ref{realstep3}.
The essential image of $\mathsf{R}_{\textup{\'et},\QQ_l}$
is contained in $\PMod_{\QQ_l}$ since
every object in $\mathsf{DM}^\otimes_{\textup{gm}}(\QQ)$
is dualizable.

Using the machinery of mixed Weil theory (see \cite[17.2]{CD3}) one can construct 
an \'etale realization functor with $\QQ_l$-coefficients
at the level of symmetric monoidal $\infty$-categories (as we did in \cite{Tan}). New pleasant features here are (1) the realization functor comes equipped with Galois action, and (2) we can work with $\ZZ_l$-coefficients.
Our construction makes use of the rigidity theorem due to Suslin-Voevodsky
\cite{MVW}, \cite{SV}, and the derived generalization of Grothendieck's
existence theorem by Lurie \cite[XII]{DAGn}.
\end{Remark}

\subsection{}
\label{realstep1}
We will begin by constructing some symmetric monoidal Quillen functors.
Let $f:R\to R'$ be a homomorphism of commutative rings.
Then it gives rise to a symmetric monoidal colimit-preserving
functor
$\widehat{f}:\Comp(\Sh_{\textup{\'et}}(\Cor,R))\to \Comp(\Sh_{\textup{\'et}}(\Cor,R'))$ which carries $R(X)$ to $R'(X)$.
Thus we have a symmetric monoidal colimit-preserving functor
$\textup{Sp}(\widehat{f}):\MSp(R)\to \MSp(R')$ which carries
the symmetric spectrum $(E_n,\sigma_n)_{\ge0}$ (here each $\sigma_n$ denotes the structure map $E_n\otimes R(1)\to E_{n+1}$) to $(\widehat{f}(E_n))_{\ge0}$ with structure maps
$\widehat{f}(E_n)\otimes R'(1)\simeq \widehat{f}(E_n\otimes R(1))\stackrel{\widehat{f}(\sigma_n)}{\to} \widehat{f}(E_{n+1})$.
Moreover, we see the following:

\begin{Lemma}
\label{basechangeQuillen}
$\textup{Sp}(\widehat{f}):\MSp(R)\to \MSp(R')$ is a left Quillen adjoint functor.
\end{Lemma}

\Proof
The functor
$\Comp(\Sh_{\textup{\'et}}(\Cor,R))\to \Comp(\Sh_{\textup{\'et}}(\Cor,R'))$
is a left adjoint by adjoint functor theorem.
Let $\widehat{g}$ be a right adjoint of $\widehat{f}$.
Moreover, $\widehat{f}$ is a left Quillen functor
with respect to the model structure given in \cite[Example 4.12]{CD1}
since it preserves
generating cofibrations and generating trivial cofibrations,
and $R(X\times \mathbb{A}^1)\to R(X)$ induced by the projection
maps to $R'(X\times \mathbb{A}^1)\to R'(X)$; see \cite[Proposition 4.9]{CD1}.
Observe that $\textup{Sp}(\widehat{f})$
is a left Quillen adjoint.
Explicitly, its right adjoint
carries $(F_n,\tau_n)_{n\ge0}$ to $(\widehat{g}(F_n))_{n\ge 0}$
with structure maps given by th composite
$\widehat{f}(\widehat{g}(F_n)\otimes R(1))\simeq \widehat{f}(\widehat{g}(F_n))\otimes R'(1)\to F_n\otimes R'(1)\stackrel{\tau_n}{\to} F_{n+1}$.
Therefore, the right adjoint preserves
termwise $\mathbb{A}^1$-local fibrations and termwise $\mathbb{A}^1$-local
equivalences.
and thus $\MSp(R)\rightleftarrows \MSp(R')$ is a Quillen
adjunction with respect to the $\mathbb{A}^1$-local projective model structures (cf. Section~\ref{Modelmot}).
Then according to \cite[Theorem 2.2]{HoG}, to see that $\MSp(R)\rightleftarrows \MSp(R')$ is a Quillen adjunction with respect to the stable model structures,
it is enough to observe that
$\textup{Sp}(\widehat{f})(s_n^{C})$ is a stable equivalence
whenever $n\ge0$ and $C$ is either a domain or codomain of
a generating cofibration in $\Comp(\ESh(\Cor,R))$ (see \cite[7.7]{HoG} for the notation).
Note that $\textup{Sp}(\widehat{f})(s_n^{C})=s_n^{\widehat{f}(C)}$. Hence our claim is
clear.
\QED

Let $l$ be a prime number which is different from
the characteristic of $k$.
For each $n\ge 1$ the natural projection $\ZZ\to \ZZ/l^n$
induces a symmetric monoidal left Quillen adjoint functor
$\MSp(\ZZ)\to \MSp(\ZZ/l^n)$ which we denote by $p_n$.
Moreover, the projective system $\cdots \to \ZZ/l^n\to \cdots \to \ZZ/l^2\to \ZZ/l$ induces a diagram $(\clubsuit)$ of stable symmetric monoidal
model categories 
\[
\MSp(\ZZ)\to \cdots \to \MSp(\ZZ/l^n)\to \cdots \to \MSp(\ZZ/l^2)\to \MSp(\ZZ/l)
\]
in which all arrows are symmetric monoidal left Quillen functors.

Consider the Quillen adjoint pair 
\[
\Sigma^\infty :\Comp(\Sh_{\textup{\'et}}(\Cor,\ZZ/l^n))\rightleftarrows \MSp(\ZZ/l^{n}\ZZ):\Omega^\infty.
\]
We here equip $\Comp(\Sh_{\textup{\'et}}(\Cor,\ZZ/l^n))$ with
the model structure given in \cite[4.12]{CD1}
such that weak equivalences (resp. fibrations) are exactly $\mathbb{A}^1$-local
equivalences (resp. $\mathbb{A}^1$-local fibrations).
We refer to it as the $\mathbb{A}^1$-local model structure.
For ease of notation,
we write $\EM(\ZZ/l^n)$
for $\Comp(\Sh_{\textup{\'et}}(\Cor,\ZZ/l^n))$.
Observe that it is a Quillen equivalence.
According to \cite[8.19]{MVW}, tensoring with $\ZZ/l^n(1)$
is invertible in the homotopy category
of $\EM(\ZZ/l^n)$, thus we deduce from \cite[Theorem 8.1]{HoG}
that the pair $(\Sigma^\infty,\Omega^{\infty})$
is a Quillen equivalence.

Let $\ESh(\Et_{/k},\ZZ/l^n)$ be the Grothendieck abelian category
of sheaves of $\ZZ/l^n$-modules on the \'etale site
$\Et_{/k}$ which consists of \'etale morphisms of finite type
to $\Spec k$.
Let $X\in \Et_{/k}$.
Let $\ZZ/l^n[X]$ be the \'etale sheaf $\Et_{/k}^{op}\to \ZZ/l^n\textup{-mod}$
associated to
the presheaf $Y\mapsto \ZZ/l^n \cdot \Hom_k(Y,X)$
where $\ZZ/l^n\cdot\Hom_k(Y,X)$ is the free $\ZZ/l^n$-module
generated by the set of $k$-morphisms $\Hom_k(Y,X)$.
This sheaf is the restriction of
$\ZZ/l^n(X)$ in $\ESh(\Cor,\ZZ/l^n)$ to $\Et_{/k}^{op}$
(by the classical Galois theory).
The abelian category
$\ESh(\Et_{/k},\ZZ/l^n)$
has a set of generators $\{\ZZ/l^n[X]\}_{X\in \Et_{/k}}$.
Recall that the symmetric monoidal
category $\ESh(\Et_{/k},\ZZ/l^n)$ is equivalent to
the symmetric monoidal category
$\ZZ/l^n\textup{-mod}_{\Gamma}$
of discrete $\ZZ/l^n$-modules with action of $\Gamma$,
where $\Gamma$ is the absolute Galois group $\Gal(\bar{k}/k)$
($\bar{k}$ is an algebraic closure of $k$).
In the symmetric monoidal category $\ESh(\Et_{/k},\ZZ/l^n)$,
$\ZZ/l^n[X]\otimes \ZZ/l^n[Y]=\ZZ/l^n[X\times_kY]$
and the commutative constraint is determined by the flip
$X\times_kY\simeq Y\times_kX$.
For a $k$-field $L$ the equivalences $\ESh(\Et_{/k},\ZZ/l^n)
\simeq \ZZ/l^n\textup{-mod}_{\Gamma}$
send $\ZZ/l^n[\Spec L]$ to $\ZZ/l^n\cdot\Hom_{k}(\Spec \bar{k},\Spec L)$ with action of $\Gamma$ induced by composition.
Here $\ZZ/l^n\cdot\Hom_{k}(\Spec \bar{k},\Spec L)$
is a free $\ZZ/l^n$-module generated by the set of morphisms from $\Spec \bar{k}$ to $\Spec L$ over $k$.
It carries $\ZZ/l^n[\Spec L]\otimes \ZZ/l^n[\Spec L']=\ZZ/l^n[\Spec L\times_k\Spec L']$
to
\[
\ZZ/l^n\cdot\Hom_{k}(\Spec \bar{k},\Spec L)\otimes_{\ZZ/l^n}\ZZ/l^n\cdot\Hom_{k}(\Spec \bar{k},\Spec L')
\]
equipped with the action of $\Gamma$ (by the tensor operation).
Since $\{\ZZ/l^n[X]\}_{X\in \Et_{/k}}$ is the set of compact generators,
and the tensor operations preserves colimits in each variable, thus
we see that the equivalences
$\ESh(\Et_{/k},\ZZ/l^n)
\simeq \ZZ/l^n\textup{-mod}_{\Gamma}$ are extended to symmetric monoidal
equivalences.

Let $\Et_{/\bar{k}}$ be the \'etale site over $\Spec \bar{k}$.
The geometric point $q:\Spec \bar{k}\to \Spec k$ given by one inclusion
$k\subset \bar{k}$ determines the exact pullback functor
\[
q_n^*:\ESh(\Et_{/k},\ZZ/l^n)\to \ESh(\Et_{/\bar{k}},\ZZ/l^n).
\]
This sends $\ZZ/l^n[X]$ to the sheaf $\ZZ/l^n[X\times_k\Spec \bar{k}]$,
that is given by $Y\mapsto \ZZ/l^n\cdot \Hom_{\Spec {\bar{k}}}(Y,X\times_k\Spec \bar{k})$. Note that there is a symmetric monoidal equivalence
$\ESh(\Et_{/\bar{k}},\ZZ/l^n)\simeq \ZZ/l^n\textup{-mod}$
which carries $F$ to $F(\Spec \bar{k})$. If one
identifies $\ESh(\Et_{/k},\ZZ/l^n)$
with $\ZZ/l^n\textup{-mod}_{\Gamma}$, then $q_n^*$ is equivalent to
the forgetful functor $\ZZ/l^n\textup{-mod}_{\Gamma}\to \ZZ/l^n\textup{-mod}$ as symmetric monoidal functors.
We equip $\Comp(\ESh(\Et_{/k},\ZZ/l^n))$ and
$\Comp(\ESh(\Et_{/\bar{k}},\ZZ/l^n))$ with the symmetric monoidal
model structures
given in \cite[Proposition 3.2, Example 2.3]{CD1}, in which weak equivalences are exactly
quasi-isomorphisms.

Let
$v_n:\ESh(\Et_{/k},\ZZ/l^n) \to \ESh(\Cor,\ZZ/l^n)$
be a left adjoint exact functor, given in \cite[6.7--6.11]{MVW}, which carries $F:\Et_{/k}^{op}\to \ZZ/l^n\textup{-mod}$ to a unique sheaf with transfers
$\overline{F}:\Cor^{op}\to \ZZ/l^n\textup{-mod}$
such that the composite
$\Et_{/k}^{op}\to \Cor^{op}\stackrel{\overline{F}}{\to} \ZZ/l^n\textup{-mod}$ is $F$. The right adjoint is determined by the composition with $\Et_{/k}^{op}\to \Cor^{op}$.
For any $X\in \Et_{/k}$, $\ZZ/l^n[X]\in \ESh(\Et_{/k},\ZZ/l^n)$
maps to $\ZZ/l^n(X)$ as an object $\ESh(\Cor,\ZZ/l^n)$.
Let $\ESh^{\textup{rep}}(\Et_{/k},\ZZ/l^n)$
be the full subcategory of $\ESh(\Et_{/k},\ZZ/l^n)$ spanned by
representable objects $\{\ZZ/l^n[X]\}_{X\in \Et_{/k}}$.
For $S\in \ESh(\Et_{/k},\ZZ/l^n)$, we have
$S\simeq \colim_{\ZZ/l^n[X]\to S}\ZZ/l^n[X]$
where $\ZZ/l^n[X]\to S$ run over the small overcategory
$\ESh^{\textup{rep}}(\Et_{/k},\ZZ/l^n)_{/S}$.
Since $v_n$ preserves small colimits, there are natural equivalences
\begin{eqnarray*}
v_n(S)\otimes v_n(T) &\simeq& (\colim'_{\ZZ/l^n[X]\to S}v_n(\ZZ/l^n[X]))\otimes (\colim'_{\ZZ/l^n[Y]\to T}v_n(\ZZ/l^n[Y])) \\
&\simeq& \colim'_{\ZZ/l^n[X]\to S}(\colim'_{\ZZ/l^n[Y]\to T}(\ZZ/l^n(X)\otimes\ZZ/l^n(Y))) \\
&\simeq& v_n(\colim_{\ZZ/l^n[X]\to S}(\colim_{\ZZ/l^n[Y]\to T}(\ZZ/l^n[X]\otimes\ZZ/l^n[Y]))) \\
&\simeq& v_n((\colim_{\ZZ/l^n[X]\to S}\ZZ/l^n[X])\otimes (\colim_{\ZZ/l^n[Y]\to T}\ZZ/l^n[Y])) \\
&\simeq& v_n(S\otimes T)
\end{eqnarray*}
where $\colim'$ stands for the colimit in $\ESh(\Cor,\ZZ/l^n)$.
Similarly, the commutative constraint $i:S\otimes T\simeq T\otimes S$
commutes with $v_n(S)\otimes v_n(T)\simeq v_n(T)\otimes v_n(S)$. Moreover, $v_n(\ZZ/l^n[\Spec k])=\ZZ/l^n(\Spec k)$. Hence we easily see that $v_n$ is
(extended to) a symmetric monoidal functor.
Thus we have an adjoint pair
\[
v_n:\Comp(\ESh(\Et_{/k},\ZZ/l^n)) \rightleftarrows \Comp(\ESh(\Cor,\ZZ/l^n)):res\]
where the left adjoint $v_n$ is symmetric monoidal, and $res$
is induced by the composition with the natural functor
$\Et_{/k}^{op}\to \Cor^{op}$.

\begin{Lemma}
\label{leftQuillen}
We abuse notation and write
\[
q_n^*:\Comp(\ESh(\Et_{/k},\ZZ/l^n))\to 
\Comp(\ESh(\Et_{/\bar{k}},\ZZ/l^n))
\]
and
\[
v_n:\Comp(\ESh(\Et_{/k},\ZZ/l^n))\to 
\Comp(\ESh(\Cor,\ZZ/l^n))
\]
for symmetric monoidal functors induced by $q_n^*:\ESh(\Et_{/k},\ZZ/l^n\to \ESh(\Et_{/\bar{k}},\ZZ/l^n)$ and $v_n:\ESh(\Et_{/k},\ZZ/l^n)\to \ESh(\Cor,\ZZ/l^n)$ respectively.
Then both $q_n^*$ and $v_n$ are left Quillen adjoint functors.
\end{Lemma}

\Proof
According to the definitions of generating (trivial) cofibrations and
\cite[Theorem 2.14]{CD1}, we see that $q_n^*$ is a left Quillen adjoint.
Next we will show that $v_n$ is a left Quillen adjoint.
We equip $\Comp(\ESh(\Cor,\ZZ/l^n))$ with the symmetric monoidal
model structure given in \cite[Example 2.4]{CD1}, in which weak equivalences
are exactly quasi-isomorphisms. Then the $\mathbb{A}^1$-local
model structure is its left Bousfield localization with respect to
$\{\ZZ/l^n(X\times\mathbb{A}^1)\to \ZZ/l^n(X)\}_{X\in \textup{Sm}_{/k}}$.
Thus it is enough to prove that $v_n$ is left Quillen when
$\Comp(\ESh(\Cor,\ZZ/l^n))$ is endowed with the
model structure given in \cite[Example 2.4]{CD1}.
To this end, note that $v_n$ preserves quasi-isomorphisms.
It remains to show that $v_n$ preserves generating cofibrations.
But it is clear from the definitions of generating cofibrations
(see \cite[Definition 2.2]{CD1})
and the fact that $v_n(\ZZ/l^n[X])=\ZZ/l^n(X)$.
\QED

We obtain the diagram of symmetric monoidal left Quillen functors
\[
\xymatrix{
\MSp(\ZZ) \ar[r] & \MSp(\ZZ/l^n) &   &      \\
\Comp(\ESh(\Cor,\ZZ/l^n))\ar[ur]^{\Sigma^\infty} & \Comp(\ESh(\textup{\'Et}_{/k},\ZZ/l^n)) \ar[l]^{v_n} \ar[r]^{q_n^*} & \Comp(\ESh(\textup{\'Et}_{/\bar{k}},\ZZ/l^n)) & 
}
\]
where $(\Sigma^\infty,\Omega^\infty)$ is a left Quillen equivalence.

\subsection{}
\label{realstep2}
Next we consider an action of the absolute Galois group on functors.
There is a commutative diagram $(\heartsuit)$
of symmetric monoidal model categories
\[
\xymatrix{
\cdots \ar[r]  & \Comp(\ESh(\Et_{/k},\ZZ/l^2)) \ar[d]^{q_2^*} \ar[r] & \Comp(\ESh(\Et_{/k},\ZZ/l)) \ar[d]^{q_1^*} \\
\cdots \ar[r] & \Comp(\ESh(\Et_{/\bar{k}},\ZZ/l^2))\ar[r] & \Comp(\ESh(\Et_{/\bar{k}},\ZZ/l))
}
\]
in which all arrows are symmetric monoidal left Quillen functors
(as in the proof of Lemma~\ref{basechangeQuillen},
for any $n\ge1$, $\Comp(\ESh(\Et_{/k},\ZZ/l^{n+1}\ZZ))\to \Comp(\ESh(\Et_{/k},\ZZ/l^{n}\ZZ))$
is a left Quillen functor since it
preserves small colimits and generating (trivial) cofibrations, see \cite[2.4]{CD1}).
For each $n\ge1$, we have the symmetric monoidal functor
\[
q_n^*:\Comp(\ESh(\Et_{/k},\ZZ/l^n))\simeq \ZZ/l^n\textup{-mod}_{\Gamma}\to
\ZZ/l^n\textup{-mod}\simeq \Comp(\ESh(\Et_{/\bar{k}},\ZZ/l^n))
\]
where the middle functor is the forgetful functor.
Then $\Gamma$ acts on the forgetful functor, i.e., $q_n^*$.
That is, if $\Aut(q_n^*)$ denotes the group of symmetric monoidal
natural equivalences of $q_n^*$, then we have the homomorphism
$\Gamma\to \Aut(q_n^*)$ of groups which carries $g\in \Gamma$
to the symmetric monoidal natural equivalence given by morphisms
$g:q_n^*(C)\to q_n^*(C)$ induced by the action of $\Gamma$ on $C\in \ZZ/l^n\textup{-mod}_{\Gamma}$.
This action commutes with the diagram in the following sense:
For any pair $m,n\in \NN$
with $m\ge n$, the action on $q_n^*$ and the vertical composition with
$\Comp(\ESh(\Et_{/k},\ZZ/l^{m}\ZZ)) \to \Comp(\ESh(\Et_{/k},\ZZ/l^n))$
determines an action of $\Gamma$ on
$\Comp(\ESh(\Et_{/k},\ZZ/l^{m}\ZZ))\to \Comp(\ESh(\Et_{/\bar{k}},\ZZ/l^n))$. On the other hand,
the action on $q_{m}^*$ and the vertical composition with
$\Comp(\ESh(\Et_{/\bar{k}},\ZZ/l^{m}\ZZ)) \to \Comp(\ESh(\Et_{/\bar{k}},\ZZ/l^n))$
also determines another action of $\Gamma$.
Then two actions coincide for any $m,n$ with $m\ge n$.

Consider the category $I$ of the form
\[
\cdots \to n\to n-1\to \cdots \to 2\to 1.
\]
Namely, objects of $I$ are natural numbers, and
the homset $\Hom_I(n,m)$ consists of one point if $n\ge m$
and $\Hom_I(n,m)$ is the empty if otherwise. We abuse notation and
write $I$ also for the nerve of $I$.
Let $\textup{WCat}_\infty$ be the $\infty$-category which consists
of pairs $(\mathcal{C},W)$ where $\mathcal{C}$ is an $\infty$-category
and $W$ is a subset of edges of $\mathcal{C}$, called a system, which are stable under
homotopy, composition and contains all weak equivalences.
Here $\mathcal{C}$ belongs to
an enlarged universe which contains model categories we treat.
The mapping space $\Map_{\textup{WCat}_\infty}((\mathcal{C},W),(\mathcal{C}',W'))$ is equivalent to the summands spanned by $f:\mathcal{C}\to\mathcal{C}'$ such that $f(W)\subset W'$; see \cite[4.1.3.1]{HA}.
Moreover, we equip $\textup{WCat}_\infty$ with the Cartesian monoidal
structure and write $\CAlg(\textup{WCat}_\infty)$
for the $\infty$-category of commutative algebra objects
with respect to this monoidal structure.
Then the diagram $(\heartsuit)$
induces $\alpha:\Delta^1\times I\to \CAlg(\textup{WCat}_\infty)$
such that the restriction $\alpha_0:\{0\}\times I\to \CAlg(\textup{WCat}_\infty)$
is given by the sequence of symmetric monoidal full subcategories spanned by cofibrant objects:
\[
\cdots \to \Comp(\ESh(\Et_{/k},\ZZ/l^2))^c \to \Comp(\ESh(\Et_{/k},\ZZ/l))^c,
\]
and the restriction $\alpha_1:\{1\}\times I\to \CAlg(\textup{WCat}_\infty)$
is given by 
\[
\cdots \to \Comp(\ESh(\Et_{/\bar{k}},\ZZ/l^2))^c\to \Comp(\ESh(\Et_{/\bar{k}},\ZZ/l))^c.
\]
Let us denote by $D_k$
and $D_{\bar{k}}$ the objects in $\Fun(I,\CAlg(\textup{WCat}_{\infty}))$
corresponding to $\alpha_0$ and $\alpha_1$
respectively.
The mapping space $\Map(D_{k},D_{\bar{k}})$
from $D_k$ to $D_{\bar{k}}$ in $\Fun(I,\CAlg(\textup{WCat}_\infty))$
is described by the Kan complex
\[
\{D_{k}\}\times_{\Fun(\{0\},\Fun(I,\CAlg(\textup{WCat}_\infty)))}\Fun(\Delta^1,\Fun(I,\CAlg(\textup{WCat}_\infty)))\times_{\Fun(\{1\},\Fun(I,\CAlg(\textup{WCat}_\infty)))}\{D_{\bar{k}}\}.
\]
Clearly, this mapping space is $1$-truncated since $D_k$ and $D_{\bar{k}}$
are sequences of symmetric monoidal $1$-categories equipped with systems.
Let $\textup{FGal}(\bar{k}/k)$ be the nerve of the category 
which consists of one object $\{*\}$ and $\Hom_{\textup{FGal}(\bar{k}/k)}(*,*)
=\Gal(\bar{k}/k)=\Gamma$ equipped with the composition determined by
the multiplication.
The action of $\Gamma$ on $\{q_n^*\}_{n\ge 1}$
induces the action
of $\Gamma$ on the morphism $f:D_k\to D_{\bar{k}}$
corresponding to $\alpha$
in $\Fun(I, \textup{WCat}_{\infty})$, which is described by
a functor $\textup{FGal}(\bar{k}/k)\to \Map(D_{k},D_{\bar{k}})$
sending $\ast$ to $f\in \Map(D_{k},D_{\bar{k}})$.

Now we will construct symmetric monoidal $\infty$-categories
by inverting weak equivalences in symmetric monoidal model categories.
There is a natural fully faithful functor
$\Cat\to \textup{WCat}_{\infty}$
which carries an $\infty$-category $\mathcal{C}$
to $(\mathcal{C},W_{eq}(\mathcal{C}))$ where $W_{eq}(\mathcal{C})$
is the collection of edges of $\mathcal{C}$.
It also
induces a fully faithful functor
$\CAlg(\Cat)\to \CAlg(\textup{WCat}_\infty)$,
where $\Cat$ is endowed with the Cartesian monoidal structure;
see \cite[4.1.3]{HA}.
According to \cite[4.1.3.4]{HA} this functor admits a
left adjoint $L:\CAlg(\textup{WCat}_\infty) \to \CAlg(\Cat)$.
By composition with this adjoint pair, we have an localization adjoint pair
(see \cite[5.2.7.2]{HTT})
\[
L^I:\Fun(I,\CAlg(\textup{WCat}_\infty)) \rightleftarrows \Fun(I,\CAlg(\Cat))
\]
by \cite[5.2.7.4]{HTT}.
Let $\alpha':\Delta^1\times I\to \CAlg(\textup{WCat}_\infty)\stackrel{L}{\to} \CAlg(\Cat)$
be the composite, and let $\alpha_0'$ and $\alpha'_1$ be the restrictions
to $\{0\}\times I$ and $\{1\}\times I$ respectively.
Let $D'_k$ and $D_{\bar{k}}'$ be the objects in $\Fun(I,\CAlg(\Cat))$
corresponding to $\alpha_0'$ and $\alpha_1'$ respectively.
The functor $\alpha'$ is informally depicted as
\[
\xymatrix{
\cdots \ar[r]  & \mathcal{D}^{\otimes}(\ESh(\Et_{/k},\ZZ/l^2)) \ar[r] \ar[d] & \mathcal{D}^{\otimes}(\ESh(\Et_{/k},\ZZ/l)) \ar[d] \\
\cdots \ar[r] & \mathcal{D}^{\otimes}(\ZZ/l^2) \ar[r] & \mathcal{D}^{\otimes}(\ZZ/l)
}
\]
where $\mathcal{D}^{\otimes}(\ESh(\Et_{/k},\ZZ/l^n))$ and $\mathcal{D}^{\otimes}(\ZZ/l^n))$ are symmetric monoidal stable presentable $\infty$-categories obtained from $\Comp^{\otimes}(\ESh(\Et_{/k},\ZZ/l^n))$
and $\Comp^{\otimes}(\ESh(\Et_{/\bar{k}},\ZZ/l^n))$ respectively, by inverting weak equivalences (we often omit the subscript $\otimes$).
The functor $L^I$ induces
\[
L^{\Delta^1\times I}:\Fun(\Delta^1,\Fun(I,\CAlg(\textup{WCat}_\infty)))\to \Fun(\Delta^1,\Fun(I,\CAlg(\textup{Cat}_\infty))).
\]
By the description of $\Map(D_k,D_{\bar{k}})$, $L^{\Delta^1\times I}$ induces
$\Map(D_k,D_{\bar{k}}) \to \Map(D_k',D_{\bar{k}}')$,
where $\Map(D_k',D_{\bar{k}}')$
is given by
\[
\{D'_{k}\}\times_{\Fun(\{0\},\Fun(I,\CAlg(\textup{Cat}_\infty)))}\Fun(\Delta^1,\Fun(I,\CAlg(\textup{Cat}_\infty)))\times_{\Fun(\{1\},\Fun(I,\CAlg(\textup{Cat}_\infty)))}\{D'_{\bar{k}}\}.
\]
By composition we have $\mathbf{t}:\textup{FGal}(\bar{k}/k)\to \Map(D_k,D_{\bar{k}}) \to \Map(D_k',D_{\bar{k}}')$ carrying $*$ to $f'$ which is the image of $f$ in
$\Map(D_k',D_{\bar{k}}')$.

Let $\widehat{\mathcal{D}}^\otimes(\ESh(\Et_{/k},\ZZ_l))$
(resp. $\widehat{\mathcal{D}}^\otimes(\ZZ_l)$)
be a symmetric monoidal stable presentable
$\infty$-category which is defined to be the limit of $\alpha_0':I\simeq \{0\}\times I\to \CAlg(\Cat)$
(resp. $\alpha_1':I\simeq \{1\}\times I\to \CAlg(\Cat)$).
Then $\alpha'$ determines a symmetric monoidal colimit-preserving functor
$\widehat{\mathcal{D}}^\otimes((\ESh(\Et_{/k},\ZZ_l))\to \widehat{\mathcal{D}}^\otimes(\ZZ_l)$.

Let $I\to \CAlg(\textup{WCat}_\infty)$ be the functor
corresponding to $\cdots \to \MSp^\otimes(\ZZ/l^2)^c\to \MSp^\otimes(\ZZ/l)^c$.
Composing with $L^I$ we have $I\to \CAlg(\Cat)$ which we described as
$\cdots \to \mathsf{DM}_{\textup{\'et}}^\otimes(\ZZ/l^2)\to \mathsf{DM}_{\textup{\'et}}^\otimes(\ZZ/l)$.
Let $\widehat{\mathsf{DM}}^\otimes_{\textup{\'et}}(\ZZ_l)$ be its limit.
Similarly, from
$\cdots \to \Comp^{\otimes}(\ESh(\Et_{/k},\ZZ/l^2))^c\to \Comp^{\otimes}(\ESh(\Et_{/k},\ZZ/l))^c$
we have a sequence $\cdots \to (\mathsf{DM}_{\textup{\'et}}^{\textup{eff}})^\otimes(\ZZ/l^2)\to (\mathsf{DM}_{\textup{\'et}}^{\textup{eff}})^\otimes(\ZZ/l)$. Let $(\widehat{\mathsf{DM}}^{\textup{eff}}_{\textup{\'et}})^\otimes(\ZZ_l)$ be its limit.
Let $\mathsf{DM}^\otimes_{\textup{\'et}}(\ZZ)$ be the symmetric monoidal
stable presentable
$\infty$-category obtained from $\MSp^\otimes(Z)^c$ by 
inverting weak equivalences.
Then the diagram $(\clubsuit)$ induces
$\mathsf{DM}^\otimes_{\textup{\'et}}(\ZZ)\to \widehat{\mathsf{DM}}^\otimes_{\textup{\'et}}(\ZZ_l)$.
Consider the symmetric monoidal left Quillen functors
\[
\Comp^{\otimes}(\ESh(\Et_{/k},\ZZ/l^n))\stackrel{v_n}{\to} \Comp^\otimes(\ESh(\Cor,\ZZ/l^n))\stackrel{\Sigma^\infty}{\to} \MSp^\otimes(\ZZ/l^n),
\]
we obtain
$\widehat{\mathcal{D}}^\otimes((\ESh(\Et_{/k},\ZZ_l)) \stackrel{\widehat{v}}{\to} (\widehat{\mathsf{DM}}^{\textup{eff}}_{\textup{\'et}})^{\otimes}(\ZZ_l)\simeq \widehat{\mathsf{DM}}_{\textup{\'et}}^\otimes(\ZZ_l)$
where the right equivalence follows from the Quillen equivalences
$\Comp^\otimes(\ESh(\Cor,\ZZ/l^n))\stackrel{\Sigma^\infty}{\to} \MSp^\otimes(\ZZ/l^n)$. Moreover, the rigidity theorem due to Suslin-Voevodsky
\cite{SV}, \cite[9.35, 7.20]{MVW}
implies:

\begin{Lemma}
\label{suslin}
The symmetric monoidal functor $\widehat{v}$ is an equivalence.
\end{Lemma}

\Proof
Since $v_n$ determines a symmetric monoidal
exact functor between symmetric monoidal stable $\infty$-categories,
we are reduced to proving that
$\mathcal{D}((\ESh(\Et_{/k},\ZZ/l^n)) \stackrel{v_n}{\to} \mathsf{DM}^{\textup{eff}}_{\textup{\'et}}(\ZZ/l^n)$ induces
an equivalence of their homotopy categories (cf. \cite[Lemma 5.8]{Tan}).
Thus our claim is a consequence of the rigidity theorem \cite[9.35, 7.20]{MVW}.
\QED

 We define some stable subcategories.
 Let $\mathsf{DM}^{\textup{eff}}_{\textup{\'etgm}}(\ZZ/l^n)$
be the smallest stable idempotent complete subcategory of $\mathsf{DM}^{\textup{eff}}_{\textup{\'et}}(\ZZ/l^n)$ which consists of
$\{\ZZ/l^n(X)\}_{X\in \textup{Sm}_{/k}}$. Let
$\widehat{\mathsf{DM}}^{\textup{eff}}_{\textup{\'etgm}}(\ZZ_l)$
be the limit $\lim_n\mathsf{DM}^{\textup{eff}}_{\textup{\'etgm}}(\ZZ/l^n)$.
We define $\widehat{\mathcal{D}}_{\textup{gm}}(\ESh(\Et_{/k},\ZZ_l))$
to be the stable subcategory of $\widehat{\mathcal{D}}(\ESh(\Et_{/k},\ZZ_l))$
that corresponds to $\mathsf{DM}^{\textup{eff}}_{\textup{\'etgm}}(\ZZ_l)$
through the equivalence $\widehat{v}$. Both categories naturally
inherit symmetric monoidal structures.

\subsection{}
\label{realstep3}
We will construct realization functors.
Let $\Comp^\otimes(\ZZ_l\textup{-mod})$ be the symmetric 
monoidal category of chain complexes of $\ZZ_l$-modules.
We equip $\Comp^\otimes(\ZZ_l\textup{-mod})$
with the (symmetric monoidal)
projective model structure (see e.g. \cite[2.3.11]{Ho}
or \cite[8.1.2.8, 8.1.4.3]{HA}), in which
weak equivalences (resp. fibrations) are exactly quasi-isomorphisms
(resp. termwise surjective maps).
Comparing the set of generating cofibrations (see \cite[2.3.3]{Ho}),
we see that the model structures on $\Comp(\ESh(\Et_{/\bar{k}},\ZZ/l^n))\simeq \Comp(\ZZ/l^n\textup{-mod})$ is the projective model
structure.
Let $\mathcal{D}^\otimes(\ZZ_l)$
and $\mathcal{D}^\otimes(\ZZ/l^n)$ be the symmetric monoidal $\infty$-categories obtained from $\Comp(\ZZ_l\textup{-mod})^c$ and $\Comp(\ZZ/l^n\textup{-mod})^c$
by inverting weak equivalences.
According to \cite[8.1.2.13]{HA}, there are natural equivalences
$\Mod_{\ZZ_l}^\otimes\simeq \mathcal{D}^\otimes(\ZZ_l)$
and $\Mod_{\ZZ/l^n}^\otimes\simeq \mathcal{D}^\otimes(\ZZ/l^n)$
of symmetric monoidal $\infty$-categories.
The base change functor $(-)\otimes_{\ZZ_l}\ZZ/l^n:\Mod^\otimes_{\ZZ_l}\to \Mod_{\ZZ/l^n}^\otimes$ gives rise to 
$\mathcal{D}^\otimes(\ZZ_l)\simeq \Mod_{\ZZ_l}^\otimes\to \lim_n\Mod_{\ZZ/l^n}^\otimes\simeq \widehat{\mathcal{D}}^\otimes(\ZZ_l)$.
Let $R$ be a noetherian commutative ring.
Let $P$ be an object in $\mathcal{D}^\otimes(R)\simeq \Mod_{HR}^\otimes$.
We say that $P$ is almost perfect if $H^n(P)=0$ for $n>>0$ and $H^i(P)$
is a finitely presented (generated) $R$-module for each $i\in\ZZ$;
see \cite[8.2.5.10, 8.2.5.11, 8.2.5.17]{HA}, \cite[VIII, 2.7.20]{DAGn}.
Let $\AMod^\otimes_{HR}$ denote the
stable subcategory of $\Mod_{HR}$
spanned by almost perfect objects.
We remark that there is a sequence of fully faithful embeddings
$\PMod_{HR}\subset \AMod_{HR}\subset \Mod_{HR}$.
We easily see that
every almost perfect object $K$ in $\mathcal{D}^\otimes(\ZZ/l^n)\simeq \Mod_{\ZZ/l^n}$ (regarded as a chain complex)
has a quasi-isomorphism $P\to K$ such that $P$
is a right bounded complex of free $\ZZ/l^n$-modules of finite rank.
Therefore
$\AMod_{\ZZ/l^n}\subset \Mod_{\ZZ/l^n}$ is closed under
tensor operation.  Similarly, $\AMod_{\ZZ_l}$ is closed under
tensor operation. Thus these come equip with symmetric monoidal structures.
Let $\lim_n\AMod_{\ZZ/l^n}$ be the full subcategory
of $\lim_n\Mod_{\ZZ/l^n}\simeq \widehat{\mathcal{D}}(\ZZ_l)$
spanned by compatible systems $\{ C(n)\in \Mod_{\ZZ/l^n}\}_{n\le1}$
such that each $C(n)$ is almost perfect.
Then thanks to the derived version of Grothendieck existence
theorem \cite[XII, 5.3.2, 5.1.17]{DAGn}, 
$\Mod_{\ZZ_l} \to \lim_n\Mod_{\ZZ/l^n}$
induces an equivalence of symmetric monoidal $\infty$-categories
$\widehat{\pi}:\AMod_{\ZZ_l} \simeq  \lim_n\AMod_{\ZZ/l^n}$.
Let $\mathsf{DM}_{\textup{\'etgm}}(\ZZ)$ be the smallest
stable idempotent complete subcategory of $\mathsf{DM}_{\textup{\'et}}(\ZZ)$
which consists of $\{\ZZ(X)\otimes \ZZ(n)\}_{X\in \textup{Sm}_{/k} \atop n\in \ZZ}$.
By Corollary~\ref{suslinho},
combined with the finiteness of \'etale cohomology and cohomological dimension \cite{SGA4.5},
we see that 
the essential image
of $\mathsf{DM}_{\textup{\'etgm}}(\ZZ)\to \widehat{\mathcal{D}}(\ZZ_l)\simeq \lim_n\Mod_{\ZZ/l^n}$ is contained in $\lim_n\AMod_{\ZZ/l^n}$.
Hence combining with Lemma~\ref{suslin}
we obtain a symmetric monoidal exact functor as the composite
\begin{eqnarray*}
\mathcal{R}_{\textup{\'et},\ZZ_l}:\mathsf{DM}_{\textup{\'etgm}}^\otimes(\ZZ) &\to&  (\widehat{\mathsf{DM}}_{\textup{\'et}})^{\otimes}(\ZZ_l) \stackrel{(\Sigma^{\infty})^{-1}}{\simeq} (\widehat{\mathsf{DM}}_{\textup{\'et}}^{\textup{eff}})^{\otimes}(\ZZ_l) \supset (\widehat{\mathsf{DM}}^{\textup{eff}}_{\textup{\'etgm}})^{\otimes}(\ZZ_l) \\
 &\stackrel{\widehat{v}^{-1}}{\simeq}& \widehat{\mathcal{D}}^\otimes_{\textup{gm}}((\ESh(\Et_{/k},\ZZ_l)) \to (\lim_n\AMod_{\ZZ/l^n})^\otimes \stackrel{\widehat{\pi}}{\simeq} \AMod_{\ZZ_l}^\otimes\subset \mathcal{D}^\otimes(\ZZ_l)
\end{eqnarray*}
where the essential image of the first line belongs to $(\widehat{\mathsf{DM}}^{\textup{eff}}_{\textup{\'etgm}})(\ZZ_l)$.
\'Etale sheafification induces
a symmetric monoidal exact functor $\Sh(\Cor,\ZZ)\to \ESh(\Cor,\ZZ)$
giving rising to a symmetric monoidal left adjoint functor
$\NSp(\ZZ) \to \MSp(\ZZ)$. The right adjoint of $\NSp(\ZZ)\to \MSp(\ZZ)$ is the forgetful functor $\MSp(\ZZ) \to \NSp(\ZZ)$.
As in the case of $\MSp(\ZZ)$ we equip $\NSp(\ZZ)$ with a stable model structure \cite[Example 7.15]{CD1}.
Then repeating the argument of Lemma~\ref{basechangeQuillen} we see that
$\NSp(\ZZ)\to \MSp(\ZZ)$ is a symmetric monoidal left Quillen adjoint functor.
Let $\mathsf{DM}^\otimes(\ZZ)$ be the symmetric monoidal stable presentable
$\infty$-category obtained from $\NSp(\ZZ)^c$ by inverting weak equivalences.
Then $\NSp(\ZZ)\to \MSp(\ZZ)$ determines a symmetric monoidal colimit-preserving functor
$\eta:\mathsf{DM}^\otimes(\ZZ) \to \mathsf{DM}_{\textup{\'et}}^\otimes(\ZZ)$.
Consider the composite of symmetric monoidal exact functors
\[
\mathsf{R}_{\textup{\'et},\ZZ_l}:\mathsf{DM}^\otimes_{\textup{gm}}(\ZZ)\stackrel{\eta}{\longrightarrow} \mathsf{DM}^\otimes_{\textup{\'etgm}}(\ZZ)\stackrel{\mathcal{R}_{\textup{\'et},\ZZ_l}}{\longrightarrow} \AMod_{\ZZ_l}^\otimes.
\]
We shall refer to $\mathsf{R}_{\textup{\'et},\ZZ_l}$ as the $l$-adic
\'etale realization functor.
Let $\mathsf{DM}^\otimes_{\textup{gm}}(\QQ)$ be
the $\QQ$-coefficient version of
$\mathsf{DM}^\otimes_{\textup{gm}}(\ZZ)$.
By Lemma~\ref{rationaldescent} below,
the composite
$\mathsf{R}_{\textup{\'et},\ZZ_l}\otimes \QQ_l:\mathsf{DM}^\otimes_{\textup{gm}}(\ZZ)\stackrel{\mathsf{R}_{\textup{\'et},\ZZ_l}}{\to}\AMod^\otimes_{\ZZ_l}\stackrel{(-)\otimes_{\ZZ_l}\QQ_l}{\to} \AMod_{\QQ_l}^\otimes$
induces a symmetric monoidal exact functor
$\mathsf{R}_{\textup{\'et},\QQ_l}:\mathsf{DM}^\otimes_{\textup{gm}}(\QQ)\to \AMod^\otimes_{\QQ_l}$, uniquely up to a contractible space of choice, such that $\mathsf{DM}^\otimes_{\textup{gm}}(\ZZ)\to \mathsf{DM}^\otimes_{\textup{gm}}(\QQ)\stackrel{\mathsf{R}_{\textup{\'et},\QQ_l}}{\to} \AMod^\otimes_{\QQ_l}$ is $\mathsf{R}_{\textup{\'et},\ZZ_l}\otimes \QQ_l$.

Next we will define Galois actions on $\mathsf{R}_{\textup{\'et},\ZZ_l}$
and $\mathsf{R}_{\textup{\'et},\QQ_l}$.
Since $\widehat{\mathcal{D}}^\otimes(\ESh(\textup{\'Et}_{/k},\ZZ_l))$
and $\widehat{\mathcal{D}}^\otimes(\ZZ_l)$ are limits of $D'_{k}$
and $D_{\bar{k}}'$ respectively, thus by taking
their limits we have a natural map of mapping
spaces
$\Map(D_k',D_{\bar{k}}')\to \Map_{\CAlg(\Cat)}(\widehat{\mathcal{D}}^\otimes(\ESh(\textup{\'Et}_{/k},\ZZ_l)),\widehat{\mathcal{D}}^\otimes(\ZZ_l))$.
The composition with
$\mathbf{t}:\textup{FGal}(\bar{k}/k)\to \Map(D_k',D_{\bar{k}}')$ constructed
in~\ref{realstep2}
gives rise to
\[
\textup{FGal}(\bar{k}/k)
\to \Map_{\CAlg(\Cat)}(\widehat{\mathcal{D}}^\otimes(\ESh(\textup{\'Et}_{/k},\ZZ_l),\widehat{\mathcal{D}}^\otimes(\ZZ_l)).
\]
If $F:\widehat{\mathcal{D}}^\otimes(\ESh(\textup{\'Et}_{/k},\ZZ_l))\to\widehat{\mathcal{D}}^\otimes(\ZZ_l)$ denotes the symmetric monoidal
functor determined by the limit of $f':D'_k\to D_{\bar{k}}'$,
the based loop space induces $\Omega_*\textup{FGal}(\bar{k}/k)=\Gal(\bar{k}/k)\to \Aut(F):=\Omega_*\Map_{\CAlg(\Cat)}(\widehat{\mathcal{D}}^\otimes(\ESh(\textup{\'Et}_{/k},\ZZ_l),\widehat{\mathcal{D}}^\otimes(\ZZ_l))$, where the target is the based loop
space with respect to $F$ that is a group
object in $\mathcal{S}$.
Recall that $\mathsf{R}_{\textup{\'et},\ZZ_l}:\mathsf{DM}_{\textup{gm}}(\ZZ)\to \AMod_{\ZZ_l}$ factors through $F:\widehat{\mathcal{D}}_{\textup{gm}}^\otimes(\ESh(\textup{\'Et}_{/k},\ZZ_l))\to\lim_n\AMod^\otimes_{\ZZ/l^n}\subset \widehat{\mathcal{D}}^\otimes(\ZZ_l)$. Thus if $\Aut(\mathsf{R}_{\textup{\'et},\ZZ_l})$
is the based loop space of $\Map_{\CAlg(\Cat)}(\mathsf{DM}^\otimes_{\textup{gm}}(\ZZ),\AMod^\otimes_{\ZZ_l})$ with respect to $\mathsf{R}_{\textup{\'et},\ZZ_l}$, then the
vertical compositions with
$\Gal(\bar{k}/k)\to \Aut(F)$ induces a map of group objects $\Gal(\bar{k}/k)\to \Aut(\mathsf{R}_{\textup{\'et},\ZZ_l})$.
If we define $\Aut(\mathsf{R}_{\textup{\'et},\QQ_l})$
in a similar way, then we also have $\Gal(\bar{k}/k)\to \Aut(\mathsf{R}_{\textup{\'et},\QQ_l})$. We will refer to $\Gal(\bar{k}/k)\to \Aut(\mathsf{R}_{\textup{\'et},\ZZ_l})$ and $\Gal(\bar{k}/k)\to \Aut(\mathsf{R}_{\textup{\'et},\QQ_l})$
as the Galois action (or the action of $\Gamma$) on $\mathsf{R}_{\textup{\'et},\ZZ_l}$ and $\mathsf{R}_{\textup{\'et},\QQ_l}$ respectively.
We have constructed the $l$-adic \'etale realization functor
$\mathsf{R}_{\textup{\'et},\ZZ_l}$ which is endowed with
the Galois action $\Gamma\to \Aut(\mathsf{R}_{\textup{\'et},\ZZ_l})$.
Furthermore, there is its rational version $\mathsf{R}_{\textup{\'et},\QQ_l}$.
The following lemmata complete the proof of Proposition~\ref{realization}.

\begin{Lemma}
\label{rationaldescent}
The composition with $\mathsf{DM}^\otimes_{\textup{gm}}(\ZZ)\to \mathsf{DM}^\otimes_{\textup{gm}}(\QQ)$ induces
a categorical equivalence
\[
\Map^\otimes_{\textup{ex}}(\mathsf{DM}^\otimes_{\textup{gm}}(\QQ),\Mod_{\QQ_l}^\otimes)\to 
\Map^\otimes_{\textup{ex}}(\mathsf{DM}^\otimes_{\textup{gm}}(\ZZ),\Mod_{\QQ_l}^\otimes),
\]
where $\Map^\otimes_{\textup{ex}}(-,-)$ denotes the full subcategory of $\Map_{\CAlg(\Cat)}(-,-)$,
spanned by those functors which preserve finite colimits, i.e., exact functors.
\end{Lemma}

\Proof
The objects in $\mathsf{DM}^\otimes_{\textup{gm}}(\ZZ)$
forms a set of compact generators, and thus
we have an equivalence $\Map^\otimes_{\textup{ex}}(\mathsf{DM}^\otimes_{\textup{gm}}(\ZZ),\Mod_{\QQ_l}^\otimes)\simeq \Map_{\textup{L}}^\otimes(\mathsf{DM}^\otimes(\ZZ),\Mod_{\QQ_l}^\otimes)$ which is given by composition
with the inclusion
$\mathsf{DM}^\otimes_{\textup{gm}}(\ZZ)\subset \mathsf{DM}^\otimes(\ZZ)$.
Here $\Map_{\textup{L}}^\otimes(-,-)$ is the full subcategory of $\Map_{\CAlg(\Cat)}(-,-)$ spanned by those functors which preserve small colimits.

Let $\mathsf{DM}(\ZZ)[\ZZ^{-1}]$ be the stable presentable
$\infty$-category obtained from $\NSp(\ZZ)$, endowed with
the model structure of left Bousfield localization
with respect to $S=\{m:F_a(\ZZ(X)[n])\to F_a(\ZZ(X)[n]);X\in \textup{Sm}_{/k},n\in \ZZ, a\ge0, m\in \NN\}$, by inverting $S$-equivalences. Here $m$ means the multiplication by $m$.
See \cite[Definition 6.3]{HoG} for the notation $F_a$.
The class of
$S$-equivalences is closed under tensoring with cofibrant objects.
Indeed, to see this, it will suffice to show that
for any cofibrant object $C\in \NSp(\ZZ)$ and the cone $T$
of $m:\ZZ(X)[n]\to\ZZ(X)[n]$, $C\otimes F_a(T)$ is $S$-equivalent to zero.
We may assume that $C$ is a relative $I$-cell complex in the sense of \cite{Ho},
where $I$ is the set of generaring cofibrations. Thus it is enough to see that
$F_a(T)\otimes D$ is $S$-equivalent to zero where $D$ is either
a domain or target of generating cofibrations; it follows from 
a direct calculation.
Therefore, according to Lemma~\ref{loccomp}, $\mathsf{DM}^\otimes(\ZZ)[\ZZ^{-1}]$
is equivalent to the symmetric monoidal $\infty$-category
obtained from $\mathsf{DM}(\ZZ)$ as the localization with respect to
$S$; see \cite[4.1.3.4]{HA}.
By the universality of localization (cf. \cite[5.5.4.20]{HTT}, \cite[4.1.3.4]{HA}), there is a natural equivalence
\[
\Map_{\textup{L}}^\otimes(\mathsf{DM}^\otimes(\ZZ)[\ZZ^{-1}],\Mod_{\QQ_l}^\otimes)\simeq 
\Map_{\textup{L}}^\otimes(\mathsf{DM}^\otimes(\ZZ),\Mod_{\QQ_l}^\otimes).
\]
Let $\mathsf{DM}_{\textup{gm}}(\ZZ)[\ZZ^{-1}]$
be the smallest stable idempotent complete subcategory of
$\mathsf{DM}(\ZZ)[\ZZ^{-1}]$
which consists of the image of
$\{\ZZ(X)\otimes \ZZ(n)\}_{X\in \textup{Sm}_{/k}, n\in \ZZ}$, that forms
a set of compact generators of $\mathsf{DM}(\ZZ)[\ZZ^{-1}]$.
Then $\Map_{\textup{ex}}^\otimes(\mathsf{DM}_{\textup{gm}}^\otimes(\ZZ)[\ZZ^{-1}],\Mod_{\QQ_l}^\otimes)\simeq \Map_{\textup{L}}^\otimes(\mathsf{DM}^\otimes(\ZZ)[\ZZ^{-1}],\Mod_{\QQ_l}^\otimes)$. Thus it will suffice to prove that
the natural symmetric monoidal functor
$\mathsf{DM}_{\textup{gm}}^\otimes(\ZZ)[\ZZ^{-1}]\to \mathsf{DM}_{\textup{gm}}^\otimes(\QQ)$
is an equivalence. 
It is enough to prove that
a categorical equivalence
$\mathsf{DM}_{\textup{gm}}(\ZZ)[\ZZ^{-1}]\simeq \mathsf{DM}_{\textup{gm}}(\QQ)$.
To this end, we let $\mathsf{DM}^{\textup{eff}}(\ZZ)[\ZZ^{-1}]$
be the localization of $\mathsf{DM}^{\textup{eff}}(\ZZ)$ (Nisnevich version
of $\mathsf{DM}_{\textup{\'et}}^{\textup{eff}}(\ZZ)$) with respect to
$T=\{m:\ZZ(X)[n]\to \ZZ(X)[n];X\in \textup{Sm}_{/k},n\in \ZZ,m\in \NN\}$.
Using the argument above or \cite[Corollary 4.11]{CD1} we see that
the collection of $T$-equivalences is closed under tensoring
with cofibrant objects. Thus by \cite[4.1.3.4]{HA}
$\mathsf{DM}^{\textup{eff}}(\ZZ)[\ZZ^{-1}]$
is equipped with a symmetric monoidal structure.
Observe that there is an equivalence of symmetric monoidal
$\infty$-categories $(\mathsf{DM}^{\textup{eff}})^\otimes(\ZZ)[\ZZ^{-1}]\simeq (\mathsf{DM}^{\textup{eff}})^\otimes(\QQ)$.
For this, consider the adjoint pair $R:\Comp(\Sh(\Cor,\ZZ))\rightleftarrows \Comp(\Sh(\Cor,\QQ)):U$, where $U$ is the forgetful functor, and $R$ is induced by
the rationalization functor $\ZZ\textup{-mod}\to \QQ\textup{-mod}$ that
is the left adjoint of the forgetful functor
$\QQ\textup{-mod}\to \ZZ\textup{-mod}$.
If one equips $\Comp(\Sh(\Cor,\ZZ))$ and $\Comp(\Sh(\Cor,\QQ))$
with the model structure in which weak equivalences (resp. cofibrations)
are quasi-isomorphisms (resp. termwise
monomorphisms) (cf. \cite{Bek}, \cite[Theorem 2.1]{CD1}), then $(R,U)$ is a Quillen adjunction
since the rationalization functor is exact.
Moreover, $U$ induces a fully faithful right derived functor:
$\mathcal{D}(\Sh(\Cor,\QQ))\to \mathcal{D}(\Sh(\Cor,\ZZ))$.
Note that the adjunction $\mathcal{D}(\Sh(\Cor,\ZZ))\rightleftarrows \mathcal{D}(\Sh(\Cor,\QQ))$
is the localization with respect to $T$.
Unwinding the definition $C\in \mathcal{D}(\Sh(\Cor,\ZZ))$
is a $T$-local object if and only if $\Ext_{\Sh(\Cor,\ZZ)}^n(\ZZ(X),C)$ is a $\QQ$-vector space
for any $n\in \ZZ$ and any $X\in \textup{Sm}_{/k}$.
Clearly, the essential image of $\mathcal{D}(\Sh(\Cor,\QQ))\to \mathcal{D}(\Sh(\Cor,\ZZ))$ lies in the full subcategory of $T$-local objects.
Conversely, let $D$ be a cofibrant-fibrant
$T$-local object in $\Comp(\Sh(\Cor,\ZZ))$. The rationalization functor
$\ZZ\textup{-mod}\to \QQ\textup{-mod}$ is exact, and the presheaf
$X\mapsto H^n(D(X))$
is the same as $X\mapsto H^n(R(D(X)))$.
If one denotes by $R(D)'$ the fibrant replacement of $R(D)$,
then the natural map $D\to U(R(D)')$ is a quasi-isomorphism.
It follows that $D$ lies in the essential image of
$\mathcal{D}(\Sh(\Cor,\QQ))$. Hence we see that
$(R,U)$ is the localization with respect to $T$.
Consequently, $\mathsf{DM}^{\textup{eff}}(\QQ)$
is the localization of $\mathcal{D}(\Sh(\Cor,\ZZ))$ with respect to $T\cup \{\ZZ(X\times \mathbb{A}^1)\to \ZZ(X)\}_{X\in \textup{Sm}_{/k}}$.
We then have $(\mathsf{DM}^{\textup{eff}})^\otimes(\ZZ)[\ZZ^{-1}]\simeq (\mathsf{DM}^{\textup{eff}})^\otimes(\QQ)$.
Next observe that $\mathsf{DM}^{\textup{eff}}(\ZZ)[\ZZ^{-1}]\to \mathsf{DM}(\ZZ)[\ZZ^{-1}]$ induced by $\Sigma^\infty:\mathsf{DM}^{\textup{eff}}(\ZZ)\to \mathsf{DM}(\ZZ)$ is fully faithful.
To see this, it is enough to show that $\Sigma^\infty$
sends $T$-local objects to
$S$-local objects. We here remark that by Voevodsky's cancellation theorem
$\Sigma^{\infty}$ is fully faithful.
Let $C$ be a $T$-local object. To check that $\Sigma^\infty(C)$
is $S$-local, it will suffice to prove that $C\otimes \ZZ(1)$ is
$T$-local. Let $\mathcal{C}\subset \mathcal{D}(\Sh(\Cor,\QQ))$ be the stable subcategory that consists of those objects
$C$ such that $C\otimes \ZZ(1)$ lies in $\mathcal{D}(\Sh(\Cor,\QQ))$, that is,
$T$-local. Then $\QQ(X)\otimes \ZZ(\mathbb{G}_m)\simeq \QQ(X\times \mathbb{G}_m)$ in $\Comp(\Sh(\Cor,\ZZ)$, and the
Suslin complex $C_*(\QQ(X\times \mathbb{G}_m))$ belongs to
$\Comp(\Sh(\Cor,\QQ))$ (see \cite[2.14]{MVW} for Suslin complexes).
Thus we deduce that $\QQ(X)\in \mathcal{C}$ for any $X\in \textup{Sm}_{/k}$.
Moreover, $\mathcal{C}$ has small coproducts such that
$\mathcal{C}\hookrightarrow \mathcal{D}(\Sh(\Cor,\QQ))$
preserves small coproducts.
Hence $\mathcal{C}=\mathcal{D}(\Sh(\Cor,\QQ))$, and $\Sigma^\infty$
sends $T$-local objects to $S$-local objects.
On the other hand,
the composite
$\mathsf{DM}^{\textup{eff}}(\ZZ)[\ZZ^{-1}]\to \mathsf{DM}(\ZZ)[\ZZ^{-1}]
\to \mathsf{DM}(\QQ)$ is fully faithful since it can be identified with
$\mathsf{DM}^{\textup{eff}}(\ZZ)[\ZZ^{-1}]\simeq \mathsf{DM}^{\textup{eff}}(\QQ)\stackrel{\Sigma^\infty}{\hookrightarrow} \mathsf{DM}(\QQ)$.
Since $\mathsf{DM}^{\textup{eff}}(\ZZ)[\ZZ^{-1}]\to \mathsf{DM}(\ZZ)[\ZZ^{-1}]$
is fully faithful, $\mathsf{DM}(\ZZ)[\ZZ^{-1}]
\to \mathsf{DM}(\QQ)$ is fully faithful when one
restricts the domain to the essential image of $\mathsf{DM}^{\textup{eff}}(\ZZ)[\ZZ^{-1}]$. Let $\mathsf{DM}^{\textup{eff}}(\ZZ)[\ZZ^{-1}]\otimes \ZZ(n)$
be the full subcategory of $\mathsf{DM}(\ZZ)[\ZZ^{-1}]$
spanned by $C\otimes \ZZ(n)$ such that 
$C$ lies in the essential image of $\mathsf{DM}^{\textup{eff}}(\ZZ)[\ZZ^{-1}]$.
We define $\mathsf{DM}^{\textup{eff}}(\QQ)\otimes \QQ(n)$ in a similar way.
Then we have $\cup_{n\ge0}\mathsf{DM}^{\textup{eff}}(\ZZ)[\ZZ^{-1}]\otimes \ZZ(-n)\simeq \cup_{n\ge0}\mathsf{DM}^{\textup{eff}}(\QQ)\otimes \QQ(-n)$.
Since $\mathsf{DM}_{\textup{gm}}(\ZZ)[\ZZ^{-1}]\subset \cup_{n\ge0}\mathsf{DM}^{\textup{eff}}(\ZZ)[\ZZ^{-1}]\otimes \ZZ(-n)$, we have
$\mathsf{DM}_{\textup{gm}}(\ZZ)[\ZZ^{-1}]\simeq \mathsf{DM}_{\textup{gm}}(\QQ)$.
\QED

\begin{Lemma}
\label{etalecoh}
Let $X$ be a smooth projective scheme over $k$ and let $\ZZ(X)$ be
the object in $\mathsf{DM}(\ZZ)$ corresponding to $X$.
Let $\ZZ(X)^\vee$ be its dual.
Put $R(X,m)=\mathsf{R}_{\textup{\'et},\ZZ_l}(\ZZ(X)^{\vee}\otimes\ZZ(m))$.
Let $\overline{X}=X\times_{k}\bar{k}$.
Then $H^s(R(X,m))$ is naturally isomorphic to the \'etale cohomology
$H^s_{\textup{\'et}}(\overline{X}, \ZZ_l(m))$,
and the Galois action on $H^s(R(X,m))$ induced by that on $\mathsf{R}_{\textup{\'et},\ZZ_l}$
coincides with that of $H^s_{\textup{\'et}}(\overline{X}, \ZZ_l(m))$.
\end{Lemma}

\Proof
Let $\mathcal{H}om(\ZZ/l^n(X),\ZZ/l^n)$
be the internal Hom object in $\mathsf{DM}_{\textup{\'et}}(\ZZ/l^n)$. 
We remark that $\mathcal{H}om(\ZZ/l^n(X),-)$ is defined to be the right
adjoint of the tensor operation
$(-)\otimes \ZZ/l^n(X):\mathsf{DM}_{\textup{\'et}}(\ZZ/l^n)\to \mathsf{DM}_{\textup{\'et}}(\ZZ/l^n)$. Since $\ZZ(X)$ is dualizable and
$\mathsf{DM}(\ZZ)\to \mathsf{DM}_{\textup{\'et}}(\ZZ/l^n)$ is symmetric monoidal,
$\ZZ(X)^\vee$ maps to $\mathcal{H}om(\ZZ/l^n(X),\ZZ/l^n)$. Now through the equivalence $\mathsf{DM}_{\textup{\'et}}(\ZZ/l^n)\simeq \mathsf{DM}_{\textup{\'et}}^{\textup{eff}}(\ZZ/l^n)$, we regard $\mathcal{H}om(\ZZ/l^n(X),\ZZ/l^n)$ as an object in $\mathsf{DM}_{\textup{\'et}}^{\textup{eff}}(\ZZ/l^n)$. Take its fibrant model $\mathbf{H}\textup{om}(\ZZ/l^n(X),\ZZ/l^n)$ in $\Comp(\ESh(\Cor,\ZZ/l^n))$.
Regard $res(\mathbf{H}\textup{om}(\ZZ/l^n(X),\ZZ/l^n))$
in $\Comp(\ESh(\textup{\'Et}_{/k},\ZZ/l^n))$
as a complex of discrete $\Gamma$-modules as follows:
Put $\bar{k}=\colim_i k_i$ where the right hand side is a filtered 
colimit of finite Galois extensions of $k$.
Then the filtered colimit
$\colim_ires(\mathbf{H}\textup{om}(\ZZ/l^n(X),\ZZ/l^n))(k_i)$
is a discrete $\ZZ/l^n$-modules with action of $\Gamma$ which is determined
by the natural
actions of $\Gal(k_i/k)$ on $res(\mathbf{H}\textup{om}(\ZZ/l^n(X),\ZZ/l^n))(k_i)$.
It
represents the image of $\ZZ(X)^\vee$ in $\mathcal{D}(\ZZ/l^n)$.
Note that
$C(k_i):=res(\mathbf{H}\textup{om}(\ZZ/l^n(X),\ZZ/l^n))(k_i)$
can be identified with
\begin{eqnarray*}
\Ext^s(\ZZ/l^n(\Spec k_i),\mathbf{H}\textup{om}(\ZZ/l^n(X),\ZZ/l^n)) &\simeq& \Ext^s(\ZZ/l^n(X\times_kk_i),\ZZ/l^n) \\
&\simeq& H^s_{\textup{\'et}}(X\times_kk_i,\ZZ/l^n)
\end{eqnarray*}
where $\Ext^s(-,-)$ is $\pi_0(\Map_{\mathsf{DM}_{\textup{\'et}}(\ZZ/l^n)}(-,-[s]))$. The second isomorphism follows from the equivalence in Lemma~\ref{suslin},
see also \cite[Theorem 10.2]{MVW}.
Through isomorphisms,
the action of $\Gal(k_i/k)$ on
$C(k_i)$
coincides with the action on
$H^s_{\textup{\'et}}(X\times_kk_i,\ZZ/l^n)$.
(Here $\Gal(k_i/k)$ acts on the $k$-scheme $X\times_kk_i$ in the obvoius way,
and it gives rise to action on $H^n_{\textup{\'et}}(X\times_kk_i,\ZZ/l^n)$.)
Since a filtered colimit commutes with
taking cohomology groups, we have isomorphisms of $\Gamma$-modules
\[
H^s(\colim_iC(k_i))\simeq \colim_iH^s(C(k_i))\simeq \colim_i H^s_{\textup{\'et}}(X\times_kk_i,\ZZ/l^n)\simeq H_{\textup{\'et}}^s(\overline{X},\ZZ/l^n).
\]
Taking account of the $t$-exactness of the equivalence
$\AMod_{\ZZ_l}\simeq \lim_n\AMod_{\ZZ/l^n}$ (see \cite[5.3.1, 5.2.12]{DAGn}),
we conclude that $H^s(R(X,0))\simeq \lim_n H_{\textup{\'et}}^s(\overline{X},\ZZ/l^n)$. (More explicitly, through the equivalence,
$(M_n)_{n\ge1}\in
\lim_n\AMod_{\ZZ/l^n}$ with $M_n\in \AMod_{\ZZ/l^n}$ corresponds
to the filtered limit $\lim_n U_n(M_n)$ in $\Mod_{\ZZ_l}$ where $U_n:\Mod_{\ZZ/l^n}\to \Mod_{\ZZ_l}$ is naturally induced by $\ZZ_l\to \ZZ/l^n$,
and thus $H^s(R(X,0))\simeq \lim_n H_{\textup{\'et}}^s(\overline{X},\ZZ/l^n)$
follows
from Milnor exact sequence, Mittag-Leffler condition and the finiteness
of \'etale cohomology.) We have
$H^s(R(X,0))\simeq H^s_{\textup{\'et}}(\overline{X}, \ZZ_l)$.

Similarly, the image of $\ZZ(m)$ in
$\mathcal{D}(\ESh(\textup{\'Et}_{/k},\ZZ/l^n))$
is $\ZZ/l^n(m)$, that is, the object
corresponding to $\ZZ/l^n(m)$ in
$\mathsf{DM}_{\textup{\'et}}(\ZZ/l^n)$.
By \cite[10.6, 10.2]{MVW}
there is the natural
equivalence $\ZZ/l^n(m)\simeq \mu_{l^n}^{\otimes m}$
in $\mathcal{D}(\ESh(\textup{\'Et}_{/k},\ZZ/l^n))$,
which corresponds to $\ZZ/l^n$-module
$\mu_{l^n}^{\otimes m}(\bar{k})$ placed in degree zero, which is endowed with the natural action
$\Gamma=\Gal(\bar{k}/k)$. Here $\mu_i$ is the sheaf given by $L\mapsto
\{a\in L|\ a^i=1\}$.
Thus we see that
$\mathsf{R}_{\textup{\'et},\ZZ_l}(\ZZ(m))$
is equivalent to $\lim_n\mu_{l^n}^{\otimes m}(\bar{k})$ placed in degree zero.
Finally, the isomorphism $H^s_{\textup{\'et}}(\overline{X}, \ZZ_l(m))\simeq H^s_{\textup{\'et}}(\overline{X}, \ZZ_l)\otimes\ZZ_l(m)$ implies also
the case of $m\neq 0$.
\QED

\begin{Corollary}
\label{suslinho}
Let $X$ be a smooth scheme over $k$, and let
$\mathsf{R}_{\textup{\'et},\ZZ_l}(\ZZ(X))(n)$ denote
the image of $\ZZ(X)$ in $\Mod_{\ZZ/l^n}$.
Then there is an isomorphism of $\ZZ/l^n$-modules
\[
\Hom_{\ZZ/l^n}(H_s(\mathsf{R}_{\textup{\'et},\ZZ_l}(\ZZ(X))(n)),\ZZ/l^n)\simeq H^s_{\textup{\'et}}(\overline{X},\ZZ/l^n).
\]
\end{Corollary}

\Proof
The same argument as in Lemma~\ref{etalecoh} shows
that
$H_s(\mathsf{R}_{\textup{\'et},\ZZ_l}(\ZZ(X))(n))$
is the algebraic singular homology $H_s^{sing}(\overline{X},\ZZ/l^n)$ \cite[10.8]{MVW} of $\overline{X}$.
Hence our assertion follows from \cite[10.11]{MVW}.
\QED

\subsection{}
\label{verytech}
We conclude this Section with technical results; Lemma~\ref{di}, \ref{loccomp}.
Let $G$ be a reductive algebraic group over a field $\mathbf{K}$ of characteristic zero.
Let $\mathcal{M}=\textup{Vect}_{\mathbf{K}}(G)$
be the category of $\mathbf{K}$-vector spaces.
It is a Grothendieck semisimple abelian category.
Let $\mathcal{G}_{\mathcal{M}}$
be the set of finite coproducts of irreducible representations of $G$.
Let $\mathcal{H}_{\mathcal{M}}=\{0\}$.
Then we easily see that the pair $(\mathcal{G}_{\mathcal{M}},\mathcal{H}_{\mathcal{M}})$ is a flat descent structure in the sense of \cite{CD1}.
We equip the category $\Comp(\mathcal{M})$ of chain complexes of objects in $\mathcal{M}$
with the $\mathcal{G}_{\mathcal{M}}$-model structure; see {\it loc. cit.}.
Let $\mathcal{D}(\mathcal{M})^\otimes$ be the
symmetric monoidal $\infty$-category
obtained from the full subcategory of cofibrant objects
of $\Comp(\mathcal{M})$ by inverting weak equivalences.
Let $\mathcal{D}_\vee(\mathcal{M})^\otimes$
denote the full subcategory of dualizable objects
in $\mathcal{D}(\mathcal{M})^\otimes$.
In \cite[Section A.6]{Tan} we define a symmetric monoidal
stable presentable $\infty$-category $\Rep^\otimes_G$.
Intuitively speaking, $\Rep^\otimes_G$
is a symmetric monoidal $\infty$-category which consists of complexes
of $\mathbf{K}$-vector spaces endowed with action of $G$.
We here recall the definition of $\Rep^\otimes_G$
by using model categories.
Put $\Spec B=G$. The group structure of $G$
gives rise to a cosimplicial diagram $\{B^{\otimes n}\}_{[n]\in \Delta}$
of commutative $\mathbf{K}$-algebras whose $n$-th term
is $B^{\otimes n}$, i.e., it comes from the \v{C}ech nerve
of the natural projection to
the classifying stack
$\pi:\Spec \mathbf{K}\to BG$; see \cite[6.1.2]{HTT}
for \v{C}ech nerves.
For a commutative algebra $A$, we let $\textup{Comp}(A)$
be the category of (not necessarily bounded)
complexes of $A$-modules.
We here equip $\textup{Comp}(A)$ with
the projective model structure (cf. \cite[2.3.3]{Ho}).
The cosimplicial diagram $\{B^{\otimes n}\}_{[n]\in \Delta}$ yields a cosimplicial diagram of (symmetric monoidal) categories
$\{\textup{Comp}(B^{\otimes n})\}_{[n]\in \Delta}$
in which each $\textup{Comp}(B^{\otimes n})\to \textup{Comp}(B^{\otimes m})$
is the base change by $B^{\otimes n}\to B^{\otimes m}$, that is a left Quillen adjoint functor.
Then it gives rise to
a cosimplicial diagram $\{\textup{Comp}(B^{\otimes n})^c\}_{[n]\in \Delta}$
of symmetric monoidal categories consisting of cofibrant objects.
Inverting quasi-isomorphisms in each category
$\textup{Comp}(B^{\otimes n})^c$ we obtain a cosimplicial
diagram $\{\NNNN_W(\textup{Comp}(B^{\otimes n})^c)\}_{[n]\in \Delta}$
of symmetric monoidal stable presentable $\infty$-categories.
Here $\NNNN_W(-)$ indicates the symmetric monoidal stable presentable $\infty$-category obtained by inverting weak equivalences.
The superscript $(-)^c$ indicates the full subcategory of cofibrant objects.
(See \cite[4.1.3]{HA})
We define $\Rep^\otimes_G$ to be a limit of $\{\NNNN_W(\textup{Comp}(B^{\otimes n})^c)\}_{[n]\in \Delta}$
among symmetric monoidal $\infty$-categories.
The limit $\Rep^\otimes_G$ is also stable and presentable.
Let $\PRep_G^\otimes\subset \Rep^\otimes_G$
be the full subcategory of dualizable objects.
The following Lemma gives a relation between
$\Rep^\otimes_G$ and $\mathcal{D}(\mathcal{M})^\otimes$.

\begin{Lemma}
\label{di}
There exists an equivalence
$\mathcal{D}(\mathcal{M})^\otimes\simeq \Rep^\otimes_G$
of symmetric monoidal $\infty$-categories.
\end{Lemma}

\Proof
We first construct a symmetric monoidal functor
$\mathcal{D}(\mathcal{M})^\otimes \to \Rep^\otimes_G$
which preserves small colimits.
Define
$\Comp(\textup{Vect}_{\mathbf{K}}(G))\to \textup{Comp}(B^{\otimes n})$
to be the functor induced by
$\Spec B^{\otimes n}\to \Spec \mathbf{K} \stackrel{\pi}{\to} BG$
for any $n\ge0$.
Consider $\Comp(\textup{Vect}_{\mathbf{K}}(G))=\Comp(\mathcal{M})$ to
be
the constant cosimplicial diagram.
These induce a map of cosimplicial diagrams
$\Comp(\mathcal{M}) \to \{\textup{Comp}(B^{\otimes n})\}_{[n]\in \Delta}$.
Note that each functor
$\Comp(\mathcal{M})\to \textup{Comp}(B^{\otimes n})$
preserves cofibrant objects since it preserves small colimits
and 
the generating cofibration
$\{S^{n+1}E\to D^{n}E\}_{n\in\ZZ,E\in \mathcal{G}_{\mathcal{M}}}$
maps to cofibrations in $\textup{Comp}(B^{\otimes n})$.
We then have the map
$\Comp(\mathcal{M})^c \to \{\textup{Comp}(B^{\otimes n})^c\}_{[n]\in \Delta}$.
By inverting weak equivalences, it gives rise to a map of cosimplicial
symmetric monoidal $\infty$-categories
\[
\NNNN_W(\Comp(\mathcal{M})^c)\to \{\NNNN_W(\textup{Comp}(B^{\otimes n})^c)\}_{[n]\in \Delta}.
\]
Since $\Rep^\otimes_G$ is the limit of $\{\NNNN_W(\textup{Comp}(B^{\otimes n})^c)\}_{[n]\in \Delta}$, we obtain a symmetric monoidal colimit-preserving
functor
$\mathcal{D}(\mathcal{M})^\otimes\to \Rep^\otimes_G$.

Next we define a $t$-structure
on $\Rep_G$.
Let $\Rep_{G,\ge0}$ (resp. $\Rep_{G,\le0}$)
be the inverse image of $\Mod_{\mathbf{K},\ge0}$ (resp. $\Mod_{\mathbf{K},\le0}$)
under the forgetful functor
$p:\Rep_G\to \Mod_{\mathbf{K}}$.
Here $C\in \Mod_{\mathbf{K}}$ belongs to $\Mod_{\mathbf{K},\ge0}$
(resp. $\Mod_{\mathbf{K},\le0}$)
if and only if $\pi_i(C)=0$ for any $i<0$ (resp. $i>0$).
The comonad $T:\Mod_{\mathbf{K}}\to \Mod_{\mathbf{K}}$
of the adjoint pair
\[
p:\Rep_G\rightleftarrows \Mod_{\mathbf{K}}:q
\]
is given by $C\mapsto B\otimes C$. Here $q$ is a right adjoint of $p$.
Identifying $\Rep_G$ with
the $\infty$-category of $T$-comodules by \cite[6.2.4.1]{HA},
we conclude by \cite[VII, 6.20]{DAGn} that
$(\Rep_{G,\ge0},\Rep_{G,\le0})$
defines a both left and right complete $t$-structure.
Let $\Rep^b_G$ (resp. $\Rep^+_G$)
denote the full subcategory of $\Rep_G$
spanned by bounded objects (resp. left bounded
objects) with respect to this $t$-structure.
Since $p$ is symmetric monoidal, they are stable under tensor product.

We claim that $\mathcal{D}(\mathcal{M})\to \Rep_G$ induces a categorical equivalence
$\mathcal{D}^+(\mathcal{M})\to \Rep^+_G$.
We first
prove that the induced functor
$w:\mathcal{D}^+(\mathcal{M})\to \Rep^+_G$
is fully faithful.
Let $C$ and $C'$ be objects in $\mathcal{D}^+(\mathcal{M})$.
We need to show that the induced map
$w_{C,C'}:\Map_{\mathcal{D}(\mathcal{M})}(C,C')\to \Map_{\Rep_G}(w(C),w(C'))$ is an equivalence in $\SSS$.
Since $\mathcal{D}(\mathcal{M})\to \Rep_G$ preserves small colimits,
the $t$-structure on $\mathcal{D}(\mathcal{M})$
is right complete and $C$ is a colimit of bounded objects, thus
we may assume that $C$ lies in $\mathcal{D}^b(\mathcal{M})$.
The full subcategory of $\mathcal{D}^b(\mathcal{M})$ spanned by
those objects $C$ such that $w_{C,C'}$ is an equivalence
for any $C'\in \mathcal{D}^+(\mathcal{M})$,
is a stable subcategory.
Hence we may and will assume that $C$ belongs to the heart $\mathcal{M}$.
To compute $\Ext^n_{\mathcal{D}(\mathcal{M})}(C,C')$,
we use the injective model structure on $\Comp(\mathcal{M})$
in which weak equivalences are quasi-isomorphisms, and cofibrations
are monomorphisms (cf. \cite{Bek}, \cite{CD1}, \cite[1.3.5]{HA}).
Since $\mathcal{M}$ has enough injective objects we suppose that $C'$ is a left bounded complex of the form
\[
\cdots\to 0\to 0 \to C^r\to C^{r+1}\to  \cdots
\]
where $C^i$ is an injective object $\mathcal{M}$
for any $i\in \ZZ$. It is a fibrant object.
To compute $\Map_{\mathcal{D}(\mathcal{M})}(C,C')$ and $\Map_{\Rep_G}(w(C),w(C'))$, 
since $C$ lies in the heart and $w$ is $t$-exact,
we may and will suppose that $C^i=0$ for $i>1$.
Let $\mathcal{I}$ be the full subcategory spanned by finite-length
complexes of injective objects.
We claim that for any $C'\in \mathcal{I}$ the map $\theta_{C,C'}^n:\Ext^n_{\mathcal{D}(\mathcal{M})}(C,C')\to \Ext^n_{\Rep_G}(w(C),w(C'))$ is an isomorphism for any $n\in \ZZ$.
If $\mathcal{P}$ is the full subcategory of $\mathcal{D}^b(\mathcal{M})$ spanned by
finite-length complexes $C'$ of injective objects
such that $\theta_{C,C'}^n$ is an isomorphism for any $n\in \ZZ$,
then $\mathcal{P}$ is stable under shifts and cones.
Thus we may and will assume that $C'$ is an injective object
in the heart $\mathcal{M}$.
When $n\le 0$, clearly it is an isomorphism.
When $n>0$, we will prove that $\Ext^n_{\mathcal{D}(\mathcal{M})}(C,C')=\Ext^n_{\Rep_G}(w(C),w(C'))=0$.
To see $\Ext^n_{\Rep_G}(w(C),w(C'))=0$ for $n>0$,
let $I=w(C')$ be an injective object
in the heart $\textup{Vect}_{\mathbf{K}}(G)$ of $\Rep_G$ and let
$p(I)\to J$ be an injective resolution,
 that is, $J$ is an injective object in the heart $\textup{Vect}_{\mathbf{K}}$ of $\Mod_{\mathbf{K}}$.
Then $q(J)$ is injective, and
$I\to q(p(I))\to q(J)$ is a monomorphism
since $p(I)\to p(q(J))\to J$ is the monomorphism where
$p(q(J))\to J$ and $I\to q(p(I))$ are a counit map and 
a unit map respectively.
Notice that $I$ is injective, thus $I$ is a retract of
$q(J)$.
Consequently, it will suffice to show that
$\Ext^n_{\Rep_G}(w(C),q(J))=0$ for $n>0$.
It follows from the adjunction that
$\Ext_{\Rep_G}^n(w(C),q(J))=\Ext^n_{\Mod_{\mathbf{K}}}(p(w(C)),J)=0$ for $n>0$.
Since $C$ is cofibrant and $C'[r]$ is fibrant
in any $r\in \ZZ$ in $\textup{Comp}(\mathcal{M})$
endowed with the injective model structure, $\Ext_{\mathcal{D}(\mathcal{M})}^n(C,C')=0$ for any $n>0$.

Next we will prove that $w$ is essentially surjective.
Let $D\in \Rep^+_G$.
We must show that there is $C\in \mathcal{D}^+(\mathcal{M})$ such that $w(C)\simeq D$.
Since $\Rep_G$ is right complete and $\mathcal{D}(\mathcal{M})\to \Rep_G$
preserves small colimits, thus by the fully faithfulness proved above,
we may and will suppose that
$D$ belongs to $\Rep^b_G$.
Let $l$ be the amplitude of $D$. We proceed by induction on $l$.
The case of $l=1$ is obvious (in this case $D$ is a shift of an object
in the heart).
Using $t$-structure one can take a distinguished triangle
\[
D_1\to D\to D_2\to D_1[1]
\]
such that
the amplitude of $D_1$ is equal or less than $l-1$,
and the amplitude of $D_2$ is equal or less than $1$.
By the inductive assumption, we have $C_1$ and $C_2$
such that $w(C_1)\simeq D_1$ and $w(C_2)\simeq D_2$.
Moreover, the fully faithfulness implies that
there exists $C_2\to C_1[1]$
such that $w(C_2)\to w(C_1[1])$
represents the homotopy class of $D_2\to D_1[1]$.
Note that $D$ is a fibre of $D_2\to D_1[1]$.
Let $C$ be a fibre of $C_2\to C_1[1]$.
By the exactness of $w$, we conclude that $w(C)\simeq D$.

It remains to show how one can  derive an equivalence $\mathcal{D}^\otimes(\mathcal{M})\simeq \Rep^\otimes_G$ from $\mathcal{D}^+(\mathcal{M})\simeq \Rep^+_G$.
We have constructed the symmetric monoidal functor
$\mathcal{D}^\otimes(\mathcal{M})\to \Rep^\otimes_G$, and thus
it suffices to prove that the underlying functor
$\mathcal{D}(\mathcal{M})\to \Rep_G$ is a categorical equivalence.
The equivalence $\mathcal{D}^+(\mathcal{M})\simeq \Rep^+_G$
induces an equivalence
$\mathcal{D}_\vee(\mathcal{M})\simeq \PRep_G$,
where $\mathcal{D}_\vee(\mathcal{M})$ denotes the full subcategory spanned by
dualizable objects.
Note that by the assumption that
$G$ is a reductive algebraic
group over a field of characteristic zero, $\mathcal{D}(\mathcal{M})$ is compactly generated, and the set of (finite dimensional)
irreducible representations is
a set of compact generators.
Thus $\mathcal{D}(\mathcal{M})\simeq \Ind(\mathcal{D}_\vee(\mathcal{M}))$.
Moreover, if $\PMod_{G}$ denotes
the full subcategory of $\Rep_G$
spanned by dualizable objects, then
$\Ind(\PRep_G)\simeq \Rep_G$; see \cite[3.22]{BFN}.
Hence we obtain $\mathcal{D}(\mathcal{M})\simeq \Ind(\mathcal{D}_\vee(\mathcal{M}))\simeq \Ind(\PRep_G)\simeq \Rep_G$.
Finally, we remark another way to deduce $\mathcal{D}(\mathcal{M})\simeq \Rep_G$. Since $\Rep_G$ is left complete and $\mathcal{D}^+(\mathcal{M})\simeq \Rep^+_G$, the functor $\mathcal{D}(\mathcal{M})\to \Rep_G$
can be viewed as a left completion \cite[1.2.1.17]{HA}
of $\mathcal{D}(\mathcal{M})$.
By using the semi-simplicity of $\mathcal{M}$ we can easily check
that $\mathcal{D}(\mathcal{M})$ is left complete.
\QED

\begin{Remark}
\label{tech}
As shown in the proof, we have
\[
\mathcal{D}(\mathcal{M})\simeq \Ind(\mathcal{D}_\vee(\mathcal{M}))\simeq \Ind(\PRep_G)\simeq \Rep_G.
\]
\end{Remark}

Let $\mathbb{M}$ be a left proper combinatorial model category.
Let $S$ be a small set of morphisms
in $\mathbb{M}$.
Then we have a new model structure of $\mathbb{M}$;
a left Bousfield localization of $\mathbb{M}$ with respect to $S$
(see e.g. \cite{Bwk}, \cite[A. 3.7.3]{HTT}), where (new) weak equivalences are
called $S$-equivalences.
We then obtain an $\infty$-category
$\NNNN_W(\mathbb{M}[S^{-1}])$  by inverting $S$-equivalences.
On the other hand, we have the $\infty$-category $\NNNN_W(\mathbb{M})$ obtained from $\mathbb{M}$
by inverting weak equivalences.
By using the localization theory at the level of $\infty$-category
\cite[5.5.4]{HTT},
one can take the localization $\NNNN_{W}(\mathbb{M})\to \NNNN_{W}(\mathbb{M})[S^{-1}]$ (see \cite[5.5.4.15]{HTT}).
Then the universality of the localization
$\mathbb{L}:\NNNN_{W}(\mathbb{M})\to \NNNN_{W}(\mathbb{M})[S^{-1}]$
\cite[5.5.4.20]{HTT} induces a functor
$F:\NNNN_{W}(\mathbb{M})[S^{-1}]\to \NNNN_{W}(\mathbb{M}[S^{-1}])$.

\begin{Lemma}
\label{loccomp}
The functor $\NNNN_{W}(\mathbb{M})[S^{-1}]\to \NNNN_{W}(\mathbb{M}[S^{-1}])$ is a categorical equivalence.
\end{Lemma}

\Proof
We have the commutative diagram
\[
\xymatrix{
  &  \NNNN_{W}(\mathbb{M}) \ar[rd]^{\mathbb{L}'} \ar[ld]_{\mathbb{L}} & \\
\NNNN_{W}(\mathbb{M})[S^{-1}]\ar[rr]^F& & \NNNN_{W}(\mathbb{M}[S^{-1}])
}
\]
that consists of left adjoint functors of presentable $\infty$-categories.
Here $\mathbb{L}'$ is the ``localization functor'' that comes from the
left Quillen functor.
Note that the right adjoint functors of $\mathbb{L}$ and $\mathbb{L}'$ are fully faithful.
We denote them by $i$ and $i'$ respectively.
Moreover, the essential image of $i$ consists of $S$-local objects,
that is, those objects $Z$ such that
$\Map_{\NNNN_{W}(\mathbb{M})}(Y,Z)\to \Map_{\NNNN_{W}(\mathbb{M})}(X,Z)$ is a weak homotopy equivalence for any $X\to Y\in S$.
Similarly, the essential image of $i'$ consists of ``model theoretic $S$-local objects'',
that is, those objects $Z$ such that
$\Map_{\mathbb{M}}(Y,Z)\to \Map_{\mathbb{M}}(X,Z)$ is a weak homotopy equivalence for any $X\to Y\in S$. Here we slightly abuse notation
and $\Map_{\mathbb{M}}(-,-)$ denotes the mapping space in $\mathbb{M}$
given by machinery of simplicial and cosimplicial frames \cite[5.4]{Ho}
or hammock localization of Dwyer-Kan (we implicitly assume suitable
cofibrant or fibrant replacements).
Thus it will suffice to prove both mapping spaces coincide.
If $\mathbb{M}$ is a simplicial model category, then
our assertion follows from \cite[1.3.4.20]{HA}.
In the general case, our assertion follows
from the simplicial case and a theorem of Dugger \cite{Dug} which
says that every combinatorial model category is Quillen equivalent to
a left proper simplicial combinatorial model category.
\QED

\end{document}